\documentclass[onefignum,onetabnum]{siamart190516}

\usepackage{lipsum}

\usepackage{cite}
\usepackage{amsfonts}
\usepackage{graphicx}
\usepackage{epstopdf}
\usepackage{algorithmic}
\ifpdf
  \DeclareGraphicsExtensions{.eps,.pdf,.png,.jpg}
\else
  \DeclareGraphicsExtensions{.eps}
\fi

\usepackage{xfrac}
\usepackage{animate}
\usepackage{enumerate}
\usepackage[shortlabels]{enumitem}
\usepackage[a4paper,margin=25mm]{geometry}    
\usepackage{setspace} 
\usepackage{amsfonts,amssymb,amsmath}         
\usepackage{graphicx}                         
\usepackage{siunitx}                          
\usepackage{hyperref}						  
\usepackage{subcaption}						  
\usepackage{wrapfig}						  
\usepackage{lipsum} 
\usepackage{xcolor}
\definecolor{myblue}{rgb}{.1, .4, .8}
\definecolor{myred}{rgb}{.8, .1, .1}
\definecolor{myyellow}{RGB}{246,235, 24}
\definecolor{mylime}{RGB}{150,201, 61}
\definecolor{mypurple}{RGB}{186,83, 159}
\definecolor{myorange}{RGB}{246,134, 59}
\definecolor{matlab1}{rgb}{0, 0.4470, 0.7410}
\definecolor{matlab2}{rgb}{0.8500, 0.3250, 0.0980}
\definecolor{matlab3}{rgb}{0.9290, 0.6940, 0.1250}
\definecolor{matlab4}{rgb}{0.4940, 0.1840, 0.5560}
\usepackage[usestackEOL]{stackengine}
\usepackage[T1]{fontenc}
\usepackage{animate}
\usepackage{booktabs}
\usepackage{longtable}
\usepackage{adjustbox}
\usepackage{mdframed}
\usepackage{todonotes}
\usepackage{xargs}
\usepackage{xstring}

\newsiamremark{remark}{Remark}
\crefname{remark}{Remark}{Remark}
\crefname{Example}{Example}{Examples}
\newsiamremark{hypothesis}{Hypothesis}
\crefname{hypothesis}{Hypothesis}{Hypotheses}
\newsiamthm{claim}{Claim}
\newsiamremark{assumption}{Assumption}
\crefname{assumption}{Assumption}{Assumption}

\headers{Regularisation of codimension-2 discontinuity sets}{ N. Cheesman, K. Uldall Kristiansen, S. J. Hogan}

\title{Regularisation of isolated codimension-2 discontinuity sets\thanks{Submitted to the editors {20/05/2021}.
\funding{{This work was funded by the UK Engineering and Physical Sciences Research Council (EPSRC) as part of the first author's PhD Studentship.}}}}

\author{
Noah Cheesman\thanks{Department of Engineering Mathematics, University of Bristol, Bristol BS8 1TW, UK
(\email{noah.cheesman@bristol.ac.uk}, \email{s.j.hogan@bristol.ac.uk}).}
\and Kristian Uldall Kristiansen\thanks{Department of Applied Mathematics and Computer Science, Technical University of Denmark, 2800 Kgs. Lyngby, Denmark  
  (\email{krkri@dtu.dk}).}
 \and S.J. Hogan\footnotemark[2]}

\usepackage{amsopn}

\newcommand{\pdiff}[2]{\frac{\partial{#1}}{\partial{#2}}}					
\renewcommand{\vec}[1]{\boldsymbol{#1}}										

\newcommand{\Real}{\mathbb{R}}

\newcommand{\sign}{\mathrm{sign}}

\newcommand{\tr}{\mathrm{tr}}

\renewcommand{\vec}[1]{{#1}}
\renewcommand{\vec}[1]{\mathbf{#1}}
\newcommand{\mat}[1]{\mathbf{#1}}
\renewcommand{\mat}[1]{\underline{\vec{#1}}}
\renewcommand{\mat}[1]{{\underline{\underline{\mathrm{\smash{#1}}}}}}
\renewcommand{\vec}[1]{\underline{\smash{#1}}}

\usepackage{overpic}

\usepackage{multirow}

\usepackage{scalerel}

\usepackage{todonotes}
\newcommandx{\john}[2][1=]{\todo[linecolor=red,backgroundcolor=red!25,bordercolor=red,#1]{#2}}
\newcommandx{\johnin}[2][1=]{\todo[inline,linecolor=red,backgroundcolor=red!25,bordercolor=red,#1]{#2}}
\newcommandx{\noah}[2][1=]{\todo[linecolor=blue,backgroundcolor=blue!25,bordercolor=blue,#1]{#2}}
\newcommandx{\noahin}[2][1=]{\todo[inline,linecolor=blue,backgroundcolor=blue!25,bordercolor=blue,#1]{#2}}
\newcommandx{\kristian}[2][1=]{\todo[linecolor=green,backgroundcolor=green!25,bordercolor=green,#1]{#2}}
\newcommandx{\kristianin}[2][1=]{\todo[inline,linecolor=green,backgroundcolor=green!25,bordercolor=green,#1]{#2}}

\newlist{thmlist}{enumerate*}{1}
\setlist[thmlist]{label=\rm{(\arabic*)}}

\crefname{thmlisti}{thm}{thms}
\creflabelformat{thmlisti}{#2#1#3}

\usepackage{cancel}

\usepackage[normalem]{ulem}

\ifpdf
\hypersetup{
  pdftitle={Regularisation of isolated codimension-2 discontinuity sets},
  pdfauthor={N. Cheesman, K. Uldall Kristiansen, S.J. Hogan}
}
\fi

\begin{document}

\maketitle

\begin{abstract}
Work on \textit{standard} piecewise-smooth (PWS) dynamical systems, with codimension-1 discontinuity sets,  relies on the Filippov framework, which does not always readily generalise to systems with higher codimension discontinuities.
These higher order degeneracies occur in applications, with spatial Coulomb friction being the prominent example. 

In this work, we consider PWS systems with \textit{isolated} codimension-2 discontinuity sets using regularization and blowup to study the dynamics. We present a general framework and give specific results for the local dynamics in a class of problems, and generalising Filippov sliding, crossing and sliding vector fields.
\end{abstract}


\begin{AMS}
  { 34E15, 37N15, 74H35
}\end{AMS}

\section{Introduction}
Piecewise-smooth (PWS) dynamical systems with discontinuity sets (along which the dynamics switches) are widely used to model systems with impact, dry friction, juddering, and buckling and in the study of relay control systems, mechanical systems etc \cite{bernardo2008piecewise,jeffrey2018hidden}.

When studying planar motion of rigid bodies with point contact, Coulomb friction results in a codimension-1 discontinuity set, where the relative velocity between two bodies is zero (see \cref{fig:fric1}).
But the existing framework for the study of these systems does not generalise to spatial motion of rigid bodies with point contact.
In this case, Coulomb friction results in a codimension-2 discontinuity set \cite{Antali2017,antali2019nonsmooth}, where both components of the relative velocity between two bodies are simultaneously zero (see \cref{fig:fric2}). Whilst friction is the primary motivation for this work, isolated codimension-2 discontinuity sets also occur in unit vector control design. For a brief introduction to unit control see \cite[\S 3.5]{utkin1999sliding}, and for an example of a resulting codimension-2 problem see \cite{levant2017quasi}.

\newlength{\twosubht}
\newsavebox{\twosubbox}
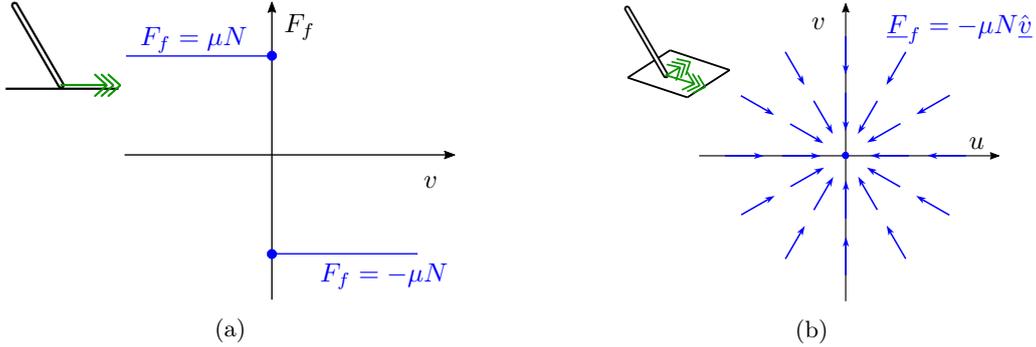
\begin{figure}[htbp]
\sbox\twosubbox{%
  \resizebox{\dimexpr.7\textwidth-1em}{!}{%
    \includegraphics[height=3cm]{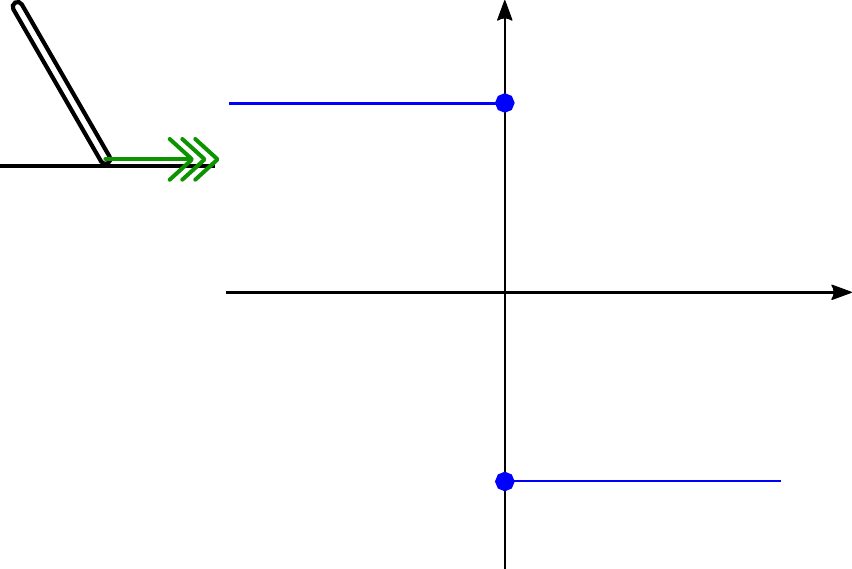}%
    \includegraphics[height=3cm]{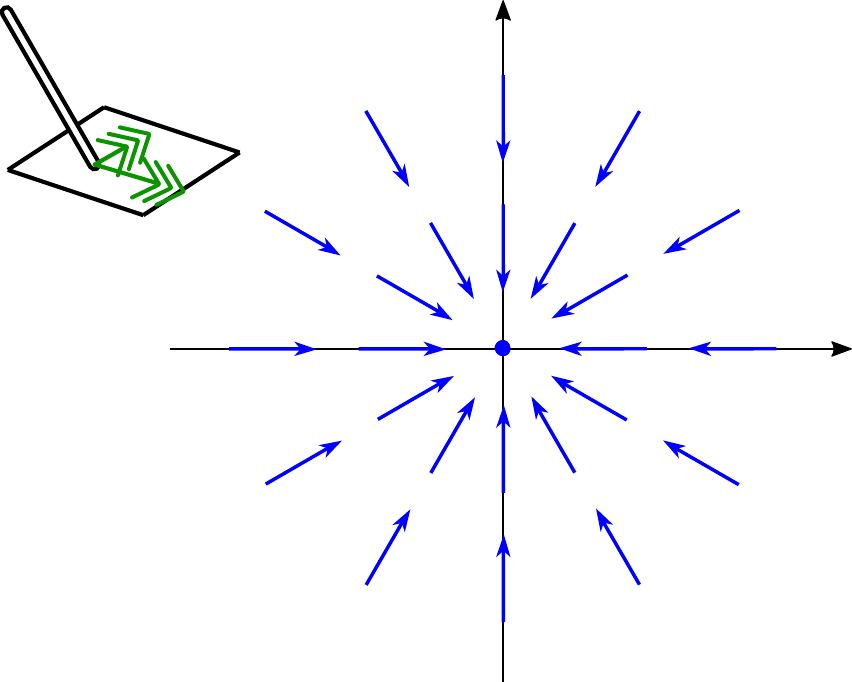}%
  }%
}
\setlength{\twosubht}{\ht\twosubbox}
\centering
\begin{subfigure}{0.42\textwidth}
\centering
	\begin{overpic}[height=\twosubht]{Figs/linear_coulomb2}
		\put(93,25){$v$}
		\put(62,60){$F_f$}
		\put(30,57){{\color{blue}{$F_f=\mu N $}}}
		\put(70,4){{\color{blue}{$F_f=-\mu N $}}}
	\end{overpic}\caption{\label{fig:fric1}
	}
\end{subfigure}
\hfil
\begin{subfigure}{0.42\textwidth}
\centering
\begin{overpic}[height=\twosubht]{Figs/spatial_coulomb2}
		\put(92,41){$u$}
		\put(50,73){$v$}
		\put(70,72){{\color{blue}{$\vec{F}_f=-\mu N \hat{\vec{v}}$}}}
	\end{overpic}
\caption{\label{fig:fric2}}
\end{subfigure}
\caption{The discontinuity sets for Coulomb friction with (a) one and (b) two degrees of freedom in the relative velocities (shown in green) between objects, where $\mu$ is the coefficient of friction and $N$ is the normal reaction force. In {(a)}, we have a planar rigid body with point contact, where Coulomb friction is given by $F_f=-\mu N \sign{(v)}=-\mu N \frac{v}{|v|}$ and there is a codimension-1 discontinuity at $v=0$. In {(b)}, we have a spatial rigid body system with point contact where Coulomb friction is given by $\vec{F}_f=-\mu N \hat{\vec{v}}=-\mu N \frac{\vec{v}}{|\vec{v}|}$ where $\vec{v}=(u,v)$ and there is a codimension-2 discontinuity at $u=v=0$.}
\label{fig:friction}
\end{figure}

We define our codimension-2 discontinuity problem as consisting of ordinary differential equations
\begin{equation}\label{eq:codim2}
	\dot{\vec{x}}=\vec{F}(\vec{x}),\,\vec{x}\in \Real^n\backslash \Sigma
\end{equation}
where $\vec{F}$ is a vector field that is smooth and well defined everywhere except on a connected \textit{codimension-2} set $\Sigma$ and $\dot{(\cdot)}=\frac{\mathrm{d}}{\mathrm{d}t}\left(\cdot\right)$ denotes differentiation with respect to time $t$. Furthermore, the vector field $\vec{F}$ has a well defined directional limit onto $\Sigma$ for each angle of approach (as with spatial Coulomb friction in \cref{fig:fric2}).
 
\begin{figure}[htbp]
	\centering
	\hfill
	\begin{subfigure}[b]{0.40\textwidth	}
		\begin{overpic}[width=\textwidth]{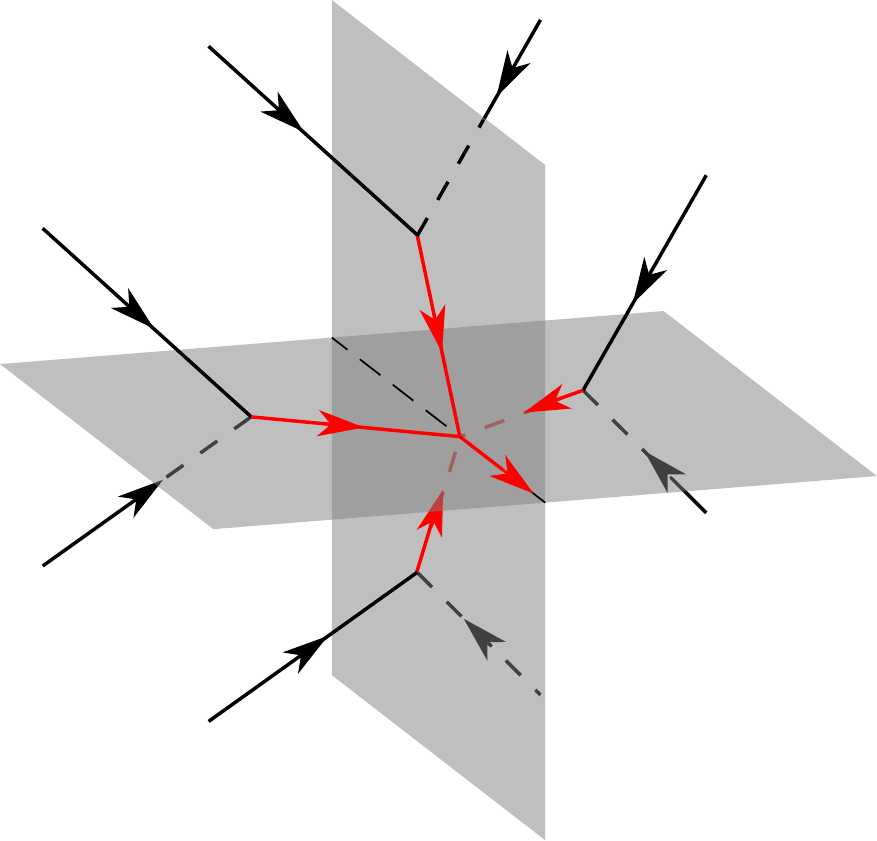}
			\put(5,85){$\mathcal{Q}_1$}
			\put(5,10){$\mathcal{Q}_3$}
			\put(83,85){$\mathcal{Q}_2$}
			\put(83,10){$\mathcal{Q}_4$}
			{\color{blue}
			\put(95,41.5){\line(-1,-3){3}}
			\put(89,27){$\Sigma_{2,4}$}
			\put(50,10){\line(-1,-1){8}}
			\put(30,0){$\Sigma_{3,4}$}
			\put(8,48){\line(-1,2){4}}
			\put(0,58){$\Sigma_{1,3}$}
			\put(62,77){\line(1,2){6}}
			\put(65,90){$\Sigma_{1,2}$}}
		\end{overpic}
		\caption{ \label{fig:2codim1}}
	\end{subfigure}
	\hfill
	\begin{subfigure}[b]{0.40\textwidth	}
		\begin{overpic}[height=\textwidth]{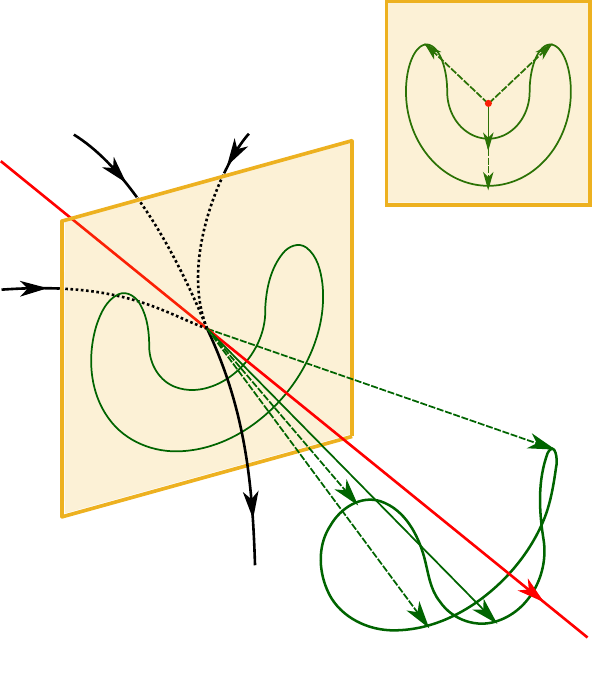}
			{\color{blue}{
			\put(8,92){$\Sigma$}\put(5.5,72){\line(1,4){4.75}}}}
		\end{overpic}
		\caption{ \label{fig:1codim2}}
	\end{subfigure}
	\hfill
	\caption{The difference between (a) the intersection of two codimension-1 discontinuity sets and (b) an isolated codimension-2 discontinuity set.
In {(a)} trajectories reach the intersection of the two discontinuities along the codimension-1 discontinuities $\Sigma_{i,j}$. In {(b)} trajectories reach the codimension-2 discontinuity set $\Sigma$ from a variety of different directions that do not correspond to codimension-1 discontinuity sets. The directional limit of $F$ onto $\Sigma$ is also shown (inset). 
\label{fig:codim1vcodim2}}
\end{figure}
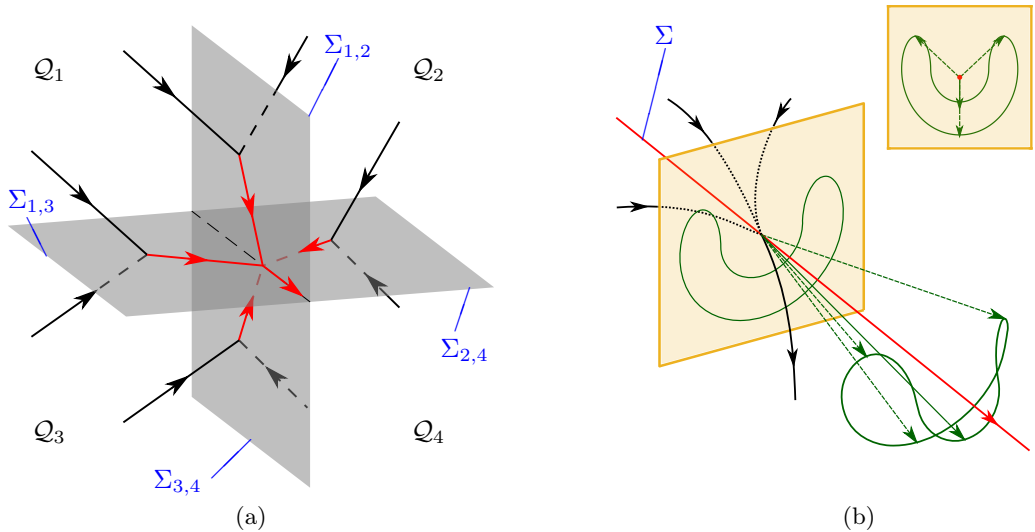

This type of {\it nonsmooth} system is related to standard PWS dynamical systems \cite{filippov2013differential}, which consist of finitely many ordinary differential equations 
\begin{equation}\label{eq:codim1}
	\dot{\vec{x}}=\vec{F}_i(\vec{x}),\,\vec{x}\in \mathcal{Q}_i\in\Real^n
\end{equation}
where each $\vec{F}_i$ is a smooth vector field.
Regions $\mathcal{Q}_i$ are open sets separated by a codimension-1 discontinuity set $\Sigma_{i,j}$ at the boundary between $\mathcal{Q}_i$ and $\mathcal{Q}_j$. 
When two of these codimension-1 discontinuity sets intersect transversally, a codimension-2 discontinuity results \cite{dieci2017moments,dieci2013filippov,dieci2011sliding,filippov2013differential,Jeffrey2014,jeffrey2018hidden,kaklamanos2019regularization}.
As noted in \cite{Antali2017}, we must distinguish between these sorts of codimension-2 problems and isolated codimension-2 problems described by \cref{eq:codim2}, as demonstrated in \cref{fig:codim1vcodim2}.

\subsection{Deficiencies of the existing framework}
\begin{figure}[htbp]
	\centering
	\begin{subfigure}[t]{0.45\textwidth}
	\centering
		\begin{overpic}[width=0.9\textwidth]{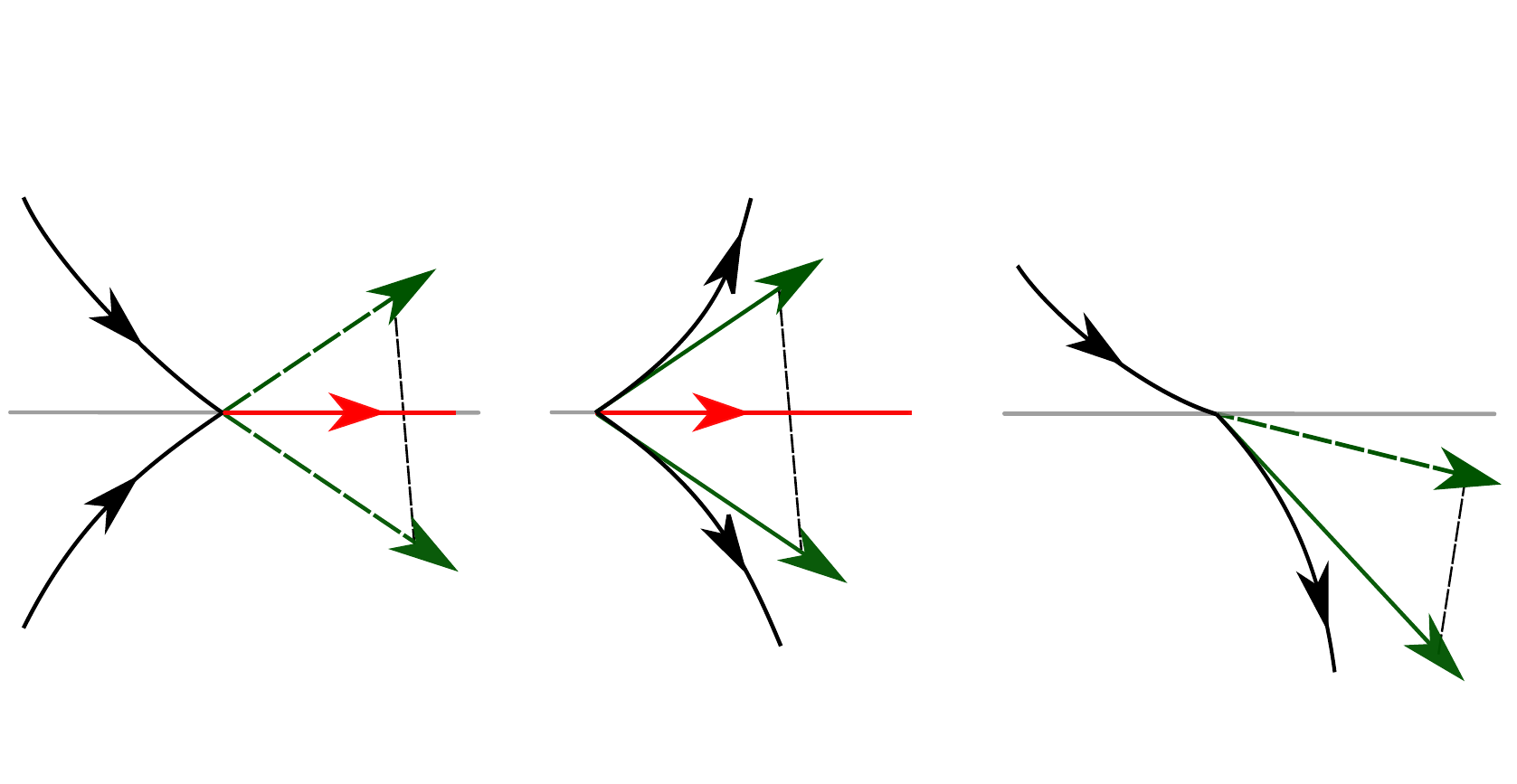}
		\put(5,4){Stable}
		\put(35,4){Unstable}
		\put(-3,23){{\color{blue}$\Sigma$}}
		\put(33,23){{\color{blue}$\Sigma$}}
		\put(100,23){{\color{blue}$\Sigma$}}
		\put(22,40){Sliding}
		\put(70,40){Crossing}
		\end{overpic}
		\caption{}\label{fig:standard_sliding}
	\end{subfigure}
	\hfil
	\begin{subfigure}[t]{0.45\textwidth}
		\centering
		\begin{overpic}[width=0.85\textwidth]{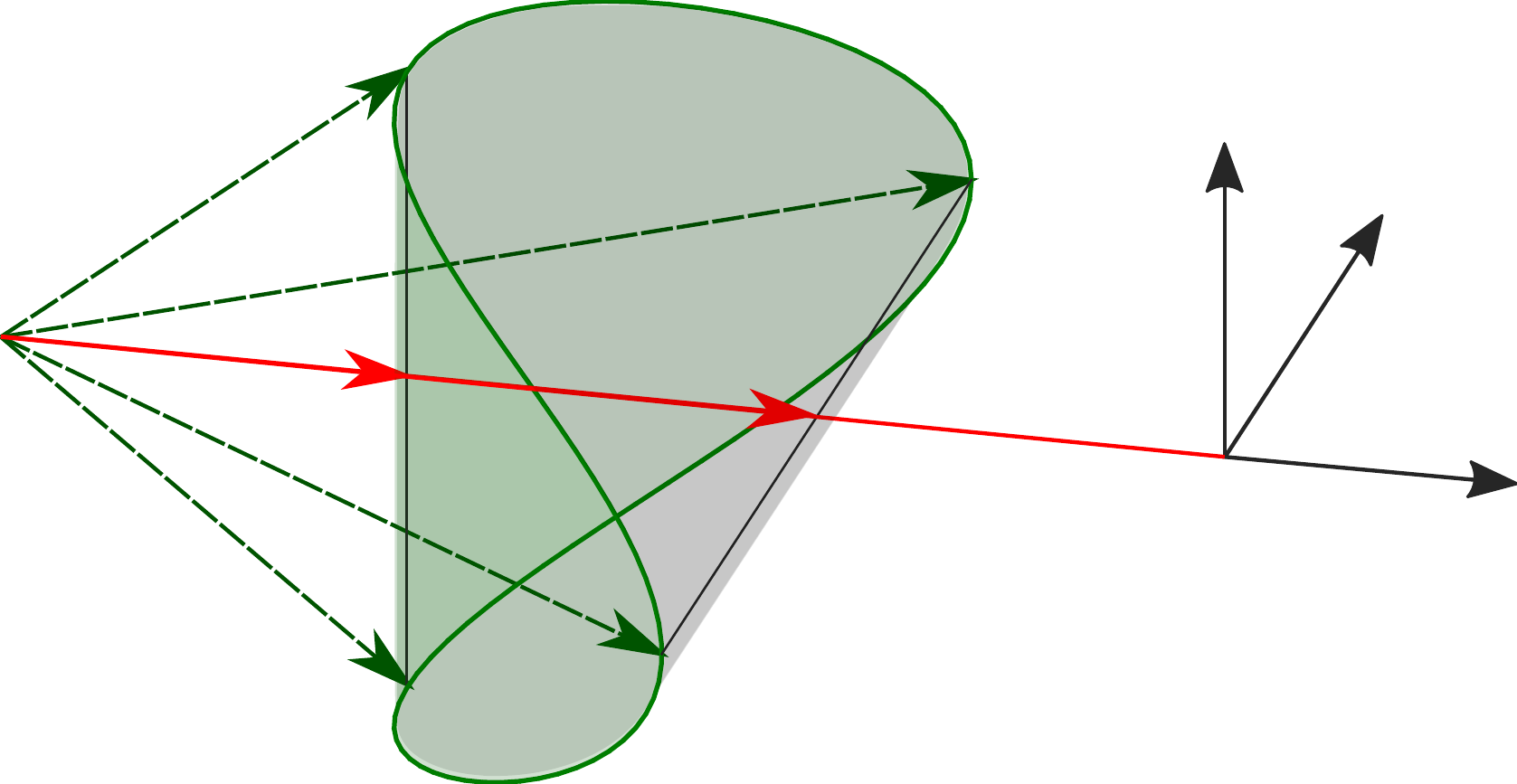}
			\put(95,15){$z$}
			\put(92,33){$x$}
			\put(75,43){$y$}
			\put(66,41){$\vec{v}_1$}
			\put(17,47){$\vec{v}_2$}
			\put(17,23){$\vec{v}_{2,4}$}
			\put(54,19){$\vec{v}_{1,3}$}
			\put(19,3){$\vec{v}_4$}
			\put(45,4){$\vec{v}_3$}
			\put(70,16){{\color{blue}$\Sigma$}}
		\end{overpic}
		\caption{ }
		\label{fig:pringle}
	\end{subfigure}
	\\
	\begin{subfigure}[t]{0.45\textwidth}
		\begin{overpic}[width=0.85\textwidth]{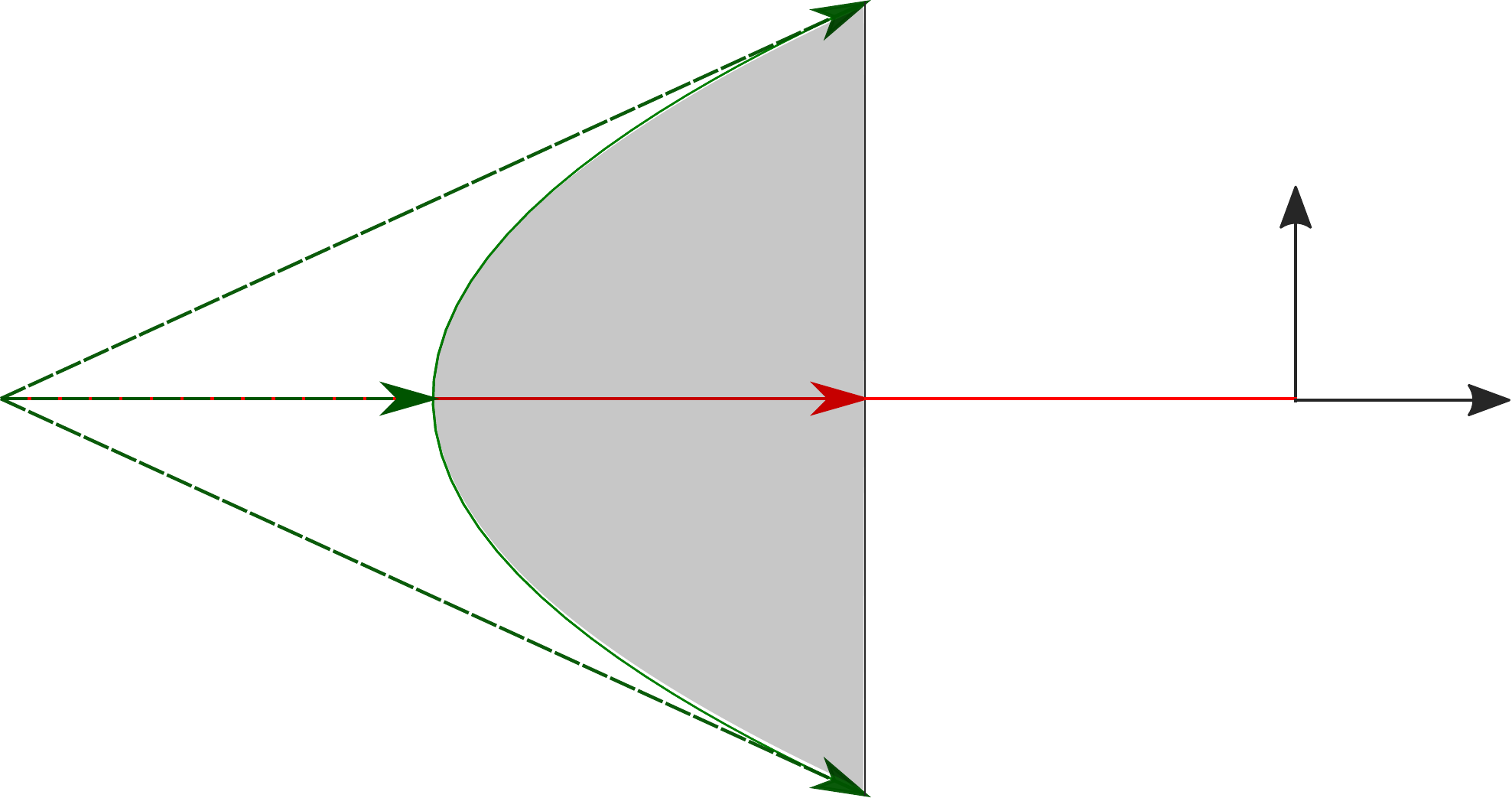}
			\put(95,28){$z$}
			\put(80,40){$x$}
			\put(60,50){$\vec{v}_1$}
			\put(60,0){$\vec{v}_3$}
			\put(60,30){$\vec{v}_{1,3}$}
			\put(60,0){$\vec{v}_3$}
			\put(20,21){$\vec{v}_2$}
			\put(20,29){$\vec{v}_4$}
			\put(73,18){{\color{blue}$\Sigma$}}
		\end{overpic}
		\caption{ }\label{fig:pringle1}
	\end{subfigure}
	\hfil
	\begin{subfigure}[t]{0.45\textwidth}
		\centering
		\begin{overpic}[width=0.85\textwidth]{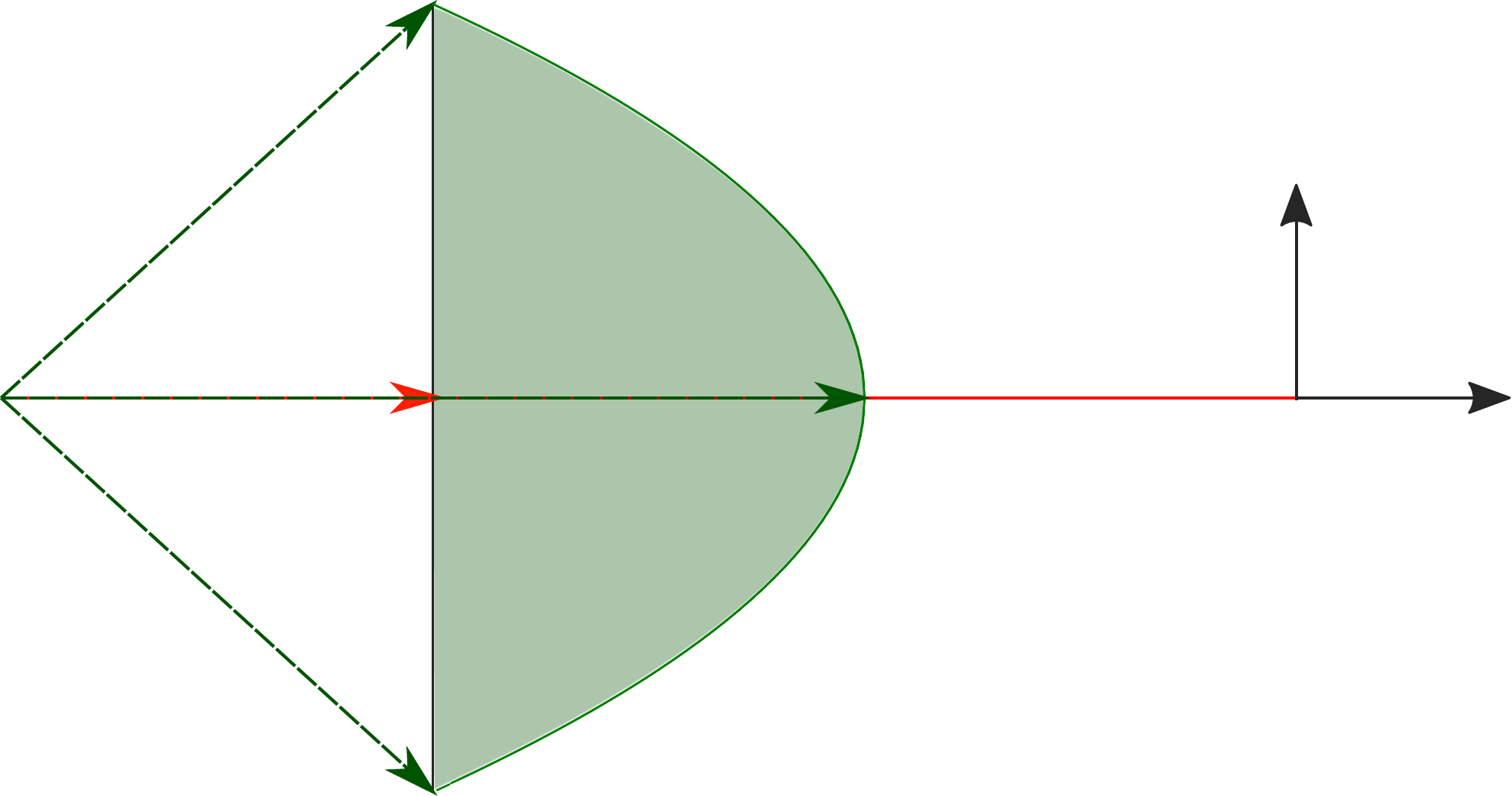}
			\put(95,28){$z$}
			\put(80,40){$y$}
			\put(19,50){$\vec{v}_2$}
			\put(19,0){$\vec{v}_4$}
			\put(18,30){$\vec{v}_{2,4}$}
			\put(60,29){$\vec{v}_3$}
			\put(73,18){{\color{blue}$\Sigma$}}
		\end{overpic}
		\caption{}\label{fig:pringle2}
	\end{subfigure}

	\caption{Demonstration of the generic non-uniqueness of the Filippov sliding vector field for codimension-2 discontinuity sets. In (a) the convex combination of the incident vectors is one-dimensional and its intersection with the codimension-1 discontinuity set $\Sigma$ (where it exists) is generically unique. In (b) we see that, the convex hull of the incident vectors is typically $n$-dimensional and its intersection with the codimension-2 discontinuity set $\Sigma$ is  $(n-2)$-dimensional. Hence, there are uncountably many possible vectors that are both a convex combination of the incident vectors $\vec{v}_i$ and tangential to $\Sigma$. By way of illustration, we show two candidate sliding vectors: $\vec{v}_{1,3}$ and $\vec{v}_{2,4}$. In (c) we see the projection of {(b)} onto the $(x,z)$ plane. The convex combination of $\vec{v}_1$ and $\vec{v}_3$ that is tangent to $\Sigma$ (the $z$-axis), is the candidate sliding vector $\vec{v}_{1,3}$. In (d) we see the projection of {(b)} onto the $(y,z)$ plane. Applying the Filippov convention to $\vec{v}_2$ and $\vec{v}_4$ gives the candidate sliding vector $\vec{v}_{2,4}\neq \vec{v}_{1,3}$. }\label{fig:generic_nonuniqueness_of_sliding}
\end{figure}
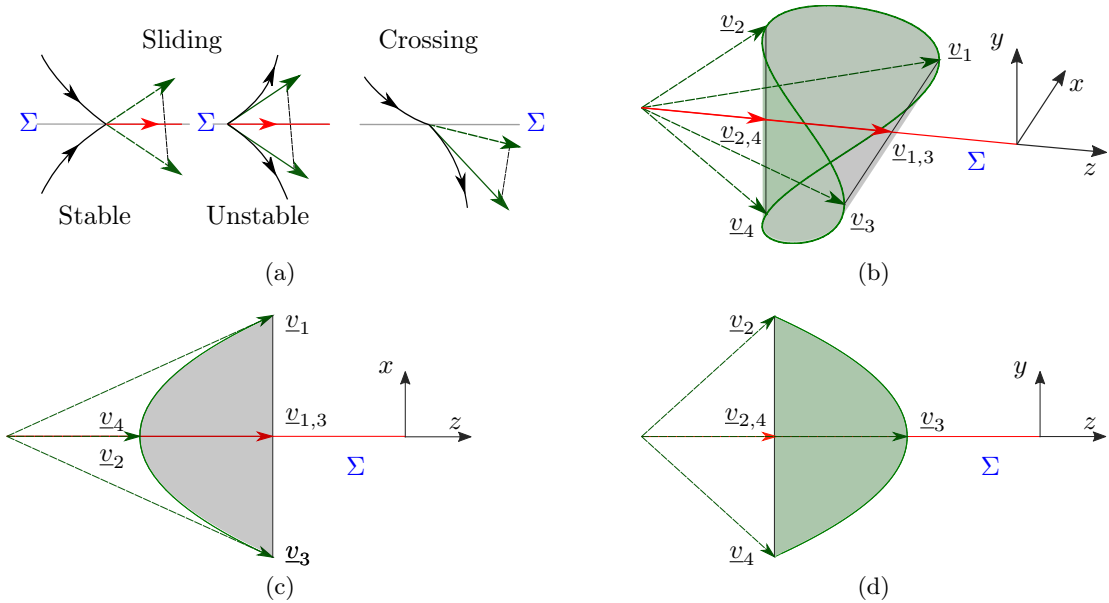

Standard PWS systems with codimension-1 discontinuity sets display two generic types of behaviour: crossing and sliding \cite{filippov2013differential}. Crossing happens when trajectories pass through the discontinuity set without remaining on it for any period of time. Sliding occurs when trajectories reach the discontinuity set and continue along it. In many applications, it is then necessary to define a \textit{sliding vector field} (the flow along this discontinuity set), using the \textit{Filippov convention}. 
These behaviours are shown in  \cref{fig:standard_sliding}. If both vector fields are pointing into (or out from) the switching surface, we have sliding. If one vector field points in whilst the other points out, we have crossing.

With codimension-2 discontinuity sets, however, these definitions do not generalise intuitively.
In \cref{fig:1codim2} it is not clear whether trajectories that reach $\Sigma$ should ``slide'' along it or ``cross'' and leave it.
Even if we can determine that there is sliding, how do we prescribe the sliding flow?

Motivated by Coulomb friction, Antali and St\'{e}p\'{a}n \cite{Antali2017} gave a generalisation of the framework for standard (codimension-1) PWS systems to (codimension-2) {\it extended Filippov systems}.
They used the Filippov convention to construct the sliding vector field: a convex combination of the vectors incident to the discontinuity set, that  is tangent to the discontinuity set. They noted that this convention does not generically give a unique sliding vector in the case of a codimension-2 problem. 

Let us first consider a standard PWS system (\cref{fig:standard_sliding}).
The Filippov convention provides a unique sliding vector; the convex combination of two incident vectors spans a one-dimensional line section and its intersection with the plane tangent to the discontinuity set $\Sigma$ is therefore a point (or the empty set in the case of crossing).
Let us now consider the same definition when applied to a codimension-2 problem (\cref{fig:pringle}). Here the codimension-2 discontinuity set $\Sigma$ coincides with the $z$-axis and we show incident vectors: $\vec{v}_i$ $i\in\{1\ldots4\}$.
The convex combination of $\vec{v}_1$ and $\vec{v}_3$ gives one candidate $\vec{v}_{1,3}$ for the sliding vector (\cref{fig:pringle1}).
But, the convex combination of $\vec{v}_2$ and $\vec{v}_4$ gives another candidate sliding vector $\vec{v}_{2,4}\neq \vec{v}_{1,3}$ (\cref{fig:pringle2}).
In fact, there are uncountably many candidates for the sliding vector field that satisfy the Filippov convention, given here by $\{\lambda \vec{v}_{1,3}+(1-\lambda)\vec{v}_{2,4}|\lambda\in [0,1] \}$. This non-uniqueness should be expected unless all the incident vectors are coplanar.
The set defined by the convex combination of incident vectors will generically be the same dimension as the space;
the intersection of that set with the discontinuity set will therefore be codimension-2.

In this paper, we use slow-fast theory, geometric singular perturbation theory (GSPT) and blowup to study extended Filippov systems. This approach not only addresses the generic non-uniqueness (\cref{fig:pringle}) and the ambiguity about the definitions of sliding  and crossing (\cref{fig:1codim2}), but also allows the use of the powerful and well-understood methods from smooth dynamical systems. We will study smooth but sharp systems as perturbations away from the nonsmooth limit. 
An understanding of the relationship between a PWS system \cite{carmona2008existence} and a (stiff) smooth system that  approximates it \cite{michelson1986steady} is important, since it can reveal the robustness (or otherwise) of the PWS system to smoothing perturbations, see \cite{kristiansen2015use}.

We regularise the nonsmooth system \cref{eq:codim2}, viewing it as the (suitably defined) limit of a smooth system \cite{teixeira2012regularization}. The resulting singularly perturbed smooth system has a hidden slow fast structure in the smoothing parameter $\varepsilon$ and GSPT can be applied. Fenichel's theorem \cite{fenichel1979geometric} states that this smooth problem lies $\varepsilon$-close to the nonsmooth limit. So we use blowup \cite{dumortier1996canard, krupa2001extending} to connect the dynamics approaching the discontinuity with the dynamics along it.
This approach has proven useful in both the study of regularisations of codimension-1 discontinuity sets \cite{teixeira2012regularization,kristiansen2015use,kristiansen2019resolution,kristiansen2020regularized}, and their intersections \cite{kaklamanos2019regularization} (\cref{fig:2codim1}). In fact, here we follow a very similar procedure to \cite{kaklamanos2019regularization}, but now for isolated codimension-2 discontinuity sets. 

\subsection{Outline}
In \cref{sec:prelims}, we introduce our notation and formalise what we mean by an isolated codimension-2 discontinuity.

Then in \cref{sec:gen}, we apply our GSPT approach to this general form of codimension-2 problem, which serves to demonstrate the procedure. In \cref{sec:lin}, we study a more restrictive class of problems where we find analytic results, allowing us to classify possible phase portraits.

In \cref{sec:examples}, we provide examples of the method's use for physical applications, and illustrate other phenomena. 

\section{Preliminaries}\label{sec:prelims}
\begin{figure}[htbp]
\centering
\begin{overpic}[width=0.9\textwidth]{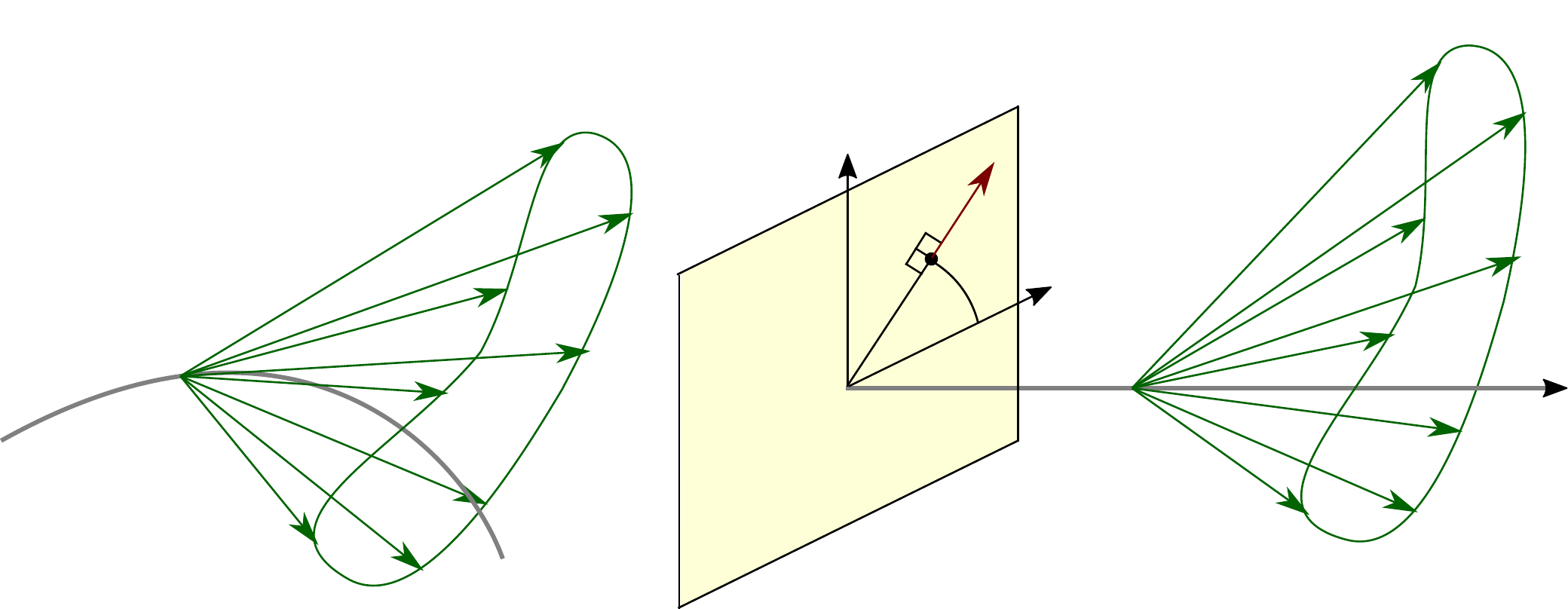}
{\color{matlab4}
\put(45,32){\bezier{0}(0,0)(15,5)(30,0)}
\put(75,32){\vector(2,-1){3}}
\put(55,36){Straightening}}
{\color{blue}
\put(0.55,14){$\Sigma$}
\put(2.8,15){\vector(2,-1){3}}
\put(69,8){$\Sigma$}
\put(70,11){\vector(-1,1){2.8}}
}
{\color{black}
\put(49,28){$y$}
\put(64,21){$x$}
\put(95,15){$\vec{z}$}
{
\put(52.5,23.5){\footnotesize$\begin{pmatrix}x\\ y\end{pmatrix}$}
}
}
{\color[RGB]{123,0,0}
\put(55.75,27.25){\rotatebox{30}{\footnotesize$\vec{e}(x,y)$}}
}
\put(55,18){$\theta$}
\put(88.5,6){$\vec{F}^*(\theta,\vec{z})$}
\end{overpic}
\caption{Straightening of $\Sigma$ so that $\Sigma=\left\{x=0,y=0,\vec{z}\in\mathbb{R}^{n-2}\right\}$ 
and the geometric interpretation of $\vec{e}$.}\label{fig:straight}
\end{figure}

We consider the system \cref{eq:codim2}, with a codimension-2 discontinuity set $\Sigma$ of a vector field $\vec{F}$. We straighten out $\Sigma$ at least locally so that $\Sigma=\{x=0,y=0,\vec{z}\in\mathbb{R}^{n-2}\}$ or a subset thereof. Hence we have $\vec{F}({x},y,\vec{z})$ where $x\in\mathbb{R}$, $y\in\mathbb{R}$ and $\vec{z}\in\mathbb{R}^{n-2}$. We assume the following. 
\begin{assumption}\label{ass:v}
$\vec{F}$ takes the form
\begin{equation}
\vec{F}(x,y,\vec{z})\equiv \vec{V}(\vec{e}(x,y),x,y,\vec{z})
\end{equation}
where $\vec{V}: S^1\times \mathbb{R}^2\times \mathbb{R}^{n-2}\rightarrow \mathbb{R}^n$ is smooth in all its entries, and $\vec{e}: \mathbb{R}^2\setminus\{0\}\rightarrow S^1$ given by 
\begin{equation}\label{eq:M}
\vec{e}(x,y):=\left(\frac{x}{\sqrt{x^2+y^2}},\frac{y}{\sqrt{x^2+y^2}}\right)^\intercal ,
\end{equation}
which corresponds to the unit vector pointing radially away from the origin at $(x,y)^\intercal\neq \vec{0}$.
\end{assumption}
The straightening of $\Sigma$ and the geometrical interpretation of $\vec{e}$ are given in \cref{fig:straight}.
\begin{remark}
There is a jump in $\vec{e}(x,y)$ at the origin. If we write $\vec{e}$ along a line, which passes through the origin at a fixed angle $\theta$, given by
\begin{equation}\label{eq:linethroughorigin}
(x,y)^\intercal=l(\cos{\theta},\sin{\theta})^\intercal,\quad  l\in\mathbb{R},
\end{equation}
then
\begin{equation}\label{eq:eline}
\vec{e}(l\cos{\theta},l\sin{\theta})=\sign{(l)}(\cos{\theta},\sin{\theta})^\intercal.
\end{equation}
Along this line $\vec{e}$ results in a standard piecewise switch. From \cref{eq:eline}, we note that for both $l>0$ and ${l\to 0^{+}} $, we have 
\begin{equation}\label{eq:cossin}
\vec{e}(l\cos{\theta},l\sin{\theta})\equiv(\cos{\theta},\sin{\theta})^\intercal.
\end{equation}
\end{remark}

From \cref{ass:v} it follows that there is a well-defined directional {\it limit vector field} $\vec{F}^*(\theta,\vec{z})$ given by
\begin{equation}\label{eq:directional_limit}
\vec{F}^*(\theta,\vec{z}):=\lim_{\epsilon\to0} \vec{F}(\epsilon\cos\theta,\epsilon\sin\theta,\vec{z})=\vec{V}((\cos{\theta}\,,\,\sin{\theta})^\intercal,0,0,\vec{z}),
\end{equation}
that is smooth in both $\theta$ and $\vec{z}$.

Hence, from now on, we consider the system
\begin{equation}\label{eq:vf}
(\dot{x},\dot{y},\dot{\vec{z}})^\intercal=
\vec{V}\left(\vec{e}(x,y),x,y,\vec{z}\right),
\end{equation} 
where  $\vec{z}\in\mathbb{R}^{m},m=n-2>0$. It will often be useful to split \cref{eq:vf} into two coupled ordinary differential equations 
\begin{subequations}\label{eq:gen_odes}
\begin{align}
(\dot{x},\dot{y})^\intercal&=\vec{U}\left(\vec{e}(x,y),x,y,\vec{z}\right),\label{eq:gen_odes1}\\
\dot{\vec{z}}&=\vec{W}\left(\vec{e}(x,y),x,y,\vec{z}\right).\label{eq:gen_odes2}
\end{align}
\end{subequations}
Equation \cref{eq:gen_odes1} corresponds to dynamics normal to $\Sigma$, and equation \cref{eq:gen_odes2} corresponds to dynamics tangent to $\Sigma$.

\subsection{Dynamics near the discontinuity set}\label{sec:nonsmoothdynamics}
Following \cite{Antali2017}, we will now study \cref{eq:gen_odes} using polar coordinates $(\rho,\theta)$, where 
\begin{equation}\label{eq:polar_coordinates}
x=\rho\cos{\theta}\quad \text{and}\quad y=\rho\sin{\theta},
\end{equation}
with the aim of describing how solutions reach $\Sigma$. Hence, using \cref{eq:cossin},  \cref{eq:gen_odes} becomes
\begin{equation}
\begin{split}\label{eq:polar_gen_odes}
\begin{pmatrix}
\dot{\rho}\\ \dot{\theta}
\end{pmatrix}&=\begin{pmatrix}
1&0\\
0&\frac{1}{\rho}
\end{pmatrix}\mat{R}(\theta)^\intercal\,\vec{U}\left( (\cos{\theta},\sin{\theta})^\intercal,\rho\cos{\theta},\rho\sin{\theta},\vec{z}\right),\\
\dot{\vec{z}}&= \vec{W}\left((\cos{\theta},\sin{\theta})^\intercal,\rho\cos{\theta},\rho\sin{\theta},\vec{z}\right),
\end{split}
\end{equation}
where $\mat{R}(\theta)$ is the rotation matrix
\begin{equation}\label{eq:rotation_matrix}
\mat{R}(\theta):=\begin{pmatrix}
\cos{\theta}& -\sin{\theta} \\
\sin{\theta}&\phantom{-}\cos{\theta}
\end{pmatrix}.
\end{equation}
Notice that along $\Sigma$, which corresponds to $\rho=0$, the dynamics is not defined. We therefore transform to a new time $\mathcal{T}$, given by $\mathrm{d}\mathcal{T}=\frac{1}{\rho}\mathrm{d}t$, so that \cref{eq:polar_gen_odes} becomes
\begin{equation}
\label{eq:polar_gen_odes_desing}
\begin{split}
\frac{\mathrm{d}}{\mathrm{d}\mathcal{T}}\begin{pmatrix}
{\rho}\\ {\theta}
\end{pmatrix}&=\begin{pmatrix}
\rho&0\\
0&1
\end{pmatrix}\mat{R}(\theta)^\intercal\,\vec{U}\left((\cos{\theta},\sin{\theta})^\intercal,\rho\cos{\theta},\rho\sin{\theta},\vec{z}\right),\\
\frac{\mathrm{d}}{\mathrm{d}\mathcal{T}}\,{\vec{z}}&= \rho \,\vec{W}\left((\cos{\theta},\sin{\theta})^\intercal,\rho\cos{\theta},\rho\sin{\theta},\vec{z}\right).
\end{split}
\end{equation}
With this transformation, \cref{eq:polar_gen_odes_desing} is well-defined along $\Sigma$ and orbits are preserved. Note that trajectories can approach $\rho=0$ in finite time in \cref{eq:polar_gen_odes}, whereas trajectories can only reach $\rho=0$ in infinite time in \cref{eq:polar_gen_odes_desing}. 

Since $\rho=0$ is an invariant manifold of \cref{eq:polar_gen_odes_desing}, equilibria of \cref{eq:polar_gen_odes_desing} exist when the $\theta$ component,
\begin{align}
\Theta\left(\rho,{\theta},\vec{z}\right)&:=(-\sin\theta,\cos\theta)\, \vec{U}\left((\cos{\theta},\sin{\theta})^\intercal,\rho\cos{\theta},\rho\sin{\theta},\vec{z}\right),\label{eq:Theta1_}\\
\intertext{is simultaneously zero. The function $\Theta\left(\rho,{\theta},\vec{z}\right)$ is smooth at $\rho=0$ and}
\Theta\left(0,{\theta},\vec{z}\right)&=(-\sin\theta,\cos\theta)\, \vec{U}\left((\cos{\theta},\sin{\theta})^\intercal,0,0,\vec{z}\right).\label{eq:Theta10_}
\end{align}

Let there exist $\theta_0,\vec{z}_0$ such that $\Theta(0,\theta_0,\vec{z}_0)=0$, then
\begin{equation}\label{eq:equator_equilibria}
P_0=(0,\theta_0,\vec{z}_0)
\end{equation}
is an equilibrium of \cref{eq:polar_gen_odes_desing}.

The importance of these equilibria along $\rho=0$ has previously been identified in \cite{Antali2017} where they were called \textit{limit directions}. These are the directions along which trajectories reach (or leave) the discontinuity set $\Sigma$. It is useful to study the stability of these equilibria. Linearising around $P_0$, we find that the Jacobian is given by

\begin{align}
\mat{J}_0(\theta_0,\vec{z}_0)&:=
\renewcommand*{\arraystretch}{1.6}
\begin{pmatrix}
\lambda_\rho(\theta_0,\vec{z_0})&0&\mat{0}_{1\times m}\\
\left.\pdiff{\Theta(0,\theta_0,\vec{z}_0)}{\rho}\right.&\left.\pdiff{\Theta(0,\theta_0,\vec{z}_0)}{\theta}\right.&\left.\pdiff{\Theta(0,\theta_0,\vec{z}_0)}{\vec{z}}\right.\\
\vec{W}((\cos\theta_0,\sin\theta_0)^\intercal,0,0,\vec{z}_0)&\mat{0}_{ m\times1}&\mat{0}_{ m\times m}
\end{pmatrix},
\intertext{which has only two non-zero eigenvalues, 
}
\lambda_\rho(\theta_0,\vec{z}_0)&:=(\cos{\theta_0} \, ,\, \sin{\theta_0} ) \, \vec{U}\left((\cos{\theta_0}, \sin{\theta_0} )^\intercal,0,0,\vec{z}_0 \right)\\
\intertext{and}
\lambda_\theta(\theta_0,\vec{z}_0)&:=\pdiff{\Theta(0,\theta_0,\vec{z}_0)}{\theta}.
\end{align}
By the stable manifold theorem \cite{perko}, if $\lambda_\rho(\theta_0,\vec{z}_0)<0$, then there is an orbit approaching $\Sigma$ along the direction $\theta=\theta_0$, and we say that the equilibrium $P_0$ is \textit{radially attracting}. Conversely, if $\lambda_\rho(\theta_0,\vec{z}_0)>0$, then there is an orbit leaving $\Sigma$ along $\theta=\theta_0$, and $P_0$ is \textit{radially repelling}. Similarly, if $\lambda_\theta(\theta_0,\vec{z}_0)<0$ then we call $P_0$ \textit{angularly attracting} and if $\lambda_\theta(\theta_0,\vec{z}_0)>0$ then we call $P_0$ \textit{angularly repelling}. If an equilibrium $P_0$ is both angularly and radially attracting, 
then there is a neighbourhood of $P_0$ that reaches $\Sigma$ at $P_0$ in finite time under the forward flow of \cref{eq:polar_gen_odes}. 

In \cite{Antali2017}, the authors define the sliding region to be the subset $\Sigma_{\mathrm{sl}}\subset \Sigma$ such that for each $\vec{z}_0\in \Sigma_{\mathrm{sl}}$, if there exist solutions $\theta_0$ such that $\Theta(0,\theta_0,\vec{z}_0)=0$, then $\lambda_\rho(\theta_0,\vec{z}_0)<0$ (they have an alternative definition when no solutions exist to $\Theta(\theta,\vec{z}_0)=0$).
We will discuss the suitability of this definition later. If there is sliding, a sliding vector field should be prescribed. Whilst \cite{Antali2017} used the Filippov convention, we will proceed with a different approach.

\subsection{Regularisation}\label{sec:regularisation}
Our strategy is to regularise the nonsmooth vector field $\vec{V}$ in \cref{eq:vf}, viewing it as the limit of a smooth one. To proceed, let us first define a regularisation function $\Psi(s)$. 
\begin{definition}\label{def:Psi}
We define a regularisation function $\Psi(s)$ that satisfies the conditions
\begin{enumerate}[label=\rm{(R\arabic*)}]
\item $\Psi(s)$ is smooth $\forall s \in [0,\infty)$ \label{R1}
\item $\Psi(0)=1$ \label{R2}
\item $\Psi(s)>0$ $\forall s\in [0,\infty)$ \label{R3}
\item $\Psi'(s)\leq 0$ $\forall s\in [0,\infty)$ \label{R4}
\item $\Psi_1(s)$ is smooth $\forall s \in [0,\infty)$, where
$$\Psi_1(s):=\begin{cases} \lim_{\sigma\to\infty}\Psi(\sigma), & s=0 \\
\Psi(1/s), & s >0 \end{cases}.$$ \label{R5}
\end{enumerate}
\end{definition}
For example, $\Psi(s)=(1+s)^{-1}$, $\Psi(s)=e^{-s}$ and $\Psi(s)=1$ all satisfy \cref{def:Psi}.
We use $\Psi$ in the regularisation of $\vec{e}$ given in \cref{eq:M}, as follows.
\begin{definition}\label{def:MPsi}
Let $D_1 \subset \mathbb{R}^2$ be the unit disc centred at the origin and let $\Psi$ be a regularisation function as in \cref{def:Psi}. Then we define $\vec{e}_\Psi:\mathbb{R}^2\rightarrow D_1$, a regularisation of $\vec{e}$, as $\vec{e}_\Psi(x,y;\varepsilon)$ given by
\begin{equation}\label{eq:ePsi}
\vec{e}_\Psi(x,y;\varepsilon):=\left(\frac{x}{\sqrt{x^2+y^2+\varepsilon^2 \Psi\left(\frac{x^2+y^2}{\varepsilon^2}\right)}},\frac{y}{\sqrt{x^2+y^2+\varepsilon^2 \Psi\left(\frac{x^2+y^2}{\varepsilon^2}\right)}}\right)^\intercal
\end{equation}
for $0<\varepsilon\ll1$.
\end{definition}

\begin{figure}[htbp]
\centering
\begin{overpic}[width=0.3\textwidth]{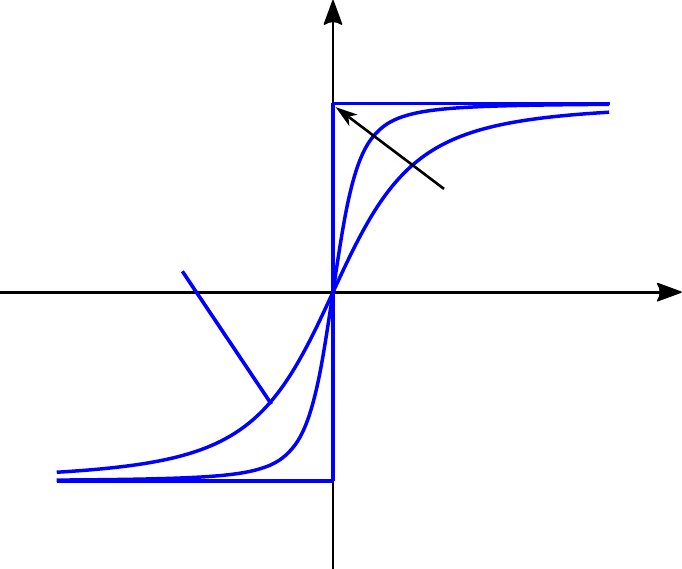}
\put(95,32){$l$}
\put(15,46){$L_\Psi(l;\varepsilon)$}
\put(67,53){$\varepsilon\to0$}
\put(41, 67){$1$}
\put(51, 10){$-1$}
\end{overpic}
\caption{The component of the map $\vec{e}_\Psi$ along any line through the origin $L_\Psi(l;\varepsilon)$ is a smooth function of $l$, whose limit is a $\sign$ function as $\varepsilon\to0$. For example, if we set $\Psi(s)=s/\sinh^2(\sqrt{s})$, then $L_\Psi(l;\varepsilon)=\tanh(l/\varepsilon)$.}\label{fig:component}
\end{figure}
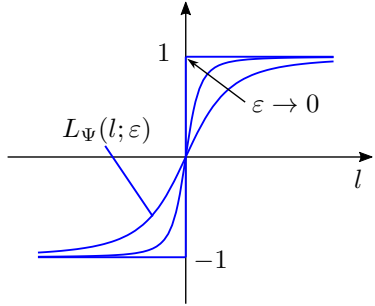

\begin{remark}
 $\vec{e}_\Psi$ does not jump discontinuously at the origin. When written along a line through the origin  \cref{eq:linethroughorigin}, we have
 \begin{align}
 \vec{e}_\Psi(l\cos{\theta},l\sin{\theta})=&L_\Psi(l;\varepsilon)(\cos{\theta},\sin{\theta})^\intercal,
 \intertext{where}
L_\Psi(l;\varepsilon):=&\frac{l}{\sqrt{l^2+\varepsilon^2\Psi\left({l^2}/{\varepsilon^2}\right)}}
\end{align}
is a smooth function whose limit is a $\sign$ function as $\varepsilon\to0$ (see \cref{fig:component}). 
\end{remark}

We note the following properties of $\vec{e}_\Psi$. 
\begin{lemma}\label{lem:ePsi}
Consider $\vec{e}_\Psi:\mathbb{R}^2\rightarrow D_1$ given by \cref{eq:ePsi}, then
\crefformat{enumi}{#2\textup{#1}#3}
\begin{enumerate}[label={\rm(\alph*)}]
\item{$\vec{e}_\Psi(x,y;\varepsilon)$ is equivariant under rotations
\begin{equation}\label{eq:rotationalequivariance}
\vec{e}_\Psi\left(\mat{R}(\phi)(x,y)^\intercal;\varepsilon\right) \equiv  \mat{R}(\phi)\vec{e}_\Psi(x,y;\varepsilon)
\end{equation}
for all $\phi\in [0,2\pi)$, where $\mat{R}(\phi)$ is a rotation matrix as in \cref{eq:rotation_matrix},\label{lem:ePsia}}
\item{$\vec{e}_\Psi(x,y;\varepsilon)$ satisfies the scaling property
\begin{equation}\label{eq:ePsiscale}
\vec{e}_\Psi(k x,k y;k \varepsilon )\equiv \vec{e}_\Psi(x,y;\varepsilon),
\end{equation}
for some $k>0$, that is, $\vec{e}_\Psi$ is a homogeneous function of degree 0 on $(x,y,\varepsilon)$,\label{lem:ePsib}}
\item{$\vec{e}_\Psi(x,y;1)$ is one-to-one,\label{lem:ePsic}}
\item{$\vec{e}_\Psi(x,y;\varepsilon)$ satisfies
\begin{equation}
\vec{e}_\Psi^\intercal \vec{e}_\Psi<1.\label{eq:ePsimod}
\end{equation}}
\end{enumerate}
\end{lemma}
The proof of \cref{lem:ePsi} can be found in \cref{sec:lemePsi}.

We now study the system
\begin{equation}
\begin{split}\label{eq:vfreg_ode}
(\dot{x},\dot{y},\dot{\vec{z}})^\intercal &= \vec{V}(\vec{e}_\Psi(x,y;\varepsilon),x,y,\vec{z})\\
\end{split}
\end{equation}
as a regularisation of \cref{eq:vf}.
In order to 
bring the singularly perturbed system
\cref{eq:vfreg_ode} 
into standard slow-fast form, we adopt the scalings 
\begin{equation}\label{eq:scalings}
x=\varepsilon x_2\quad \text{and}\quad y=\varepsilon y_2.
\end{equation}
As a result we find
\begin{equation}
\begin{split}\label{eq:vfreg_ode_sf2}
\varepsilon( \dot{x}_2,\dot{y}_2)^\intercal &= \vec{U}(\vec{e}_\Psi( x_2, y_2;1), \varepsilon x_2,\varepsilon y_2,\vec{z}),\\
\dot{\vec{z}}&=\vec{W}(\vec{e}_\Psi( x_2, y_2;1), \varepsilon x_2,\varepsilon y_2,\vec{z}),
\end{split}
\end{equation}
using the scaling property in \cref{lem:ePsi}\cref{lem:ePsib}. \cref{eq:vfreg_ode_sf2} is a standard slow-fast system written in slow time.
Setting $\varepsilon=0$ in \cref{eq:vfreg_ode_sf2} gives the reduced problem 
\begin{subequations}\label{eq:RP}
\begin{align}
\vec{0}&= \vec{U}(\vec{e}_\Psi( x_2, y_2;1), 0,0,\vec{z}),\label{eq:RP1}\\
\dot{\vec{z}}&=\vec{W}(\vec{e}_\Psi( x_2, y_2;1), 0,0,\vec{z})\label{eq:RP2}.
\end{align}
\end{subequations}
If we can solve \cref{eq:RP1}, finding
\begin{equation}\label{eq:RPCM}
\vec{e}_\Psi(x_2,y_2;1)=\left(c(\vec{z}),s(\vec{z})\right)^\intercal
\end{equation}
such that $\vec{U}(\left(c(\vec{z}),s(\vec{z})\right)^\intercal, 0,0,\vec{z})=\vec{0}$, $\forall \vec{z} \in \mathcal{U}$, then there exists a critical manifold given by \cref{eq:RPCM}. The slow flow, found by substituting \cref{eq:RPCM} into \cref{eq:RP2},
\begin{equation}\label{eq:RPSF}
\dot{\vec{z}}=\vec{W}\left(\left(c(\vec{z}),s(\vec{z})\right)^\intercal, 0,0,\vec{z}\right)
\end{equation}
is independent of the regularisation $\Psi$. 

It is evident from the limit $\varepsilon\to 0$ of \cref{eq:scalings} that this slow flow corresponds to a notion of a sliding vector field along $\Sigma$. This definition differs from that in \cite{Antali2017}, where the Filippov convention was used. We show in \cref{sec:lin} that these definitions are equivalent, provided $\vec{V}$ depends linearly on $\vec{e}$.
Nevertheless, for a general $\vec{V}$, there is nothing to suggest that there is a unique solution \cref{eq:RPCM} for any given $\vec{z}$, and therefore no guarantee of a unique sliding flow (see \cref{sec:nonunique} below). 

We can also find the layer problem of \cref{eq:vfreg_ode_sf2},
\begin{equation}
\begin{split}\label{eq:vfreg_ode_lp}
( {x}_2',{y}_2')^\intercal &= \vec{U}(\vec{e}_\Psi( x_2, y_2;1), 0,0,\vec{z}),\\
{\vec{z}}'&=0,
\end{split}
\end{equation}
where  $(\cdot)'=\frac{\mathrm{d}}{\mathrm{d}\tau}$ and the fast time $\tau$ is given by $\mathrm{d}\tau=\frac{1}{\varepsilon}\mathrm{d}t$. From Fenichel's theorem \cite{fenichel1979geometric}, the flow of \cref{eq:vfreg_ode_sf2} can be approximated by combining the flow of \cref{eq:RP} with the flow of \cref{eq:vfreg_ode_lp} for $0<\varepsilon\ll1$ where the critical manifold \cref{eq:RPCM} is normally hyperbolic with respect to \cref{eq:vfreg_ode_lp}.
We discuss the layer problem and its implications in \cref{sec:geneps1}.

\subsection{Summary}
In \cref{sec:nonsmoothdynamics}, using polar coordinates, we have shown how trajectories of the nonsmooth system \cref{eq:vf} approach or leave $\Sigma$. In \cref{sec:regularisation}, we have shown that the slow flows of  \cref{eq:vfreg_ode} limit onto $\Sigma$ as $\varepsilon\to0$. However, we have not yet been able to connect the solutions of the nonsmooth system \cref{eq:vf} in \cref{sec:nonsmoothdynamics} with those of the general regularised system \cref{eq:vfreg_ode} in \cref{sec:regularisation}. 
In the next section, using blowup, we will describe these dynamical systems in two charts, which can be connected by changes of coordinates in the regions where the charts overlap.

\section{Blowup for the general regularised system}\label{sec:gen}
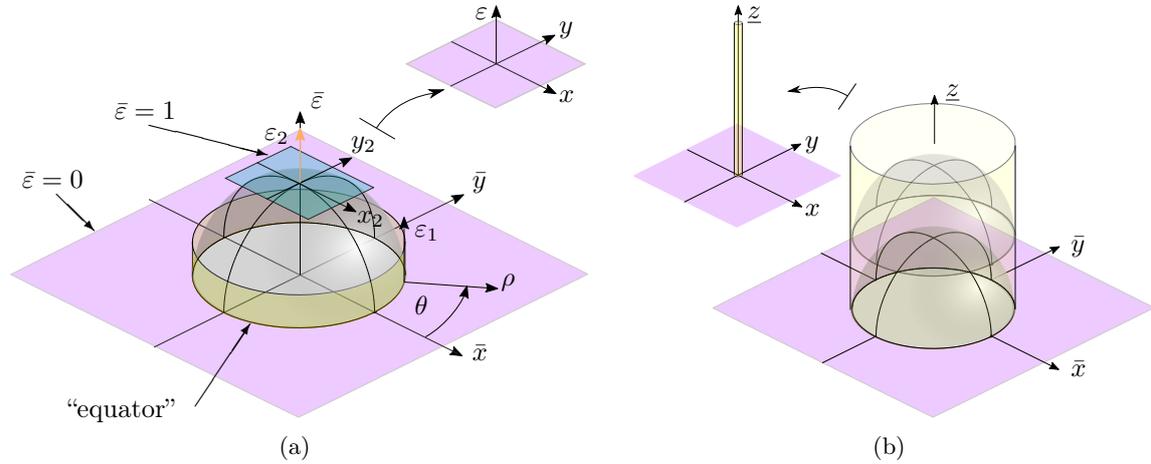
\begin{figure}[htbp]
\sbox\twosubbox{%
  \resizebox{\dimexpr.99\textwidth}{!}{%
   \begin{overpic}[height=3cm]{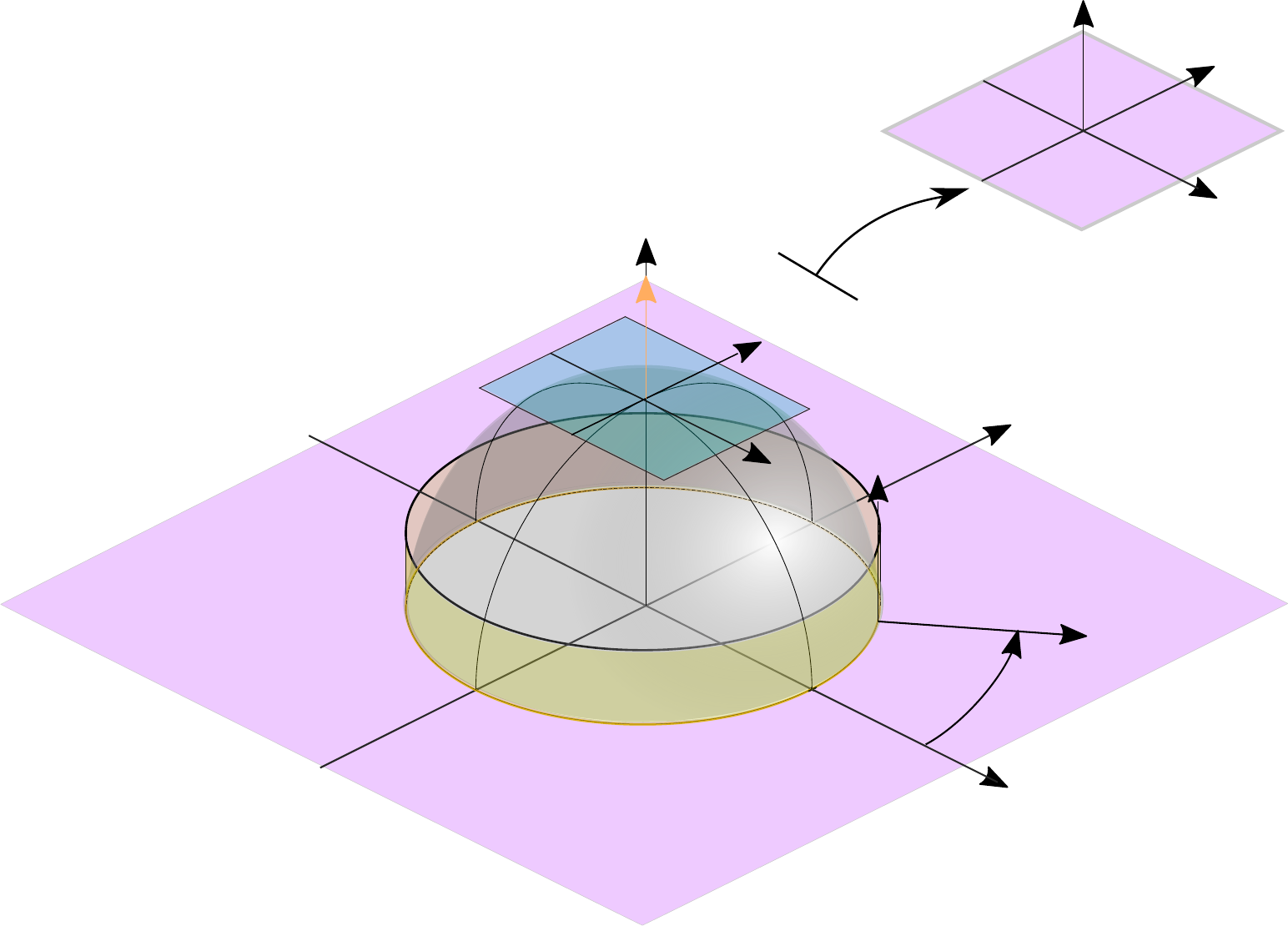}
   \end{overpic}\quad
   \begin{overpic}[height=3cm]{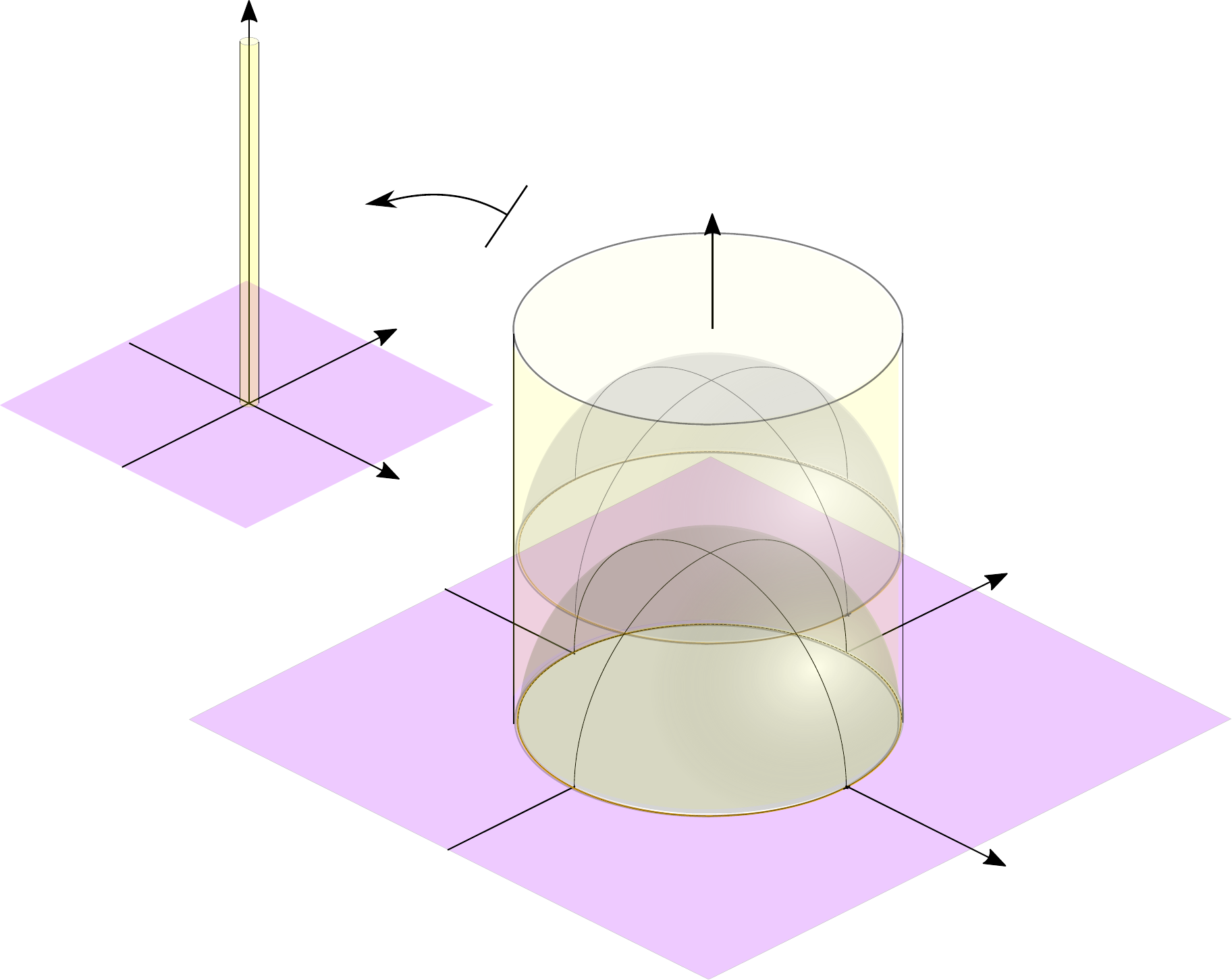}
   \end{overpic}\quad
  }%
}
\setlength{\twosubht}{\ht\twosubbox}
\centering
\hfil
\subcaptionbox{\label{fig:onesphere}}{
	\centering
	\begin{overpic}[height=\twosubht]{Figs/new_blowup_sketch_transform2_coords}
    	\put(80,10){$\bar{x}$}
        \put(80,40){$\bar{y}$}
        \put(52,54){$\bar{\varepsilon}$}
        \put(8,39){\vector(1,-1){6.5}}
        \put(2,40){$\bar{\varepsilon}=0$}
        \put(24,51){\vector(4,-1){20.5}}
        \put(18,52){$\bar{\varepsilon}=1$}
        \put(85,23){$\rho$}
        \put(60,34){{$x_2$}}
        \put(59,47){$y_2$}
        \put(44,48.5){$\varepsilon_2$}
        \put(70,31.5){$\varepsilon_1$}
        \put(70,18){$\theta$}
        \put(32,2){\vector(2,3){9.5}}
        \put(10,0){``equator''}
    	\put(95,55){${x}$}
        \put(95,67){${y}$}
        \put(80.5,69.5){${\varepsilon}$}
	\end{overpic}
}
\hfil
\subcaptionbox{\label{fig:manyspheres}}{
	\centering
	\begin{overpic}[height=\twosubht]{Figs/new_blowup_sketch_4d_urm5}
		\put(84,8){$\bar{x}$}
        \put(84,32){$\bar{y}$}
        \put(60,62){$\vec{z}$}
        \put(33,39){$x$}
        \put(33,52){$y$}
        \put(22,77){$\vec{z}$}
	\end{overpic}
}
\caption{Sketches of the blowup geometry, given by \cref{eq:blowup}: (a) blowup of a point $\vec{z}$ on $x=y=\varepsilon=0$ to the unit hemisphere (note the geometric interpretations of the scaling chart $\bar{\varepsilon}=1$, the plane $\bar{\varepsilon}=0$ and the coordinates in the charts $\kappa_1, \kappa_2$). (b) how blowup extends the dimension of \textit{each} point $\vec{z}$ in the implicit set $x=y=\varepsilon=0$ (here we collapse the $m$-dimensional $\vec{z}$ onto a single axis).\label{fig:insert_sphere}}
\end{figure}

In this section we study a general regularised system of the form \cref{eq:vfreg_ode} using blowup.
This analysis acts as a demonstration of our approach in \cref{sec:lin} where $\vec{V}$ depends linearly on $\vec{e}$.
First, we consider the system 
\begin{equation}
\begin{split}\label{eq:vfreg_extend}
	({x}',{y}',\vec{z}')^\intercal &=\varepsilon \vec{V}(\vec{e}_\Psi(x,y;\varepsilon),x,y,\vec{z}),\\
	{\varepsilon}'&=0,
\end{split}
\end{equation}
found by rescaling time in \cref{eq:vfreg_ode} and treating the parameter $\varepsilon$ as a variable.
In \cref{eq:vfreg_extend}, the set $\varepsilon=0$ is a hyperplane of equilibria and the intersection of
$\varepsilon=0$ with $\Sigma$ is singular due to the nonsmoothness of $\vec{e}_\Psi$ there.
But we can gain smoothness through the use of the blowup transformation
\begin{equation}\label{eq:blowup}
	(\rho,(\bar{x},\bar{y},\bar{\varepsilon}),\vec{z})\rightarrow(x,y,\varepsilon,\vec{z})
\end{equation}
defined by
\begin{equation}\label{eq:blowup_def}
	(x,y,\varepsilon)=\rho(\bar{x},\bar{y},\bar{\varepsilon}),\quad (\rho,(\bar{x},\bar{y},\bar{\varepsilon}))\in \lbrack 0, \infty)\times S_+^2 ,
\end{equation}
where
\begin{equation}\label{eq:sphere}
	S_+^2:=\{(\bar{x},\bar{y},\bar{\varepsilon})\,|\, \bar{x}^2+\bar{y}^2+\bar{\varepsilon}^2=1,\bar{\varepsilon}\geq0\}
\end{equation}
is the unit hemisphere. 
Informally, we are inserting a hemisphere at each point on $x=y=\varepsilon=0$ in order to make the vector field well behaved along the discontinuity set $\Sigma$ in the limit $\varepsilon\to0$ (see \cref{fig:insert_sphere}).
With  \cref{eq:blowup},  $\bar{\varepsilon}$ is then a common factor of the transformed vector field and so can be divided out. We study the resulting desingularised system in the sequel.

Although blowup allows us to gain smoothness in the vector field, we cannot easily study the dynamics everywhere simultaneously. Instead we use \textit{charts} to study multiple separate systems that are simpler but valid  only in certain regions. 
Here, we adopt an atypical approach to charts. Whilst we study the {\it scaling} chart (the directional chart found by setting $\bar{\varepsilon}=1$) in the usual way, the trigonometric terms associated with $\vec{e}$ mean that it is natural to use a single {\it entry/exit} chart\footnote{Hence we do not consider \textit{multiple directional} entry/exit charts (found by setting $\bar{x}=\pm1$ and $\bar{y}=\pm1$ here), despite Krupa \& Szmolyan's observation  \cite{krupa2001extending} that it is often ``almost mandatory'' to do so.} (found by setting $\bar{\varepsilon}=0$).
In particular, near the equator $\{S_+^2|\bar{\varepsilon}=0\}$ we re-parameterise the blowup using polar-like coordinates $(\rho,\bar{\varepsilon},\theta)$, to write 
\begin{equation}\label{eq:equator_parameterisation}
	\bar{x}=\sqrt{1-{\bar{\varepsilon}}^2}\cos{\theta},\quad \bar{y}=\sqrt{1-{\bar{\varepsilon}}^2}\sin{\theta},
\end{equation}
and so obtain the local version of \cref{eq:blowup}, the {\it entry/exit chart}
\begin{equation}\label{eq:entrychart}
	\kappa_1:\quad 	x=\rho \sqrt{1-\varepsilon_1^2}\cos{\theta},\quad y=\rho\sqrt{1-\varepsilon_1^2}\sin{\theta},\quad\varepsilon=\rho\varepsilon_1,
\end{equation}
with the chart-specific coordinates ($\rho,\varepsilon_1,\theta$).
Note that setting $\varepsilon_1=0$ in \cref{eq:entrychart} gives \cref{eq:polar_coordinates}.

The \textit{scaling chart} is found by setting $\bar{\varepsilon}=1$ in \cref{eq:blowup_def}
\begin{equation}\label{eq:scalingchart}
	\kappa_2:\quad x=\rho_2 x_2,\quad y=\rho_2 y_2, \quad \varepsilon=\rho_2,
\end{equation}
with the chart specific-coordinates ($x_2,y_2,\rho_2$).
Note that \cref{eq:scalingchart} is equivalent to \cref{eq:scalings}.

The change of coordinates between $\kappa_2$ and $\kappa_1$ is given by 
\begin{equation}\label{eq:coordinatechange}
	\kappa_{12}:(x_2,y_2,\rho_2)\mapsto
	\left\{\begin{aligned}
	\rho&=\frac{\rho_2}{\sqrt{1+x_2^2+y_2^2}}\\
	\varepsilon_1&=\frac{1}{\sqrt{1+x_2^2+y_2^2}}\\
	\theta&=\arctan\left(\frac{y_2}{x_2}\right)
	\end{aligned}\right. .
\end{equation}
These charts and coordinates are shown in 
\cref{fig:insert_sphere}.

Now we shall demonstrate that the dynamics in the plane $\varepsilon_1=0$ of the chart $\kappa_1$ correspond precisely to the dynamics studied in \cref{sec:nonsmoothdynamics}, and that the dynamics in the scaling chart $\kappa_2$ correspond precisely to the dynamics of the scaled regularised system in \cref{sec:regularisation}.

\subsection{Entry chart\texorpdfstring{ $\kappa_1$}{}: dynamics near the discontinuity set}\label{sec:geneps0}
Let us study the dynamics in the chart $\kappa_1$. Re-writing \cref{eq:vfreg_extend} using \cref{eq:entrychart} and transforming time, we divide the right hand side by the common factor $\varepsilon_1$, to find
\begin{equation}\label{eq:polarising_desingularising}
	\begin{split}
\frac{\mathrm{d}}{\mathrm{d}\mathcal{T}}\begin{pmatrix}
			{\rho}\\{\theta}
		\end{pmatrix}&=
		\begin{pmatrix}
			\rho\xi&0\\
			0&\xi^{-1}
		\end{pmatrix}
		\mat{R}^\intercal(\theta)\vec{U}\left(\tilde{\vec{e}}_\Psi(\theta,\varepsilon_1),\rho \xi\cos\theta,\rho\xi\sin\theta,\vec{z}\right),\\
		\frac{\mathrm{d}}{\mathrm{d}\mathcal{T}}{\vec{z}}&=\rho \vec{W}\left(\tilde{\vec{e}}_\Psi(\theta,\varepsilon_1),\rho \xi\cos\theta,\rho\xi\sin\theta,\vec{z}\right),\\
		\frac{\mathrm{d}}{\mathrm{d}\mathcal{T}}{\varepsilon}_1&=-\varepsilon_1\xi(\cos\theta,\sin\theta)\vec{U}\left(\tilde{\vec{e}}_\Psi(\theta,\varepsilon_1),\rho \xi\cos\theta,\rho\xi\sin\theta,\vec{z}\right),
	\end{split}
\end{equation}
where
$\mathrm{d}\mathcal{T}=\frac{\varepsilon_1}{\varepsilon}\mathrm{d}t=\frac{1}{\rho}\mathrm{d}t$,
\begin{equation}
\xi=(1-\varepsilon_1)^{\frac{1}{2}},
\label{eq:xi}
\end{equation} 
$\vec{U}$ is given in \cref{eq:gen_odes1}, $\vec{W}$ in \cref{eq:gen_odes2}, 
\begin{equation}
	\begin{split}
		\tilde{\vec{e}}_\Psi(\theta,\varepsilon_1)&:=\vec{e}_\Psi\left(\rho\xi\cos\theta,\rho\xi\sin\theta;\rho\varepsilon_1\right)\\
		&:=\left({\frac{{1-\varepsilon_1^2}}{{1-\varepsilon_1^2\left(1-\Psi_1\left(\frac{\varepsilon_1^2}{1-\varepsilon_1^2}\right)\right)}}}\right)^\frac{1}{2}(\cos\theta,\sin\theta)^\intercal
	\end{split}
\end{equation}
and $\Psi_1$ is defined as in \cref{R5} of \cref{def:Psi} in \cref{sec:regularisation}.

The singularity at $\varepsilon_1=1$  does not concern us since we will only work with \cref{eq:polarising_desingularising} near $\varepsilon_1=0$. Note that the plane $\varepsilon_1=0$ is an invariant manifold
and that the dynamics within the plane are identical to \cref{eq:polar_gen_odes_desing}.
Accordingly, we shall follow the analysis in \cref{sec:nonsmoothdynamics} closely.

\begin{theorem}\label{thm:geneps0}
	Consider \cref{eq:polarising_desingularising}, then
	\begin{enumerate}[\rm{(\alph*)}]
		\item{	If $\theta=\theta_0$ is a solution to 
				\begin{equation}\label{eq:Theta}
					\Theta(0,\theta,\vec{z}):=(-\sin{\theta},\cos{\theta})\vec{U}((\cos{\theta},\sin{\theta})^\intercal,0,0,\vec{z})=0
				\end{equation}
				for a given $\vec{z}=\vec{z}_0$, then there exists an equilibrium of \cref{eq:vfreg_extend} along the equator of the blown up hemisphere at 
				\begin{equation}
					Q_0=(0,\theta_0,\vec{z}_0,0).
				\end{equation}\label{thm:geneps0a}}
		\item{	Consider such an equilibrium point $Q_0$: if there exists an eigenvalue $\lambda_\rho(\theta_0,\vec{z}_0)$
				of the linearisation of \cref{eq:polarising_desingularising},
				corresponding to radial attractiveness as in \cref{sec:nonsmoothdynamics}, then there exists a 
				corresponding eigenvalue  $-\lambda_\rho(\theta_0,\vec{z_0})$
				with an eigenvector solely in the $\varepsilon_1$ direction.
				Hence, if there is a stable manifold of $Q_0$ within $\varepsilon_1=0$ extending into $\rho>0$, 
				then there is a corresponding unstable manifold of $Q_0$ within $\rho=0$ extending into $\varepsilon_1>0$, and vice versa.\label{thm:geneps0b}}
		\item{	If no solutions for $\theta$ exist to \cref{eq:Theta} for a given $\vec{z}=\vec{z}_0$, 
				then there are no equilibria along the equator of the blown-up hemisphere inserted at that 
				particular $\vec{z}$ (the case of no \textit{limit directions} in \cite{Antali2017}). There is therefore a limit cycle on the equator which which has a stable manifold in $\varepsilon_1=0$ and an unstable manifold in $\rho=0$ if $\int_0^{2\pi} \lambda_\rho(\theta,\vec{z})/|\Theta(0,\theta,\vec{z})|\mathrm{d}\theta<0$, and vice versa if $\int_0^{2\pi} \lambda_\rho(\theta,\vec{z})/|\Theta(0,\theta,\vec{z})|\mathrm{d}\theta>0$.  \label{thm:geneps0c}}
		\item{	On the sphere, where $\rho=0$, $\vec{z}$ remains constant for the flow of \cref{eq:polarising_desingularising}.\label{thm:geneps0d}} 
	\end{enumerate}
\end{theorem}

\begin{figure}[htbp]
  \centering
   \begin{overpic}[width=0.4\textwidth]{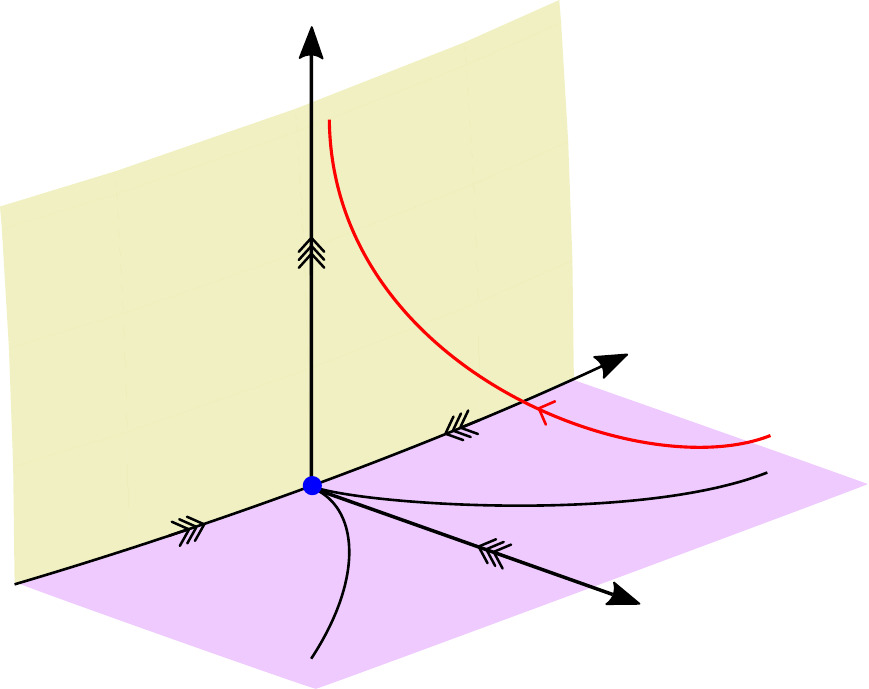}
   		\put(29,26){$Q_0$}
    	\put(75,5){$\rho$}
	    \put(27,72){$\varepsilon_1$}
	    \put(70,40){$\theta$}
	    \put(17,21.5){$\lambda_\theta$}
	    \put(48,34){$\lambda_\theta$}
	    \put(50,10){$\lambda_\rho$}
	    \put(22,50){$-\lambda_\rho$}
   \end{overpic}
  \caption{Sketch demonstrating how orbits travel up the side of the sphere. The yellow surface corresponds to the hemisphere of the blowup. The red orbit is for $0<\varepsilon\ll1$. }\label{fig:side_of_sphere}
\end{figure}

\begin{proof}
	Following \cref{sec:nonsmoothdynamics} closely,
	\begin{enumerate}[\rm{(\alph*)}]
		\item{	We have already noted that $\varepsilon_1=0$ is an invariant manifold of \cref{eq:polarising_desingularising} 
				and that dynamics within it are given by \cref{eq:polar_gen_odes_desing}. Therefore, if there exists an equilibrium $P_0$ of
				\cref{eq:polar_gen_odes_desing} then $Q_0$ is an equilibrium of \cref{eq:polarising_desingularising}.}
		\item{	The Jacobian of \cref{eq:polarising_desingularising} at $Q_0$, is given by
				\begin{equation}
					\mat{K}_0(\theta_0,\vec{z}_0):=
					\begin{pmatrix}
						\mat{J}_0(\theta_0,\vec{z}_0)&\mat{0}_{ m \times1}\\
						\mat{0}_{1\times m }&-\lambda_\rho(\theta_0,\vec{z}_0)
					\end{pmatrix}.
				\end{equation}
				It is then straightforward to show that there are now three non-zero eigenvalues: 
				$\lambda_\rho$, $\lambda_\theta$ and $-\lambda_\rho$. 
				The eigenvectors corresponding to $\lambda_\rho$ and $\lambda_\theta$ are as before, with 0 in the $\varepsilon_1$ component. The eigenvector corresponding to
				$-\lambda_\rho$ points solely in the $\varepsilon_1$ direction and hence represents travelling up 
				or down the sphere (\cref{fig:side_of_sphere}).
				Then by the stable manifold theorem \cite{perko}, there are corresponding stable and unstable manifolds of $Q_0$ tangent to the eigenvectors associated with either $\lambda_\rho$ and $-\lambda_\rho$ or vice versa}.
		\item{	If no solutions exist to \cref{eq:Theta} for a fixed $\vec{z}=\vec{z}_0$, then along the invariant manifold of the equator ($\rho=\varepsilon_1=0$), 
				the $\theta$ component of \cref{eq:polarising_desingularising} does not change sign and there exists a limit cycle along the equator. In the plane $\varepsilon_1=0$, if we take a Poincar\'{e} section at some constant $\theta$, we can find that the local stability of this limit cycle is given by $\int_0^{2\pi} \lambda_\rho(\theta,\vec{z})/|\Theta(0,\theta,\vec{z})|\mathrm{d}\theta$. As with the stability of equilibria along the equator, the stability of the closed orbit is opposite within $\rho=0$.}
		\item{	On the sphere $\rho=0$, we find $\dot{\vec{z}}=\vec{0}$ from \cref{eq:polarising_desingularising}. 
				Hence the dynamics on the sphere of trajectories that reach the equator at $Q_0$ are effectively parameterised by $\vec{z}=\vec{z}_0$ 
				in \cref{eq:polarising_desingularising}.}
	\end{enumerate}
\end{proof}

From \cref{thm:geneps0}, we note that if $\lambda_\rho(\theta_0,\vec{z}_0)<0$, then there exists an orbit of \cref{eq:polarising_desingularising} that approaches
the equator along $\theta=\theta_0$ (tangent to a stable manifold of the equilibrium $Q_0$) and then travels along the corresponding unstable manifold of $Q_0$ up the sphere.
If both $\lambda_\rho(\theta_0,\vec{z}_0)<0$ and $\lambda_\theta(\theta_0,\vec{z}_0)<0$ then there exists an open set of orbits of the nonsmooth system that reach $Q_0$ and travel up the sphere (see \cref{fig:side_of_sphere}). 
Manifolds that appear in the dynamics in the chart $\kappa_1$ that extend into $\varepsilon_1>0$ can be recovered in the scaling chart $\kappa_2$ using \cref{eq:coordinatechange}, where the system is standard slow-fast. 

\begin{figure}[htbp]
	\centering
	\begin{subfigure}{0.3\textwidth}
        \begin{overpic}[width=0.9\textwidth]{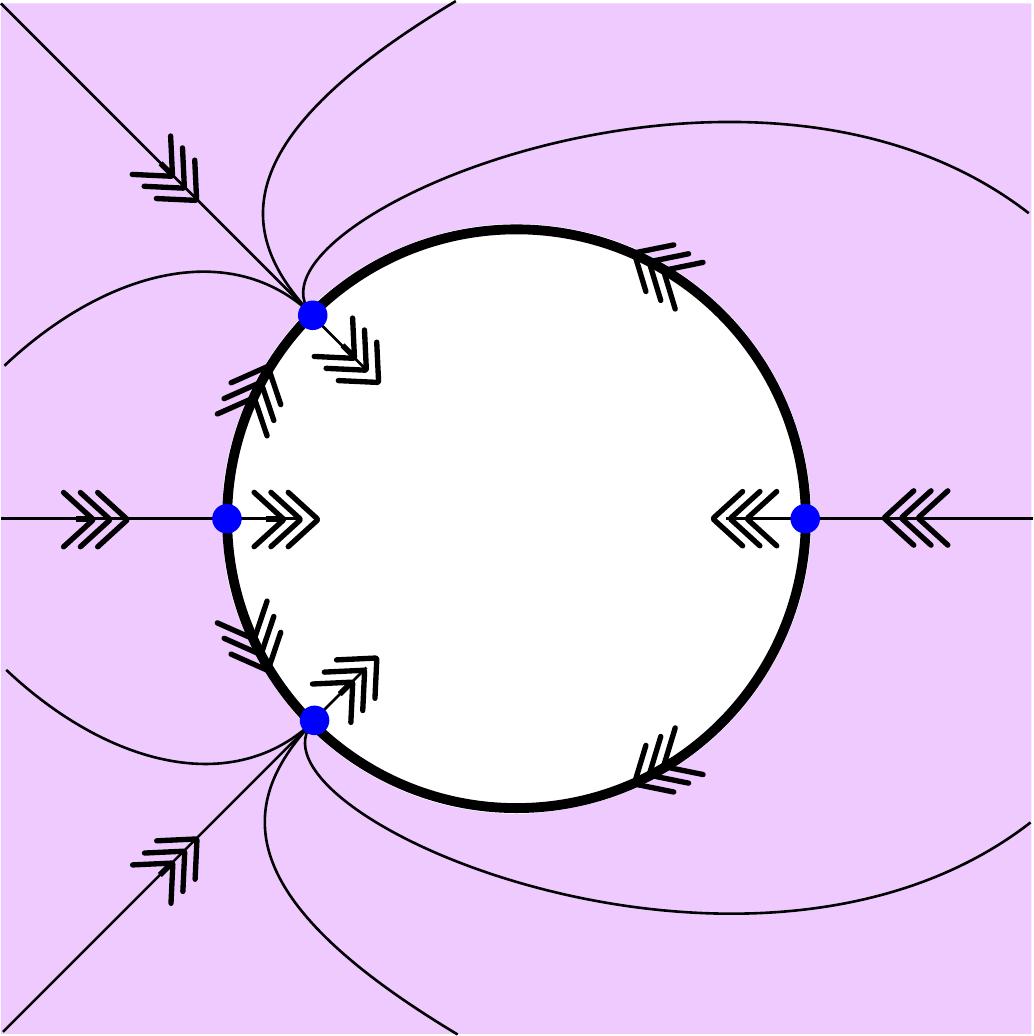}
    	    \put(68,54.5){$Q_1$}
	        \put(25.5,54.5){$Q_3$}
	        \put(35,68){$Q_2$}
	        \put(35,28){$Q_4$}
        \end{overpic}
        \caption{}\label{fig:eps0_node}
	\end{subfigure}
	\hfil
	\begin{subfigure}{0.3\textwidth}
        \begin{overpic}[width=0.9\textwidth]{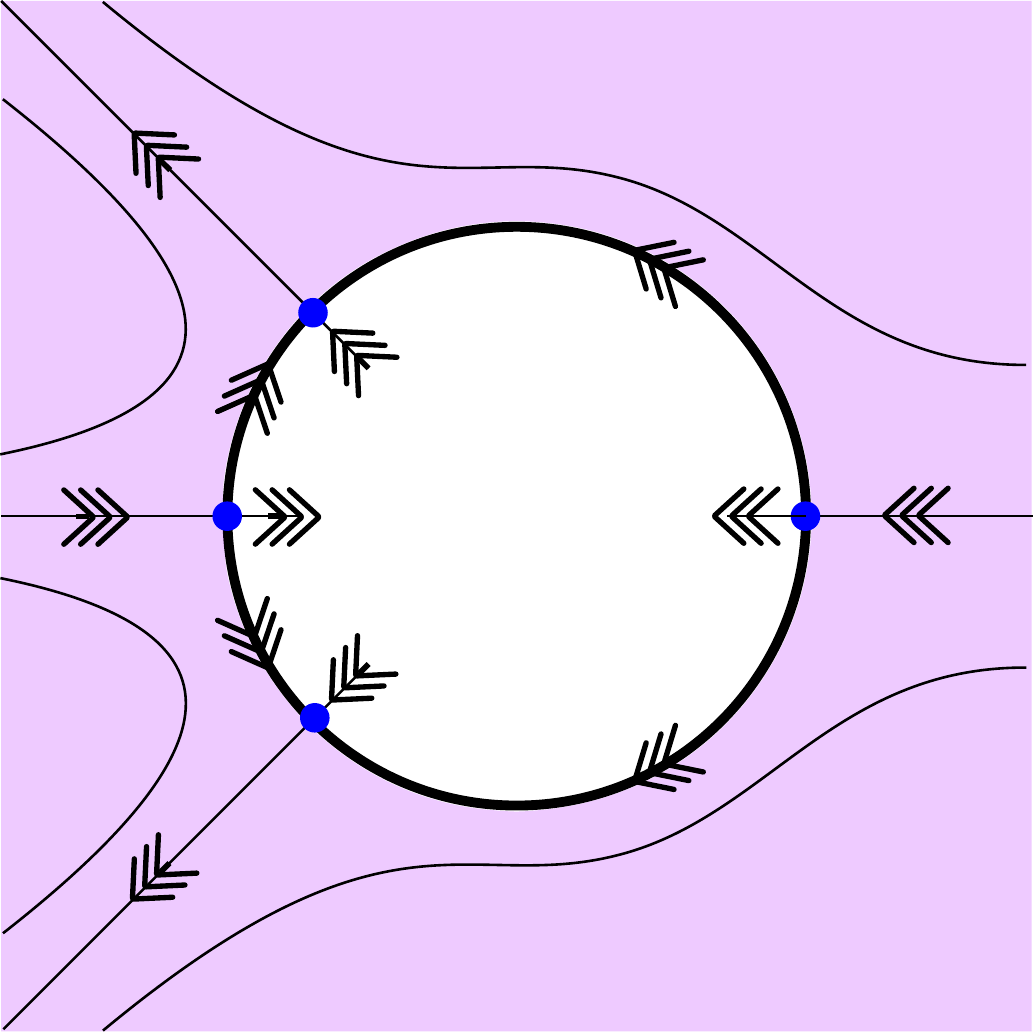}
    	    \put(68,54.5){$Q_1$}
	        \put(25.5,54.5){$Q_3$}
	        \put(35,68){$Q_2$}
	        \put(35,28){$Q_4$}
        \end{overpic}
		\caption{}\label{fig:eps0_saddle}
	\end{subfigure}
	\hfil
	\begin{subfigure}{0.3\textwidth}
        \begin{overpic}[width=0.9\textwidth]{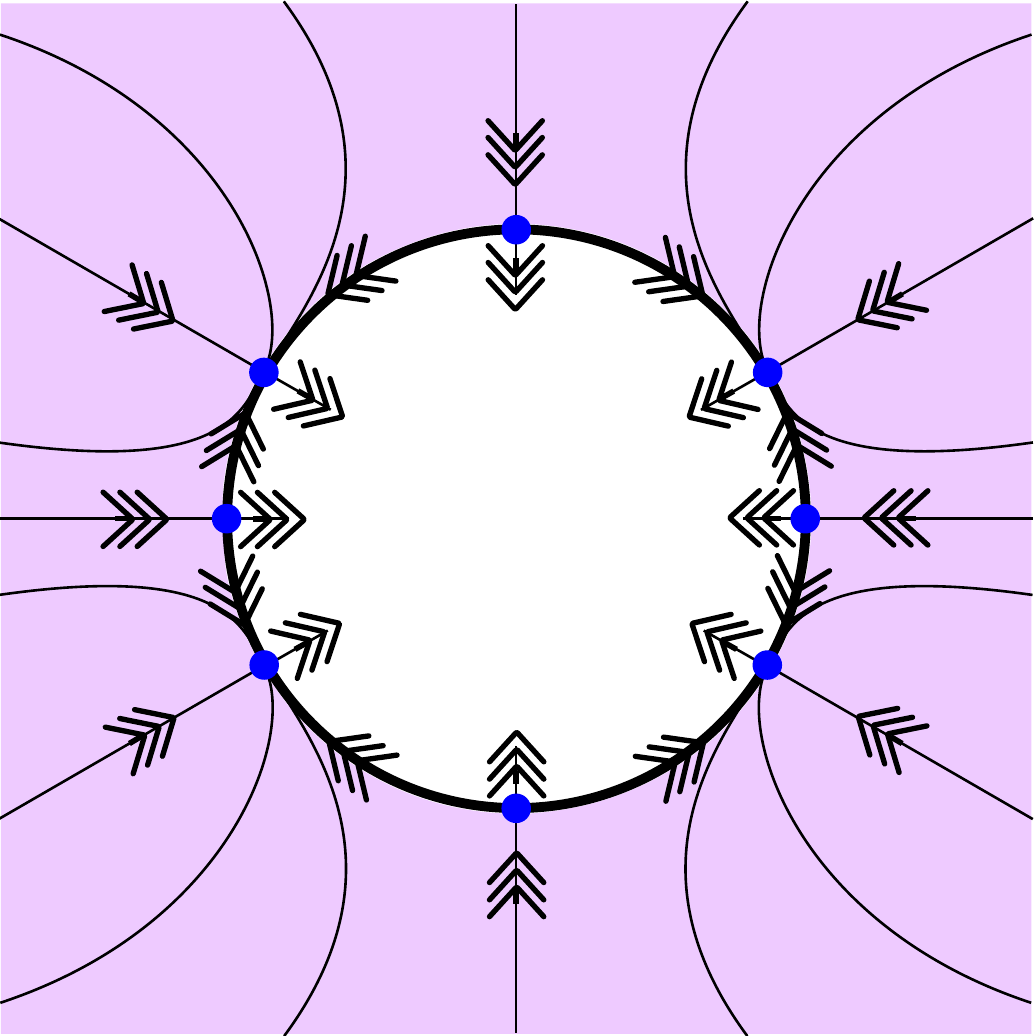}
			{\scriptsize
			\put(67.5,54){$Q_1$}
			\put(63,65.5){$Q_2$}
			\put(53,72){$Q_3$}
			\put(29.5,65){$Q_4$}
			\put(25,43){$Q_5$}
			\put(29.5,31.5){$Q_6$}
			\put(53,25.5){$Q_7$}
			\put(66.5,42.5){$Q_8$}}
        \end{overpic}
		\caption{}\label{fig:eps0_cubic}
	\end{subfigure}
	\hfil 
    \caption{Sketches of the dynamics in the chart $\kappa_1$ for $\varepsilon_1=0$, projected into the $(\bar{x},\bar{y})$ plane.
    	{(a)} all equilibria on the equator are radially attracting ($\lambda_r|_{Q_i}<0,\forall i$),
    	whilst $Q_1$ and $Q_3$ are angularly repelling and $Q_2$ and $Q_4$ are angularly attracting
    	($\lambda_\theta|_{Q_i}>0,i\in\{1,3\}$ and $\lambda_\theta|_{Q_i}<0,i\in\{2,4\}$ respectively).
    	{(b)} $Q_1$ and $Q_3$ are radially attracting and angularly repelling ($\lambda_\theta|_{Q_i}>0,\lambda_r|_{Q_i}<0,i\in\{1,3\}$)
    	and $Q_2$ and $Q_4$ are radially repelling and angularly attracting ($\lambda_\theta|_{Q_i}<0,\lambda_r|_{Q_i}>0,i\in\{2,4\}$).
    	{(c)} each $Q_i$ is radially attracting ($\lambda_\rho|_{Q_i}<0$ $\forall i$), and $Q_i$ are angularly attracting
    	for $i$ even ($\lambda_\theta|_{Q_i}<0,i\in\{2,4,6,8\}$) and angularly repelling for $i$ odd ($\lambda_\theta|_{Q_i}<0,i\in\{1,3,5,7\}$).}
    \label{fig:eps0}
\end{figure}
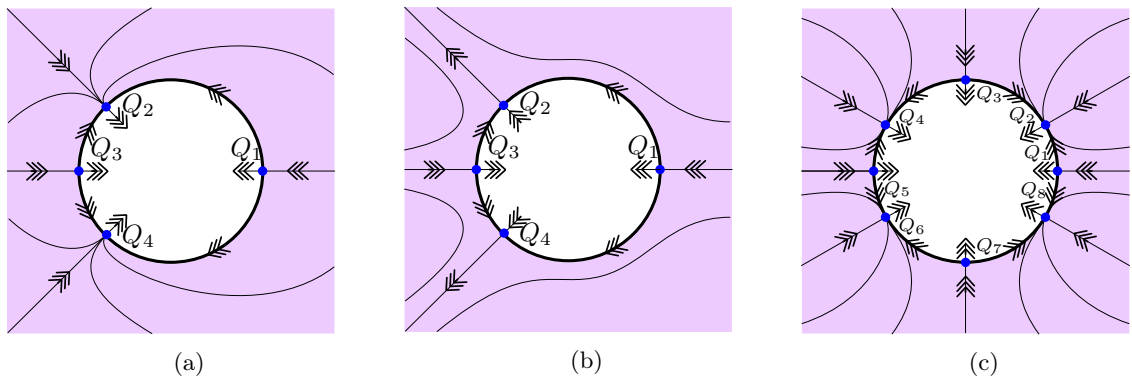

In \cref{fig:eps0}, we see sketches of dynamics in the chart $\kappa_1$ for $\varepsilon_1=0$, projected into the $(\bar{x},\bar{y})$ plane.
The sphere and manifolds extending from equilibria along its equator are also projected down into the same plane, ignoring the $\vec{z}$ dynamics.

Odd numbers of equilibria along the equator occur at bifurcations, when
\begin{equation}\label{eq:eqbif}
  \Theta(0,\theta^*,\vec{z})=\pdiff{}{\theta}\Theta(0,\theta^*,\vec{z})=0.  
\end{equation}
Where there are no solutions to \cref{eq:Theta}, and hence no equilibria along the equator (no limit-directions),
the equator is a closed orbit of the system.
Therefore, from index theory \cite{perko}, there must exist at least one critical set in the region enclosed by the equator.
This region is described by the scaling chart $\kappa_2$, which we will now discuss.

\subsection{Scaling chart\texorpdfstring{ $\kappa_2$}{}: dynamics on the discontinuity set}\label{sec:geneps1}
In the previous section, using the entry chart $\kappa_1$, we studied how orbits can reach or leave $\Sigma$. To study the dynamics \textit{along} $\Sigma$ we use the scaling chart $\kappa_2$.
Substituting \cref{eq:scalingchart} into \cref{eq:vfreg_extend}, we find
\begin{subequations}\label{eq:geneps1}
	\begin{align}
			(x_2',y_2')^\intercal&=\vec{U}\left(\vec{e}_\Psi(x_2,y_2;1),\varepsilon x_2,\varepsilon y_2,\vec{z}\right),\label{eq:geneps1.1} \\
			\vec{z}'&=\varepsilon \vec{W}\left(\vec{e}_\Psi(x_2,y_2;1),\varepsilon x_2,\varepsilon y_2,\vec{z}\right), \label{eq:geneps1.2}\\
			\varepsilon'&=0.\label{eq:geneps1.3}
	\end{align}
\end{subequations}
Equations \cref{eq:geneps1.1} and \cref{eq:geneps1.2} are equivalent to \cref{eq:vfreg_ode_sf2} but written with respect to a fast time.

In \cref{sec:regularisation}, we identified that if there is a solution 
$\vec{e}_\Psi(x_2,y_2;1)=\left(c(\vec{z}),s(\vec{z})\right)^\intercal$
such that $\vec{U}(\left(c(\vec{z}),s(\vec{z})\right)^\intercal,0,0,\vec{z})=0$,
then there exists a critical manifold in the scaling chart $\kappa_2$ given by \cref{eq:RPCM}, and a corresponding slow flow \cref{eq:RPSF}.
However, we have yet to study the layer problem \cref{eq:vfreg_ode_lp} and its implications.

First, it will be useful to calculate the derivative $D_{(x_2,y_2)}\vec{e}_\Psi(x_2,y_2;1)$ to find
\begin{equation}
		D_{(x_2,y_2)}\vec{e}_\Psi(x_2,y_2;1)
		=\frac{1}{\sqrt{\zeta^2+\Psi(\zeta^2)}}\mat{I}-\frac{1+\Psi'(\zeta^2)}{\sqrt{\zeta^2+\Psi(\zeta^2)}^3}
		\begin{pmatrix}
			x_2\\y_2
		\end{pmatrix}
		\begin{pmatrix}
			x_2&y_2
		\end{pmatrix}, 
\end{equation}
where $\zeta^2:=x_2^2+y_2^2$.

\begin{remark}\label{rem:jaco_properties}
We note certain properties of the matrix $D_{(x_2,y_2)}\vec{e}_\Psi(x_2,y_2;1)$. 
Firstly, by \cref{def:Psi}, both the trace and determinant are positive:
\begin{align}
	\tr\left(D_{(x_2,y_2)}\vec{e}_\Psi(x_2,y_2;1)\right)&=\frac{\zeta^2\left(1-\Psi'(\zeta^2)\right)+2\Psi(\zeta^2)}{\left(\zeta^2+\Psi\left(\zeta^2\right)\right)^{3/2}}>0
	\intertext{and}
	\det(D_{(x_2,y_2)}\vec{e}_\Psi(x_2,y_2;1))&=\frac{\Psi(\zeta^2)-\Psi'(\zeta^2)\zeta^2}{(\zeta^2+\Psi(\zeta^2))^2}>0
	\intertext{respectively. Secondly $D_{(x_2,y_2)}\vec{e}_\Psi(x_2,y_2;1)$ is positive definite, since}
	\vec{\eta}^\intercal D_{(x_2,y_2)}\vec{e}_\Psi(x_2,y_2;1)\, \vec{\eta}&=
	\frac{(\eta_1y_2-\eta_2x_2)^2 +\vec{\eta}^\intercal\vec{\eta}\Psi(\zeta^2)-(\eta_1x_2+\eta_2y_2)^2 \Psi'(\zeta^2)}{\left(\zeta^2 + \Psi(\zeta^2)\right)^{\frac{3}{2}}}>0
\end{align}
for any $\vec{\eta}=(\eta_1,\eta_2)^\intercal$.
\end{remark}

\begin{theorem}\label{thm:geneps1}
Given a system of the form \cref{eq:vfreg_extend}, the following holds. 
\begin{enumerate}[\rm{(\alph*)}]
	\item{	Suppose there exists smooth functions $ c(\vec{z}),s(\vec{z}):c(\vec{z})^2+s(\vec{z})^2<1$ for $\vec{z}\in\mathcal{U}$ where 
			\begin{equation}\label{eq:csdef}
				\vec{U}\left((c(\vec{z}),s(\vec{z}))^\intercal,0,0,\vec{z}\right)=\vec{0},
			\end{equation}
			then
			\begin{equation}\label{eq:critical_set}
				C=\left\{(x_2,y_2,\vec{z})\left| \vec{e}_\Psi(x_2,y_2;1)=(c(\vec{z}),s(\vec{z}))^\intercal,\vec{z}\in \mathcal{U}\right.\right\}
			\end{equation}
			is a critical set of the equations in the scaling chart $\kappa_2$, that is, a set of equilibria of the layer problem.\label{thm:geneps1a}}
	\item{	The attractiveness and normal hyperbolicity of the critical set \cref{eq:critical_set} are given by the eigenvalues of the Jacobian of the fast subsystem \cref{eq:geneps1.1} about $C$
			\begin{equation}\label{eq:scjaco}
				\mat{J}:=\left.D_{(c,s)}\vec{U}((c,s)^\intercal,0,0,\vec{z})D_{(x_2,y_2)}\vec{e}_\Psi(x_2,y_2;1)\right|_C.
			\end{equation}
			Then it follows that if $\det{\left(D_{(c,s)}\vec{U}((c,s)^\intercal,0,0,\vec{z})\right)}<0$ then the critical set $C$
			is a saddle with respect to the fast flow, irrespective of the regularisation function $\Psi$. On the other hand, if 
			\begin{equation}
				\det{\left(D_{(c,s)}\vec{U}((c,s)^\intercal,0,0,\vec{z})\right)}>0
			\end{equation}						
			then the critical set $C$ can be a stable or unstable node or focus, possibly depending on the regularisation.\label{thm:geneps1b}}
	\item{	When a critical set $C$ exists, the slow flow along it is given by
			\begin{equation}\label{eq:slow_flow}
				\dot{\vec{z}}=\vec{W}((c(\vec{z}),s(\vec{z}))^\intercal,0,0,\vec{z}),
			\end{equation}
			which is independent of the regularisation function $\Psi$.\label{thm:geneps1c}}
	\item{	Where a  critical manifold $C$ is compact and normally hyperpoblic, 
			it perturbs to a corresponding \textit{slow manifold} $C_\varepsilon$ in the smoothed system for $0<\varepsilon\ll1$.
			This slow manifold $C_\varepsilon$ lies $\varepsilon$-close to $C$ and the slow flow \cref{eq:slow_flow} gives the first order
			approximation to the flow along $C_\varepsilon$. \label{thm:geneps1d}}
\end{enumerate}
\end{theorem}
\begin{proof}
The proof is as follows.
\begin{enumerate}[\rm{(\alph*)}]
	\item{	Treating $\varepsilon$ as a parameter, let us study the layer problem of \cref{eq:geneps1}, found by setting $\varepsilon=0$,
			\begin{subequations}\label{eq:LP}
				\begin{align}
					(x_2',y_2')^\intercal &= \vec{U}(\vec{e}_\Psi( x_2, y_2;1), 0,0,\vec{z})\label{eq:LP1}\\
					\vec{z}'&=0\label{eq:LP2}.
				\end{align}
			\end{subequations}
			Evidently sets of the form \cref{eq:critical_set} are sets of equilibria of \cref{eq:LP}, also described by \cref{eq:RPCM}.}
	\item{	In order to describe the local stability of these critical set we study the layer problem.
			Linearising  \cref{eq:LP1} around a critical set $C$, we find the Jacobian \cref{eq:scjaco}. 
			Since the determinant of $D_{(x_2,y_2)}\vec{e}_\Psi(x_2,y_2;1)$
			is strictly  positive (\cref{rem:jaco_properties}), we can determine that 
			\begin{align*}
				\sign(\det{\mat{J}})&\equiv\sign\left(\det\left(D_{(c,s)}\vec{U}((c,s)^\intercal,0,0,\vec{z})\,D_{(x_2,y_2)}\vec{e}_\Psi(x_2,y_2;1)\right)\right)\\
				&\equiv\sign\left(\det\left(D_{(c,s)}\vec{U}((c,s)^\intercal,0,0,\vec{z})\right)\right).
			\end{align*}
			Therefore, if $\det\left(D_{(c,s)}\vec{U}((c,s)^\intercal,0,0,\vec{z})\right)<0$ then the critical set $C$ is a saddle,
			irrespective of the regularisation function. It is not so straightforward to comment on the stability on the critical set when
			$\sign\left(\det\left(D_{(c,s)}\vec{U}((c,s)^\intercal,0,0,\vec{z})\right)\right)>0$.
			In fact, in Case III of \cref{sec:cases}, we prove by example that the stability of $C$ can depend upon the regularisation function $\Psi$.}
	\item{	As in \cref{sec:regularisation}, 
	we find the slow flow along a critical set $C$  is given by \cref{eq:slow_flow}, subject to the algebraic constraint that $(x_2,y_2)\in C$, which is trivially independent of $\Psi$.}
	\item{	By Fenichel's theorem \cite{fenichel1979geometric}, where a critical manifold $C$ is compact and normally hyperbolic, it perturbs to a slow manifold
			$C_\varepsilon$ which lies Hausdorff distance $\mathcal{O}(\varepsilon)$ from $C$. In addition the slow flow along $C_\varepsilon$ is
			smoothly $\mathcal{O}(\epsilon)$ close to \cref{eq:slow_flow} from the same theorem.}
\end{enumerate}
\end{proof}

\begin{figure}[htbp]
	\centering
	\begin{subfigure}{0.25\textwidth}
	\centering
        \begin{overpic}[width=0.9\textwidth]{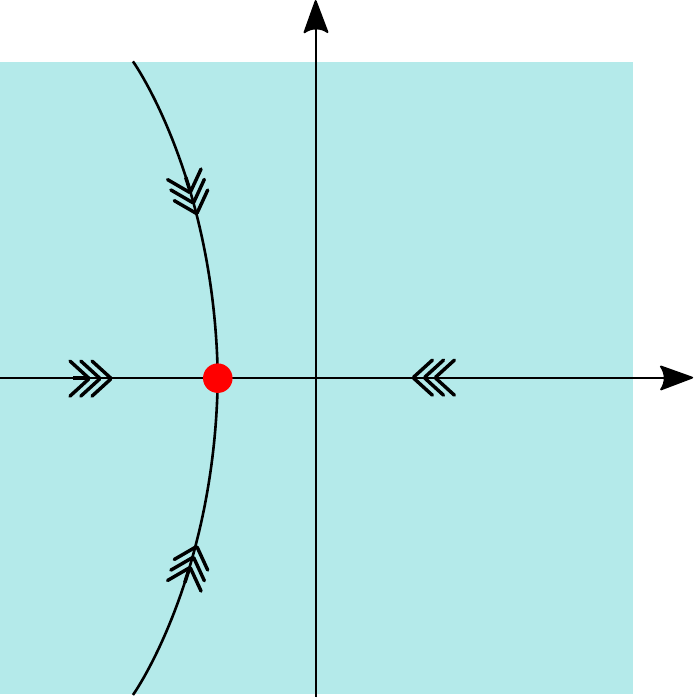}
    	    \put(93,49){$x_2$}
	        \put(48,94){$y_2$}
	        \put(23,48){$C$}
        \end{overpic}
		\\ \vspace{8pt}
        \begin{overpic}[width=0.99\textwidth]{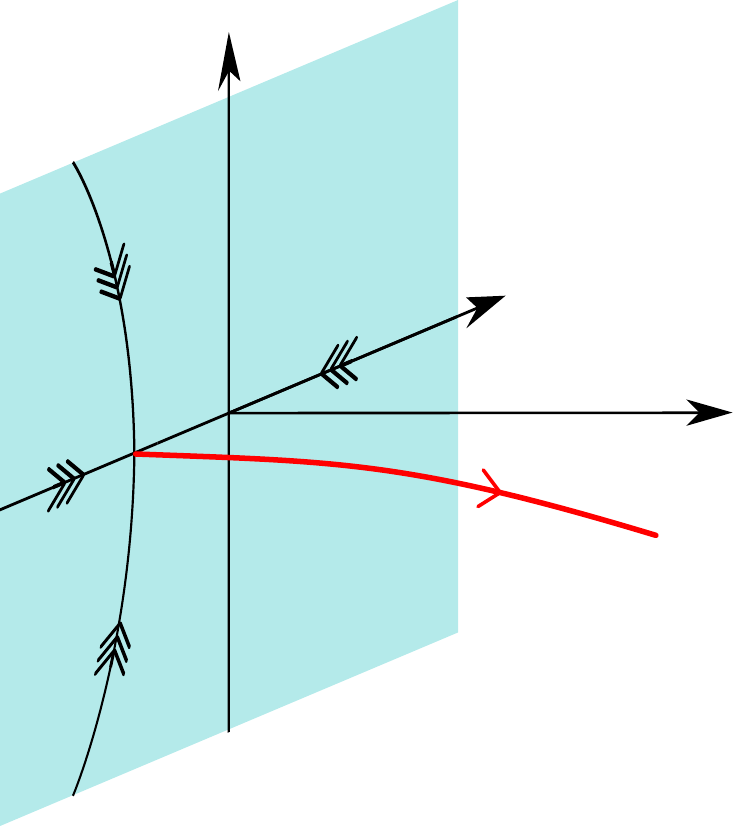}
        	\put(58,66){$x_2$}
	        \put(19,90){$y_2$}
	        \put(17,38){$C$}
	        \put(83,54){$\vec{z}$}
        \end{overpic}
        \caption{}\label{fig:eps1node_combi}
	\end{subfigure}
	\hfil
	\begin{subfigure}{0.25\textwidth}
		\centering
        \begin{overpic}[width=0.9\textwidth]{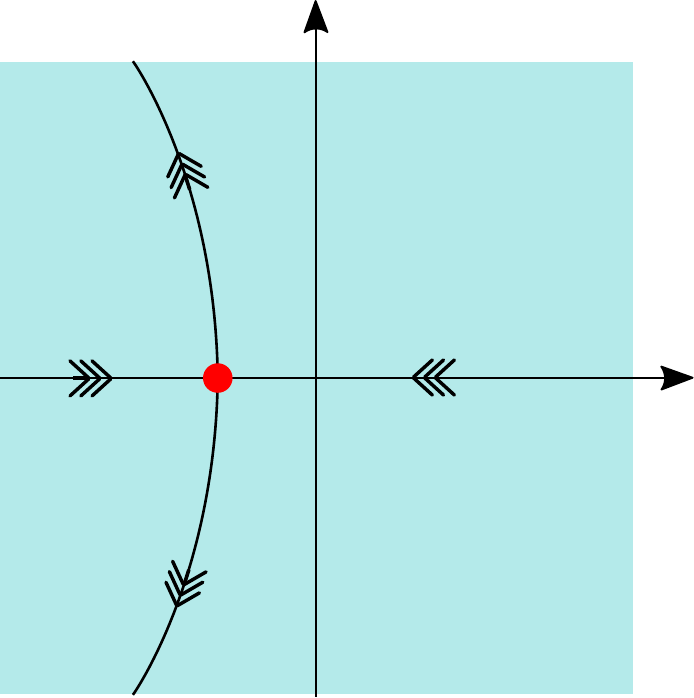}
	        \put(93,49){$x_2$}
	        \put(48,94){$y_2$}
	        \put(23,48){$C$}
        \end{overpic}
		\\ \vspace{8pt}
        \begin{overpic}[width=0.99\textwidth]{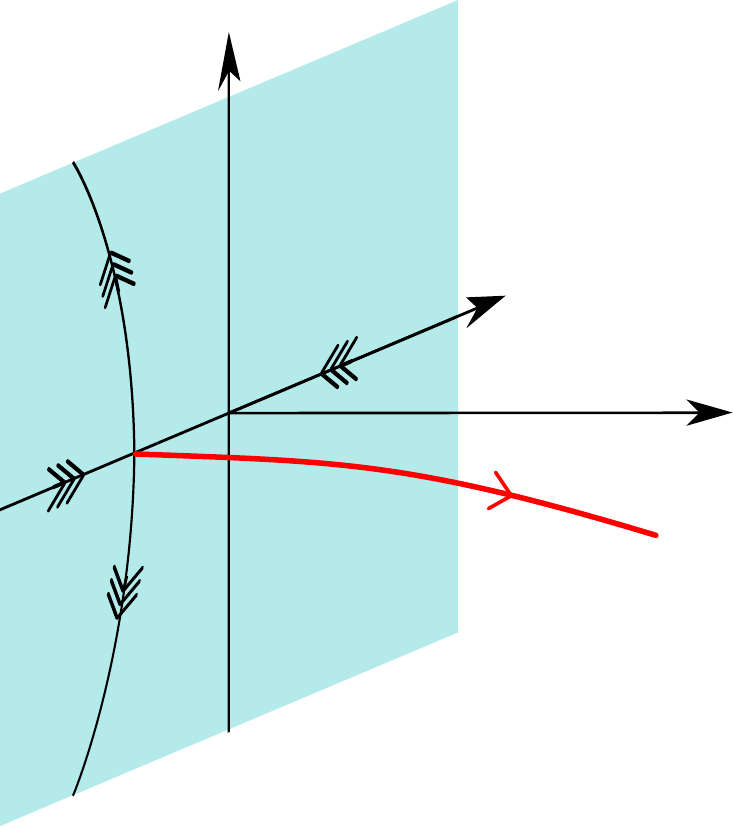}
        	\put(58,66){$x_2$}
	        \put(19,90){$y_2$}
	        \put(17,38){$C$}
	        \put(83,54){$\vec{z}$}
        \end{overpic}
		\caption{}\label{fig:eps1saddle_combi}
	\end{subfigure}
	\hfil
	\begin{subfigure}{0.25\textwidth}        
		\centering
		\begin{overpic}[width=0.9\textwidth]{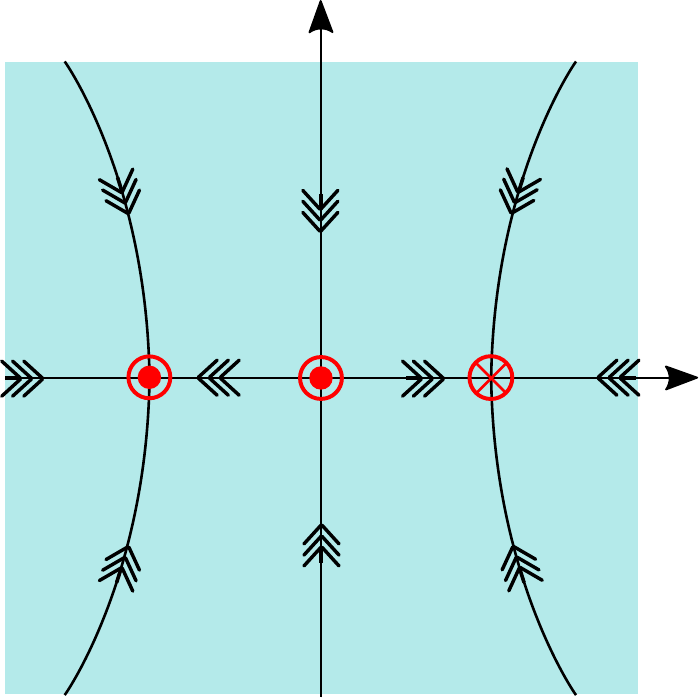}
	        \put(93,49){$x_2$}
	        \put(48,94){$y_2$}
	        \put(8,49){$C_1$}
	        \put(48,49){$C_2$}
	        \put(73,49){$C_3$}
        \end{overpic}
		\\ \vspace{8pt}
		\begin{overpic}[width=0.99\textwidth]{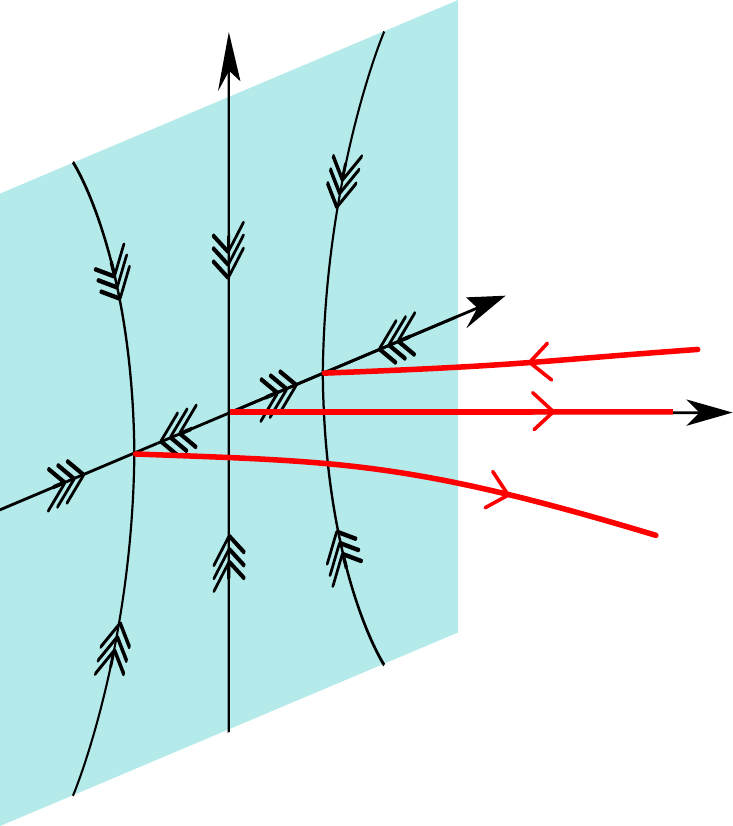}
	        \put(58,66){$x_2$}
	        \put(19,90){$y_2$}
	        \put(83,43){$\vec{z}$}
	        \put(17,38){$C_1$}
	        \put(29,56.5){$C_2$}
	        \put(42,63){$C_3$}
        \end{overpic}
		\caption{}\label{fig:eps1cubic_combi}
	\end{subfigure}
    \caption{	Sketches of possible limiting dynamics in the scaling chart $\kappa_2$, corresponding to $\bar{\varepsilon}=1$ in \cref{eq:blowup_def}.
			    {(a)} there exists one critical set $C$ which is a stable in the layer problem.
			    {(b)} there exists one critical set -- a saddle node in the layer problem.
			    (c) there are three critical sets -- two stable nodes and one saddle with respect to the fast flow, each with a different slow flow.}\label{fig:eps1_combi}
\end{figure}
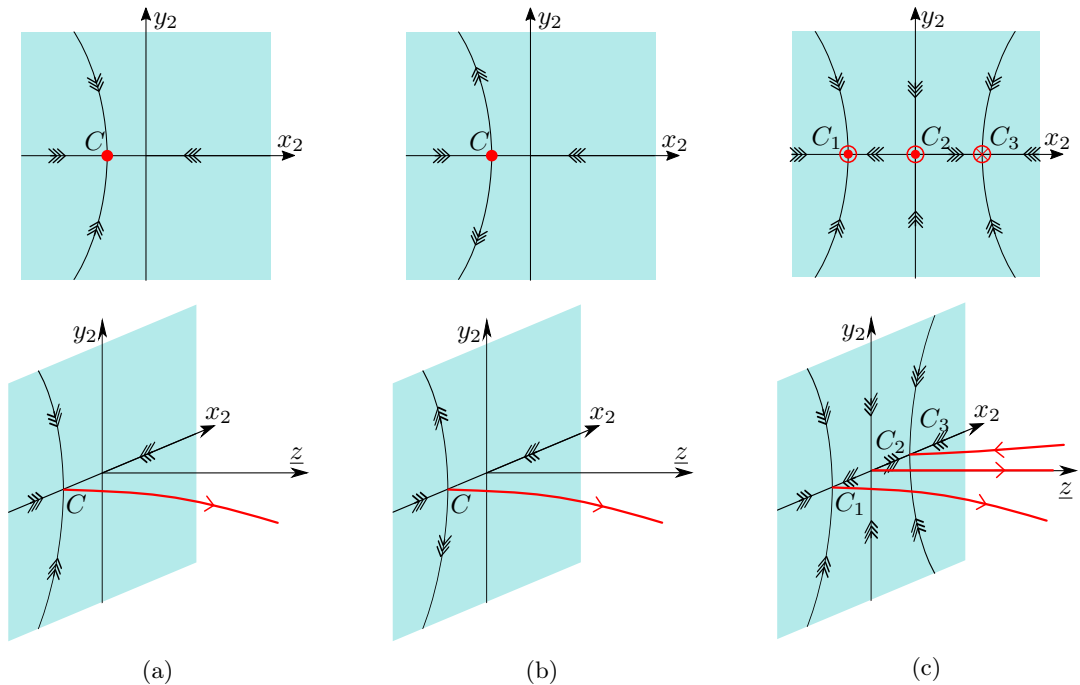

We recall from \cref{sec:regularisation} that the slow flow \cref{eq:slow_flow} corresponds to sliding dynamics. Examples of dynamics in the scaling chart $\kappa_2$ are shown in \cref{fig:eps1_combi}. \cref{fig:eps1node_combi} corresponds to stable sliding, \cref{fig:eps1saddle_combi} to a type of unstable sliding, and \cref{fig:eps1cubic_combi} to non-unique sliding.

\subsection{Summary}
We can piece together solutions to the full system \cref{eq:vfreg_extend} from the dynamics in charts $\kappa_1$ and $\kappa_2$.
In the regions where charts $\kappa_1$ and $\kappa_2$ overlap, we can track how orbits in the entry/exit chart appear in the scaling chart and vice versa, using the change of coordinates $\kappa_{12}$ \cref{eq:coordinatechange}. 

\begin{figure}[htbp]
	\centering
	\hfil
	\begin{subfigure}{0.3\textwidth}
        \begin{overpic}[width=0.9\textwidth]{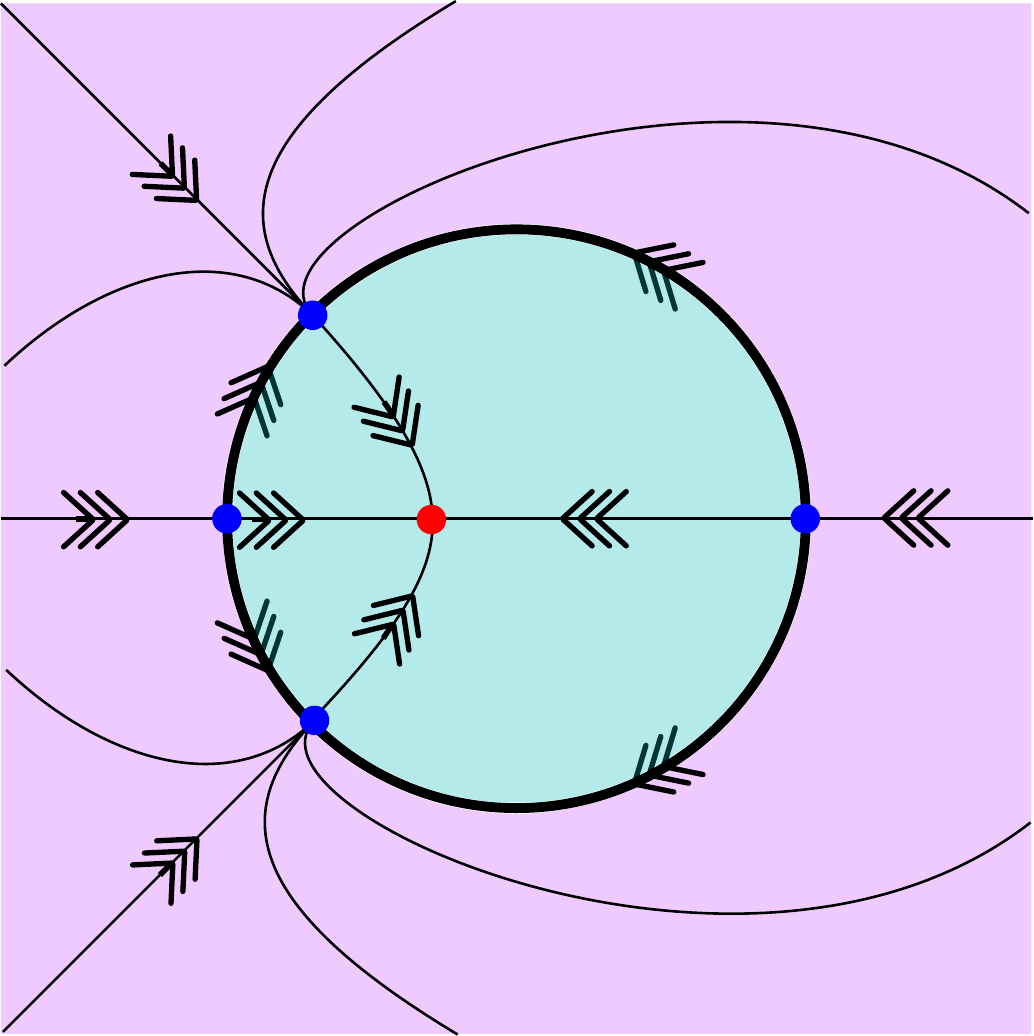}
	        \put(68,54.5){$Q_1$}
	        \put(25.5,54.5){$Q_3$}
	        \put(35,68){$Q_2$}
	        \put(35,28){$Q_4$}
	        \put(43,51){$C$}
        \end{overpic}
        \caption{}\label{fig:epsboth_node}
	\end{subfigure}
	\hfil
	\begin{subfigure}{0.3\textwidth}
        \begin{overpic}[width=0.9\textwidth]{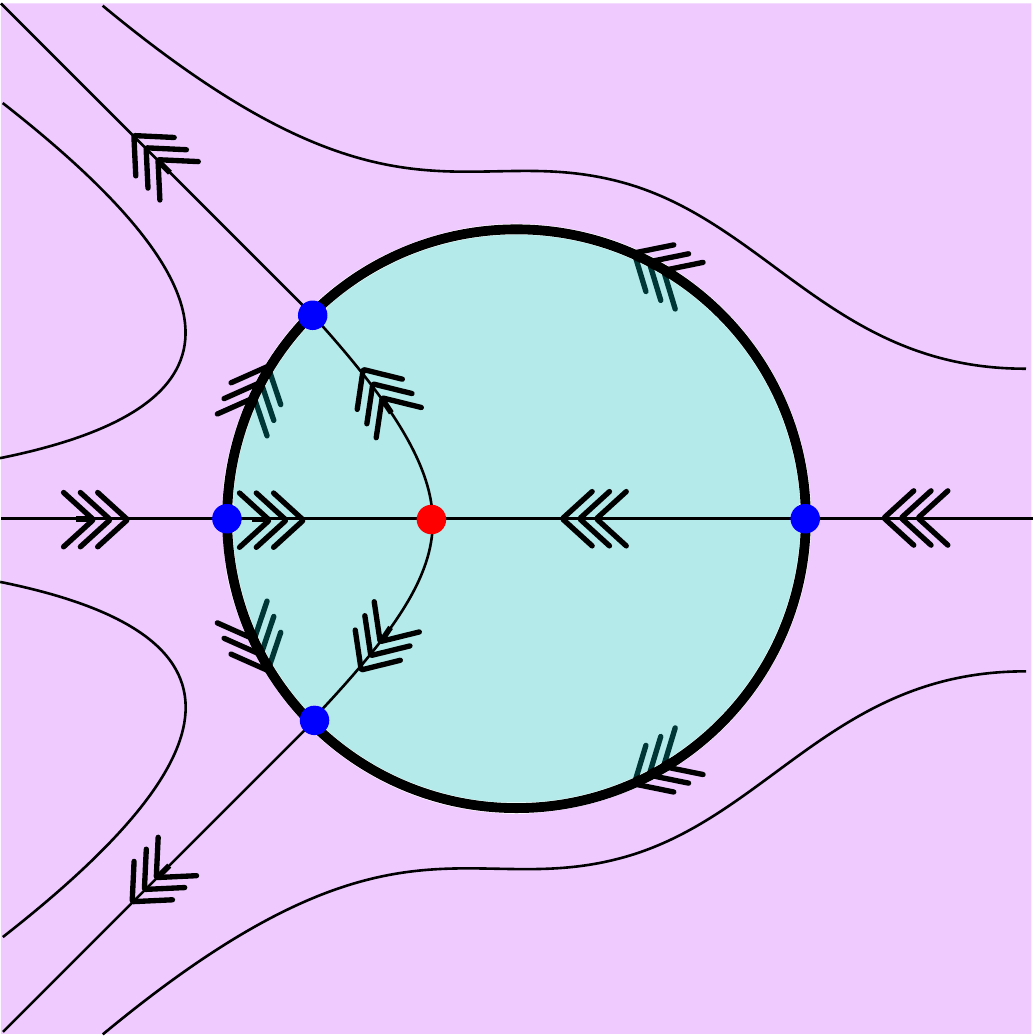}
	        \put(68,54.5){$Q_1$}
	        \put(25.5,54.5){$Q_3$}
	        \put(35,68){$Q_2$}
	        \put(35,28){$Q_4$}
	        \put(43,51){$C$}
        \end{overpic}
		\caption{}\label{fig:epsboth_saddle}
	\end{subfigure}
	\hfil
	\begin{subfigure}{0.3\textwidth}
        \begin{overpic}[width=0.9\textwidth]{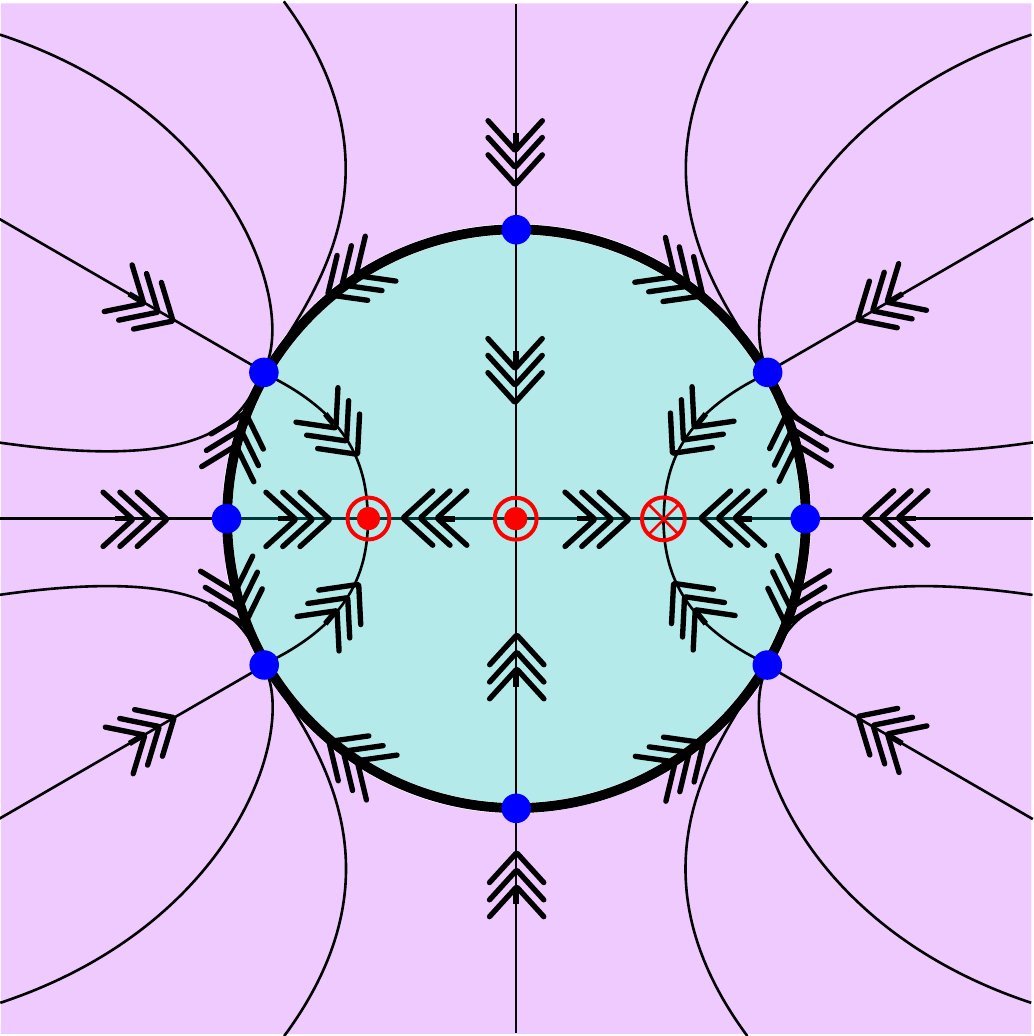}
        \end{overpic}
		\caption{}\label{fig:epsboth_cubic}
	\end{subfigure}
	\hfil
    \caption{	Sketches of examples of the dynamics of \cref{eq:vfreg_extend} after             the blowup in \cref{eq:blowup} projected into the $(\bar{x},\bar{y})$ plane.
    			Each subfigure pieces together the respective dynamics in the plane ${\varepsilon_1}=0$ of the chart $\kappa_1$
    			(\cref{fig:eps0}) and in the limit $\varepsilon\to0$ of the scaling chart $\kappa_2$ (\cref{fig:eps1_combi}).
    			Labels are omitted in (c), where {\color{red}{$\boldsymbol{\odot}$}} signifies that the slow flow along the critical manifold is moving out of the plane whilst {\color{red}{$\boldsymbol{\otimes}$}} signifies that it is moving into the plane.
    			We see that the sliding vector field is dependent upon the direction in which trajectories approach the discontinuity set.
    			Any trajectory that starts in the righthand plane tends to the rightmost critical manifold $C_3$, and hence travels `into the page',
    			yet any trajectory that starts in the lefthand plane tends to the leftmost critical manifold $C_1$, and hence travels `out of the page'.}\label{fig:epsboth}
\end{figure}

\begin{figure}[htbp]
	\centering
	\hfil
	\begin{subfigure}{0.3\textwidth}
        \begin{overpic}[width=0.99\textwidth]{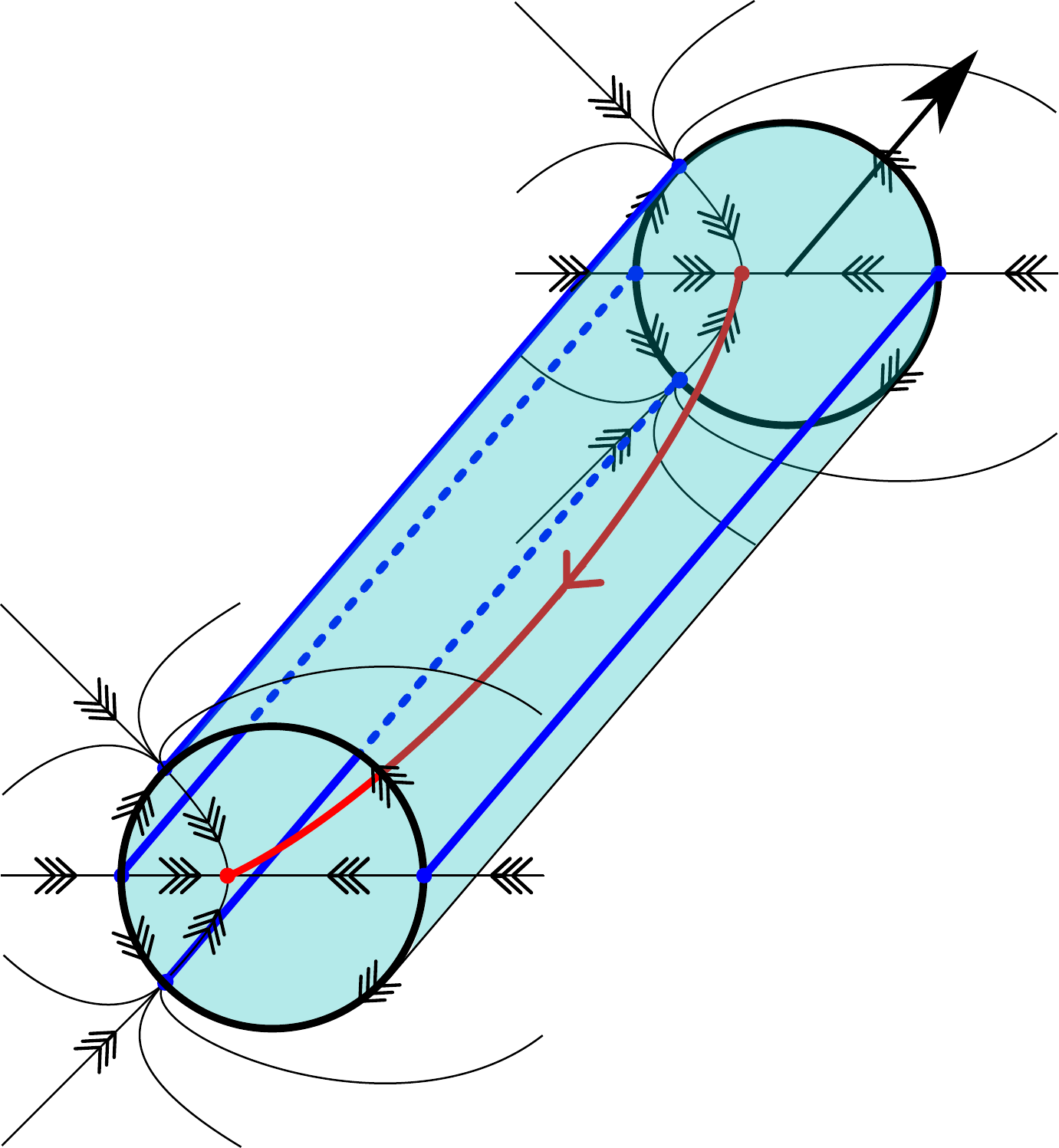}
        \end{overpic}
        \caption{}\label{fig:epsboth_node_3D}
	\end{subfigure}
	\hfil
	\begin{subfigure}{0.3\textwidth}
        \begin{overpic}[width=0.99\textwidth]{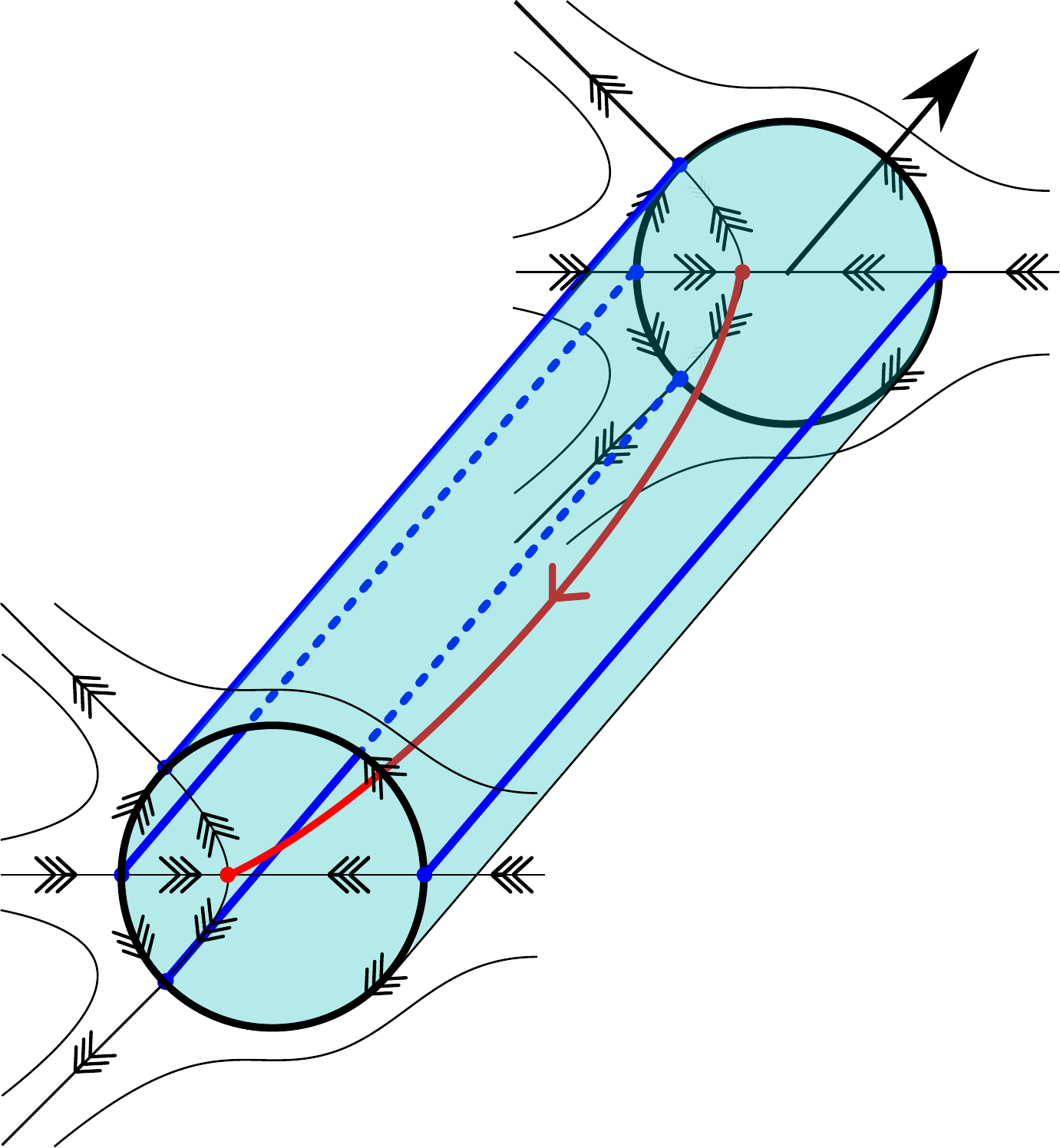}
        \end{overpic}
		\caption{}\label{fig:epsboth_saddle_3D}
	\end{subfigure}
	\hfil
	\begin{subfigure}{0.3\textwidth}
        \begin{overpic}[width=0.99\textwidth]{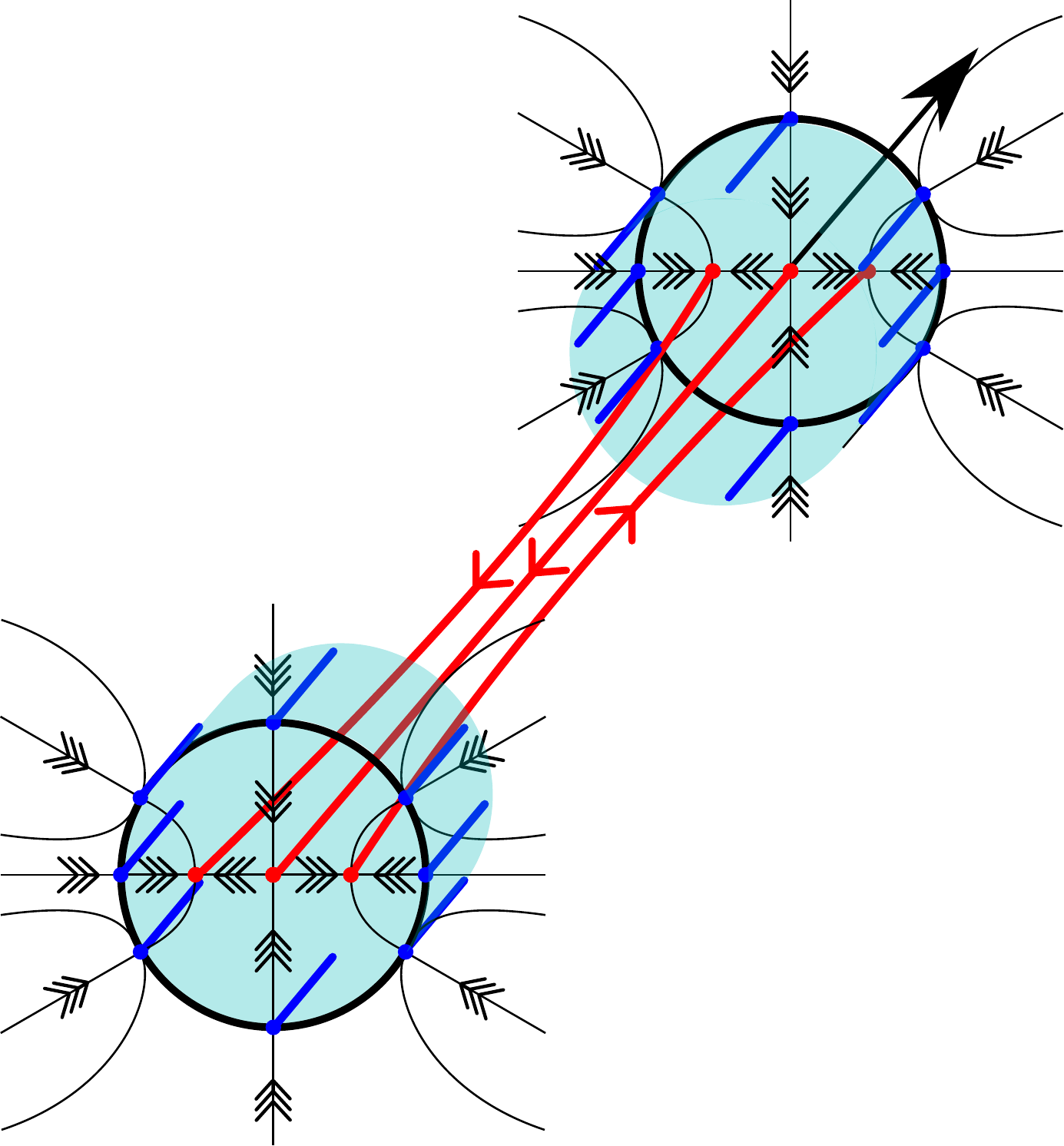}
        \end{overpic}
		\caption{}\label{fig:epsboth_cubic_3D}
	\end{subfigure}
	\hfil
    \caption{	Sketches of the dynamics projected onto $(\bar{x},\bar{y},\vec{z})$ space.
    			Each subfigure pieces together the respective dynamics in the plane ${\varepsilon_1}=0$
    			of the chart $\kappa_1$ (\cref{fig:eps0}) and in the limit $\varepsilon\to0$ of the scaling chart $\kappa_2$ (\cref{fig:eps1_combi}).
    			In (c), the projection of the equator is cut away, so that the slow flow is visible.}\label{fig:epsboth_3D}
\end{figure}

In \cref{fig:epsboth,fig:epsboth_3D}, we give sketches of the phase portraits of the dynamics of the full system \cref{eq:vfreg_extend} after the blowup in the nonsmooth limit $\varepsilon\to0$ for each of the three examples in  \cref{fig:eps0,fig:eps1_combi}, corresponding to stable sliding, unstable sliding and non-unique sliding respectively.

For a general $\vec{V}$, it is difficult to comment any further. Depending on the non-linearity of the vector field,
finding solutions to \cref{eq:Theta} will typically be difficult, even impossible,
and the dynamics in the scaling chart can become complex, including non-unique sliding
(as in \cref{fig:eps1cubic_combi,fig:epsboth_cubic,fig:epsboth_cubic_3D})
and limit cycles.
We have already noted that these dynamics may depend upon the regularisation function $\Psi$.
For that reason we shall analyse a simpler setting, where the vector field $\vec{V}$ is linear in $\vec{e}$, the so-called ``\texorpdfstring{$\vec{e}$}{e}-linear'' system.

\section{Classification of the ``\texorpdfstring{$\vec{e}$}{e}-linear'' system}\label{sec:lin}
Let us consider a normal form of \cref{eq:vf} that is linear in $\vec{e}$
\begin{equation}
	\begin{split}\label{eq:ode_normal_form}
		\begin{pmatrix}
			\dot{x}\\\dot{y}\\\dot{\vec{z}}
		\end{pmatrix}=
		\begin{pmatrix}
			\mat{A}(\vec{z})\\\mat{B}(\vec{z})
		\end{pmatrix}
		\vec{e}(x,y)+
		\begin{pmatrix}
			\vec{f}(x,y,\vec{z})\\
			\vec{g}(x,y,\vec{z})
		\end{pmatrix},
	\end{split}
\end{equation}
where 
$\mat{A}(\vec{z})\in\mathbb{R}^{2\times2}$, $\mat{B}(\vec{z})\in \mathbb{R}^{ m \times2}$, $\vec{f}(x,y,\vec{z})\in\mathbb{R}^2$
and $\vec{g}(x,y,\vec{z})\in\mathbb{R}^{m}$. Here $\mat{A}$, $\mat{B}$, $\vec{f}$ and $\vec{g}$ are assumed to be well defined at $\Sigma$. 
Without loss of generality (see \cref{sec:lemnormalform}), we can consider $\mat{A}(\vec{z})$ to be of the form
\begin{equation}\label{eq:matrix_normal_form}
	\mat{A}(\vec{z})=
	\begin{pmatrix}
		a(\vec{z})& -b(\vec{z})\\
		b(\vec{z})& \phantom{-}d(\vec{z})
	\end{pmatrix}.
\end{equation}

\begin{figure}[htbp]
	\centering
	\begin{subfigure}{0.45\textwidth}
		\centering
		\begin{overpic}[width=0.85\textwidth]{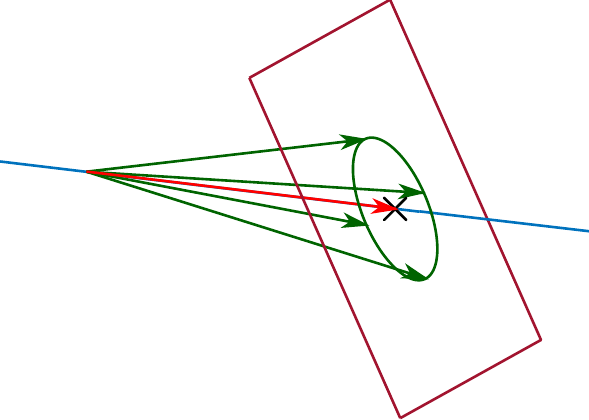}
		\put(7,46){$\Sigma$}
				\put(70,15){\color[RGB]{0,100,0}{$\vec{F}^*(\theta,\vec{z})$}}
				\put(71,65){\color[RGB]{162,20,47}{$\vec{F}_C(\lambda_1,\lambda_2,\vec{z})$}}
		\end{overpic}
		\caption{Sliding}\label{fig:whylinear_slide}
	\end{subfigure}
	\hfil
	\begin{subfigure}{0.45\textwidth}
		\centering
			\begin{overpic}[width=0.85\textwidth]{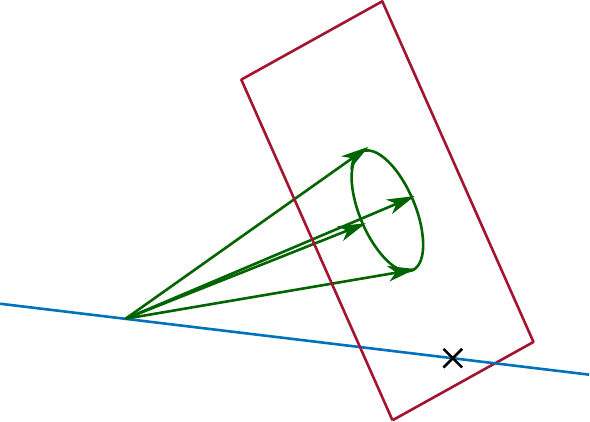}
			\put(7,22){$\Sigma$}
				\put(50,48){\color[RGB]{0,100,0}{$\vec{F}^*(\theta,\vec{z})$}}
				\put(70,65){\color[RGB]{162,20,47}{$\vec{F}_C(\lambda_1,\lambda_2,\vec{z})$}}
			\end{overpic}
			\caption{Crossing}\label{fig:whylinear_cross}
	\end{subfigure}
	\caption{Sliding and crossing in the $\vec{e}$-linear system \cref{eq:ode_normal_form,eq:matrix_normal_form}. In (a), there is sliding, when a convex combination of a pair of limit vectors is tangent to $\Sigma$. In (b), there is crossing, when no affine combination is convex. $\vec{F}^*(\theta,\vec{z})$ is the limit vector field, given in \cref{eq:lin_limit_vector}. $\vec{F}_C(\lambda_1,\lambda_2,\vec{z})$ is the plane of affine combinations of pairs of these limit vectors,  given in \cref{eq:lin_limit_plane}. A cross ({\bf $\times$}) marks the unique intersection of $\vec{F}_C(\lambda_1,\lambda_2,\vec{z})$ with $\Sigma$, given by \cref{eq:Fsigma}. }\label{fig:whylinear}
\end{figure}
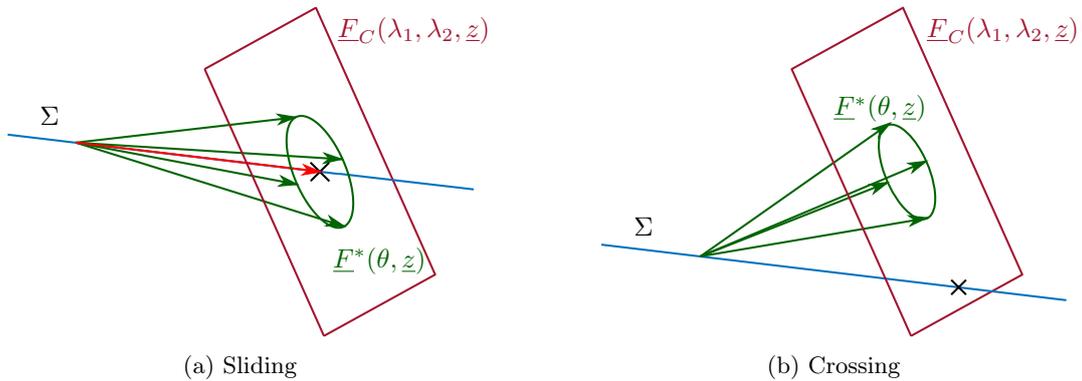

We study this particular class of systems for a number of reasons. Most importantly, the physical systems that have motivated this paper fit into this form.
Rigid body mechanical systems in 3D with Coulomb friction \cite{antali2019nonsmooth} can be written as \cref{eq:ode_normal_form}, with $b=0$ \cite{antali2019nonsmooth}.
In addition, the Filippov convention is well-posed here. With \cref{eq:matrix_normal_form}, combinations of pairs of incident vectors are coplanar, so there is a unique {\it affine} combination of incident vectors that is tangent to the discontinuity (see \cref{fig:whylinear}).
Where this affine combination is also a convex combination, there is sliding (\cref{fig:whylinear_slide}),
and where the combination is not convex, there is crossing (\cref{fig:whylinear_cross}).
Moreover, \cref{eq:ode_normal_form} is analogous to standard Filippv systems: the linear dependence on $\vec{e}$
relates to the linear dependence of standard Filippov problems on $\sign$ functions (unlike the nonlinear dependence of problems on $\sign$ functions in \cite{jeffrey2018hidden}).

Using the Filippov approach, the limit vector field \cref{eq:directional_limit} is given by
\begin{equation}\label{eq:lin_limit_vector}
	\vec{F}^*(\theta,\vec{z})=
	\begin{pmatrix}
		\mat{A}(\vec{z})\\\mat{B}(\vec{z})
	\end{pmatrix}
	\begin{pmatrix}
		\cos\theta\\\sin\theta
	\end{pmatrix}+
	\begin{pmatrix}
		\vec{f}(0,0,\vec{z})\\\vec{g}(0,0,\vec{z})
	\end{pmatrix}.
\end{equation}
The plane of affine combinations of pairs of these limit vectors $\vec{F}_C(\lambda_1,\lambda_2,\vec{z})$ is given by
\begin{equation}\label{eq:lin_limit_plane}
\begin{split}
	\vec{F}_C(\lambda_1,\lambda_2,\vec{z})&=\alpha_1 F^*(\theta_1,\vec{z})+\alpha_2 F^*(\theta_2,\vec{z}),\quad \theta_1\neq\theta_2,\quad\theta_1,\theta_2\in [0,2\pi),\quad \alpha_1+\alpha_2=1\\
	&=
	\begin{pmatrix}
		\mat{A}(\vec{z})\\\mat{B}(\vec{z})
	\end{pmatrix}
	\begin{pmatrix}
		\lambda_1\\\lambda_2
	\end{pmatrix}+
	\begin{pmatrix}
		\vec{f}(0,0,\vec{z})\\\vec{g}(0,0,\vec{z})
	\end{pmatrix}
	\end{split}
\end{equation}
where $(\lambda_1,\lambda_2)\in\mathbb{R}^2$. The intersection of this plane with $\Sigma$ is unique and given by 
\begin{equation}\label{eq:Fsigma}
	\vec{F}_\Sigma(\vec{z})=
	\begin{pmatrix}
		0\\0\\
		-\mat{B}(\vec{z})\mat{A}^{-1}(\vec{z})\vec{f}(0,0,\vec{z})+\vec{g}(0,0,\vec{z})
	\end{pmatrix}.
\end{equation}
$\vec{F}_C(\lambda_1,\lambda_2,\vec{z})$ is \textit{convex} when $(\lambda_1,\lambda_2)\in D_1$, the unit disc centred on the origin (see \cref{def:MPsi}), or
\begin{equation}
	(\mat{A}^{-1}(\vec{z})\vec{f}(0,0,\vec{z}))^\intercal (\mat{A}^{-1}(\vec{z})\vec{f}(0,0,\vec{z})) \leq 1.
\end{equation}
As already noted, the Filippov convention is well-posed here. Nevertheless it is worth studying the linear case \cref{eq:ode_normal_form,eq:matrix_normal_form} using regularisation in order to understand the robustness of the convention to perturbations and the stability of sliding.

So we regularise and augment \eqref{eq:ode_normal_form}, to give
\begin{equation}
\begin{split}\label{eq:odes_smoothed}
	\begin{pmatrix}
			{x}'\\{y}'\\{\vec{z}}'
	\end{pmatrix}&=\varepsilon\left(
	\begin{pmatrix}
		\mat{A}(\vec{z})\\\mat{B}(\vec{z})
	\end{pmatrix}\vec{e}_\Psi(x,y;\varepsilon)+
	\begin{pmatrix}
		\vec{f}(x,y,\vec{z})\\ \vec{g}(x,y,\vec{z})
	\end{pmatrix}\right),\\
	{\varepsilon}'&=0\
\end{split}
\end{equation}
and then study the dynamics of \cref{eq:odes_smoothed,eq:matrix_normal_form} in the charts $\kappa_1$ and $\kappa_2$.

\subsection{Entry chart\texorpdfstring{ $\kappa_1$}{}: dynamics near the discontinuity set}
As in the general case, we consider the augmented system and blowup $x=y=\varepsilon=0$ to a sphere using \cref{eq:blowup}.
We first study the equations in chart $\kappa_1$,
so that \cref{eq:polarising_desingularising} becomes
\begin{equation}
	\begin{split}\label{eq:normal_reg_polar}	
		\frac{\mathrm{d}}{\mathrm{d}\mathcal{T}}\begin{pmatrix}
			{\rho}\\{\theta}
		\end{pmatrix}&=
		\begin{pmatrix}
			\rho\xi&0\\
			0&\xi^{-1}
		\end{pmatrix}
		\mat{R}^\intercal(\theta)\left(\mat{A}(\vec{z})\vec{\hat{e}}_\Psi(\theta,\varepsilon_1)
		+\vec{f}(\rho\xi\cos\theta,\rho\xi\sin\theta,\vec{z})\right),\\
		\frac{\mathrm{d}}{\mathrm{d}\mathcal{T}}{\vec{z}}&=\rho\left(\mat{B}(\vec{z})\vec{\hat{e}}_\Psi(\theta,\varepsilon_1)+
		\vec{g}(\rho\xi\cos\theta,\rho\xi\sin\theta,\vec{z})\right),\\
		\frac{\mathrm{d}}{\mathrm{d}\mathcal{T}}{\varepsilon}_1&=-\frac{\varepsilon_1}{\rho}\dot{\rho}.
	\end{split}
\end{equation}
The equator $\rho=\varepsilon_1=0$ is once again an invariant manifold, with equilibria when
\begin{equation}\label{eq:equator_zeros}
	\begin{split}
		\Theta(0,\theta,\vec{z})&=(0,1) \mat{R}^\intercal(\theta) \left( \mat{A}(\vec{z})
		\begin{pmatrix}
			\cos{\theta}\\ \sin{\theta}
		\end{pmatrix}
		+ \vec{f}(0,0,\vec{z})\right)\\
		&=(d(\vec{z})-a(\vec{z}))\sin{\theta}\cos{\theta}+f_2(0,0,\vec{z}) \cos{\theta}-f_1(0,0,\vec{z}) \sin{\theta}+b(\vec{z})=0.
	\end{split}
\end{equation}

Solutions of \cref{eq:equator_zeros} can be viewed as the intersection of the hyperbola
\begin{align}
	\begin{pmatrix}
		\bar{x}&\bar{y}
	\end{pmatrix}
	\begin{pmatrix}
		0 &\frac{(d(\vec{z})-a(\vec{z}))}{2}\\
		\frac{(d(\vec{z})-a(\vec{z}))}{2}&0
	\end{pmatrix}
	\begin{pmatrix}
		\bar{x}\\\bar{y}
	\end{pmatrix}+
	\begin{pmatrix}
		f_2(0,0,\vec{z})&-f_1(0,0,\vec{z})
	\end{pmatrix}
	\begin{pmatrix}
		\bar{x}\\\bar{y}
	\end{pmatrix}
	+b(\vec{z})=0 \label{eq:conic11}
\end{align}
with the unit circle (the equator)
\begin{align}
	\begin{pmatrix}
		\bar{x}&\bar{y}
	\end{pmatrix}
	\begin{pmatrix}
		1&0\\0&1
	\end{pmatrix}
	\begin{pmatrix}
		\bar{x}\\ \bar{y}
	\end{pmatrix}-1=0. \label{eq:conic12}
\end{align}

\begin{figure}[htbp]
	\centering
	\hfil
	\begin{subfigure}[t]{0.23\textwidth}
		\centering
		\begin{overpic}[width=0.94\textwidth]{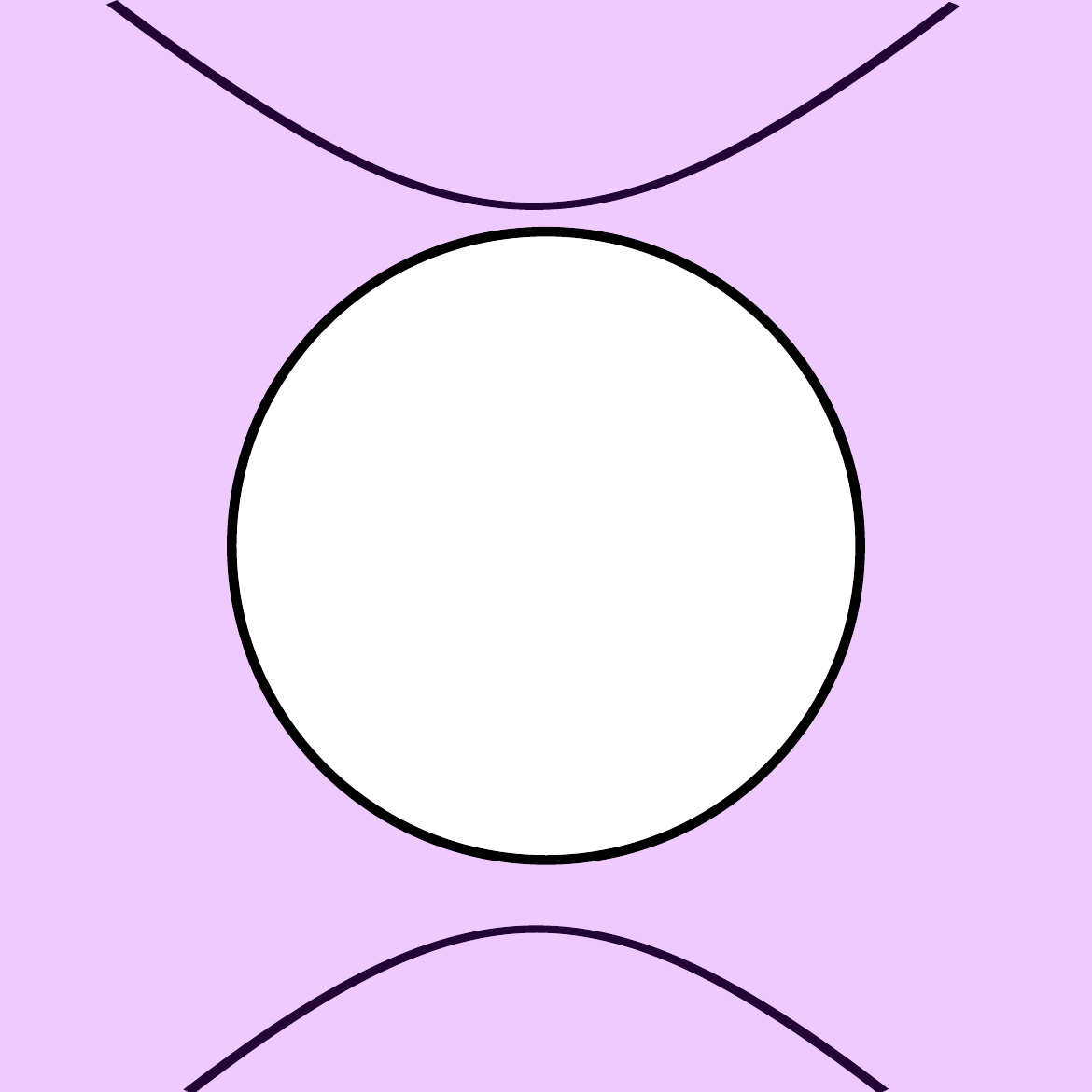}
		\end{overpic}
		\caption{0 intersections}
	\end{subfigure}
	\hfil
	\begin{subfigure}[t]{0.23\textwidth}
		\centering
		\begin{overpic}[width=0.94\textwidth]{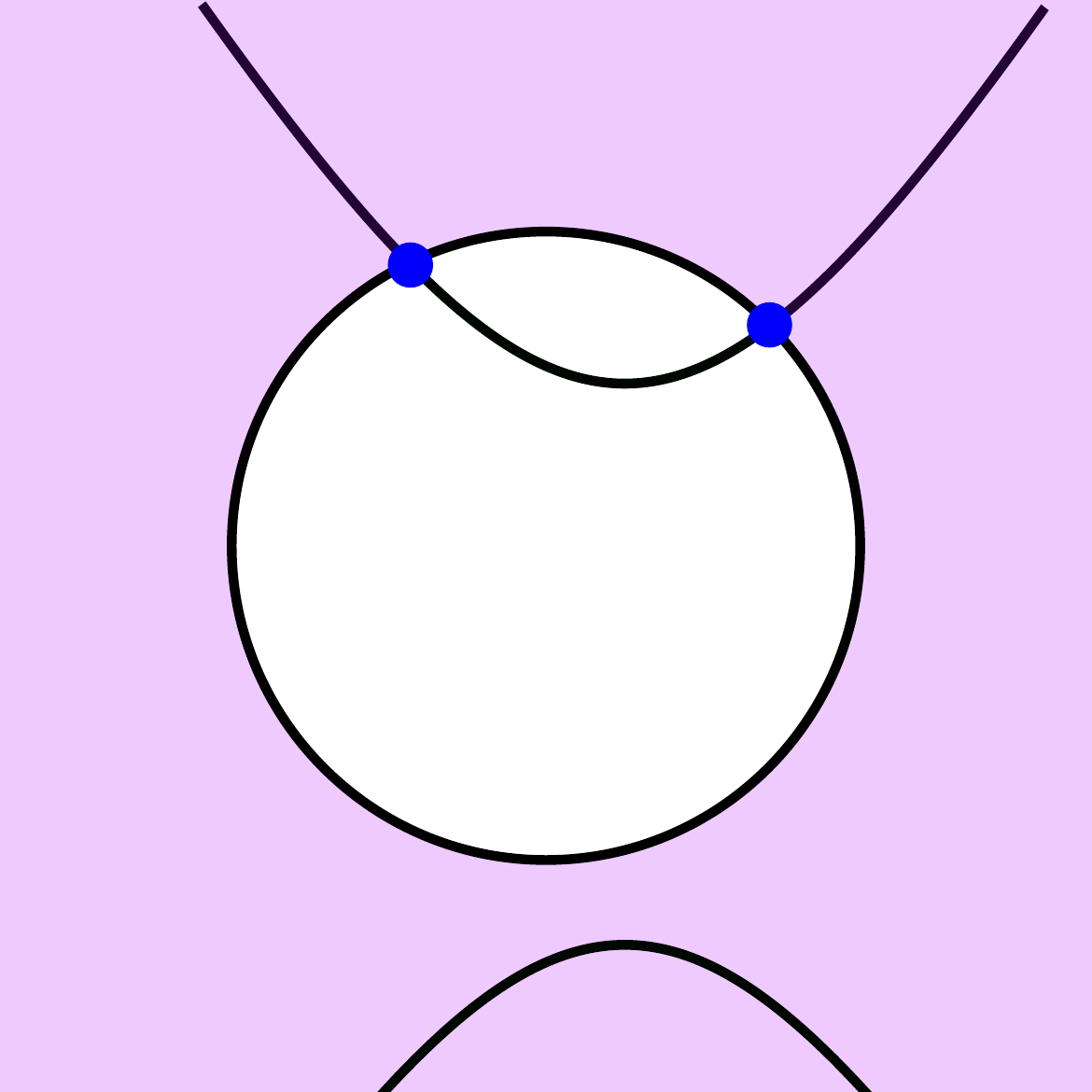}
		\end{overpic}
		\caption{2 intersections}
	\end{subfigure}
	\hfil
	\begin{subfigure}[t]{0.46\textwidth}
		\centering
		\begin{overpic}[width=0.47\textwidth]{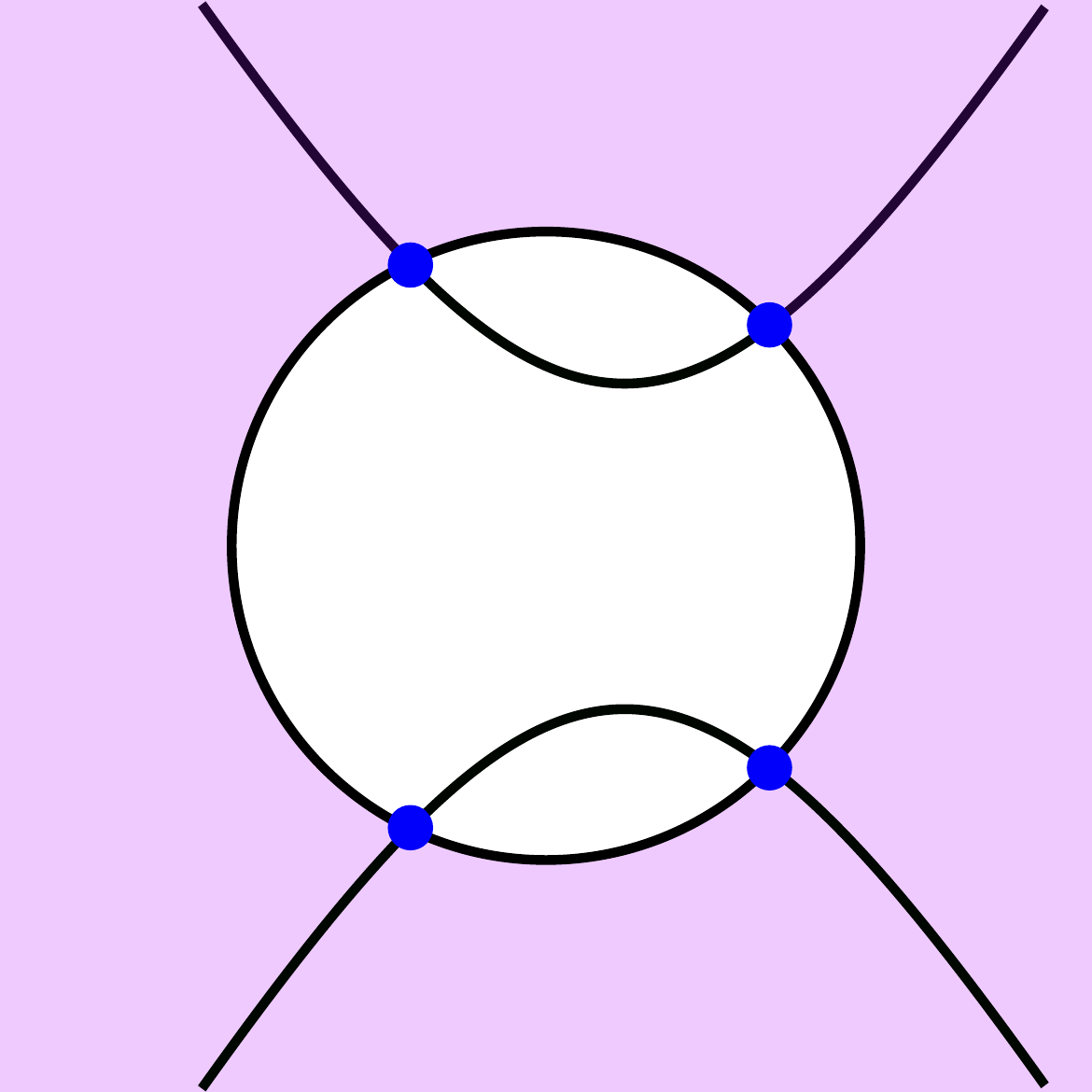}
		\end{overpic}
		\begin{overpic}[width=0.47\textwidth]{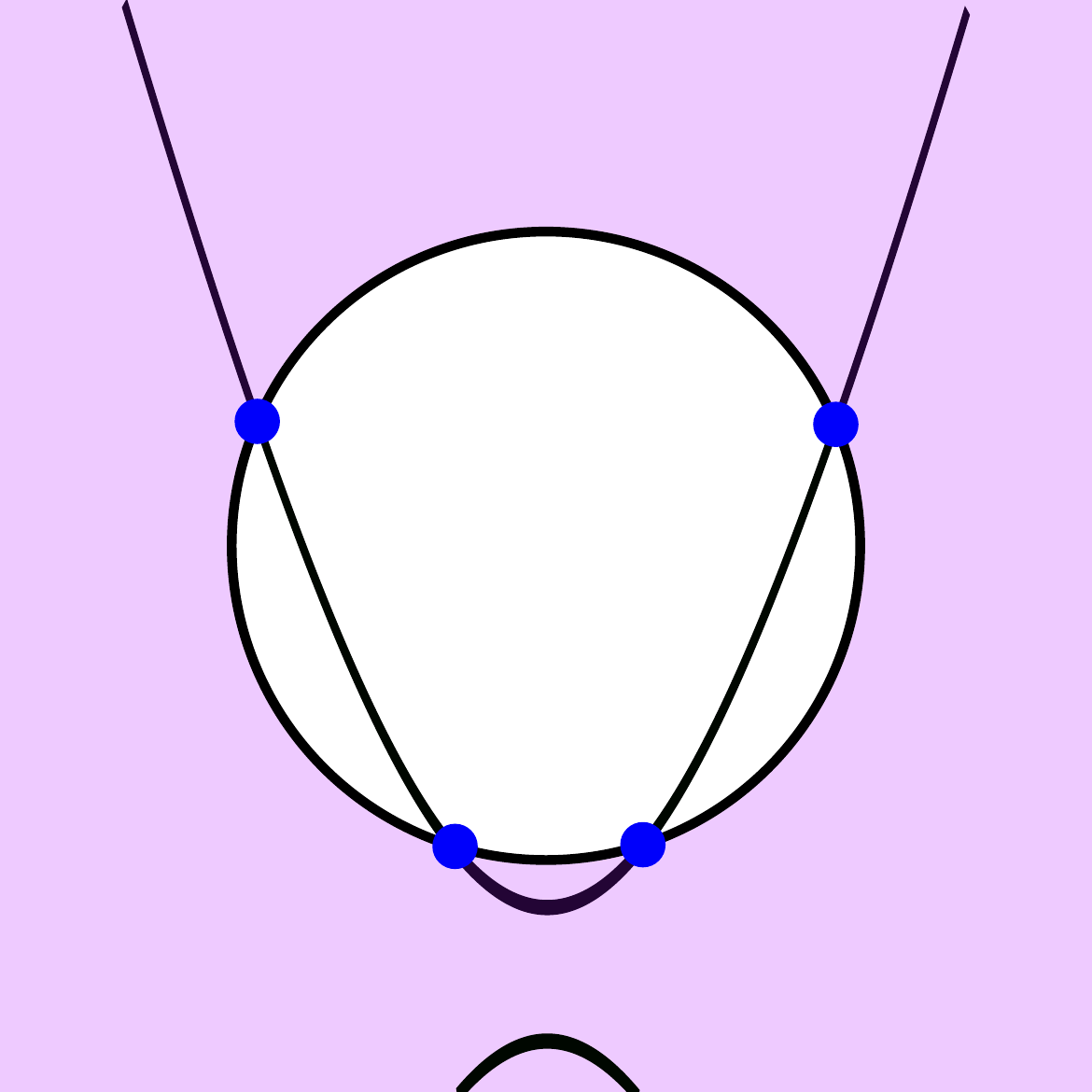}
		\end{overpic}
		\caption{4 intersections}
	\end{subfigure}
	\hfil
	\caption{	Intersections of conics: the ways that the hyperbola \cref{eq:conic11} and the unit circle \cref{eq:conic12} can intersect generically. There can be only 0, 2 or 4 intersections of these conics and hence  only 0, 2 or 4
				equilibria on the equator (or ``limit directions'') for our linear system \cref{eq:ode_normal_form}. Each subfigure depicts the ($\bar{x},\bar{y}$) plane and the unit circle corresponds to the equator of the blowup.}\label{fig:intersections}
\end{figure}

Hence there can be at most 4 unique equilibria. This is clear either algebraically when writing \cref{eq:conic11} and \cref{eq:conic12} as a single quartic in either $\bar{x}$ or $\bar{y}$, or geometrically in \cref{fig:intersections}. Generically, equilibria on the equator are created in saddle node bifurcations (\cref{sec:bifurcations}). These occur where \cref{eq:eqbif} is satisfied, resulting in
\begin{equation}\label{eq:tangencies}
    \begin{split}
        f_1&=b\sin\theta- (a-d)\cos^3\theta\\
        f_2&=(a-d)\sin^3\theta-b\cos\theta.\\
    \end{split}
\end{equation}
However, closed form solutions to \cref{eq:tangencies} are too lengthy to be of much use.
In our geometric interpretation, this bifurcation occurs when the hyperbola \cref{eq:conic11} is tangent to the unit circle \cref{eq:conic12}.

There are also degenerate cases of \cref{eq:equator_zeros} to consider.
\begin{enumerate}[label=(D\arabic*)]
	\item{	$a(\vec{z})=d(\vec{z})$:\\
			In this case, unique equilibria are given by
			$$\theta=\theta^*(\vec{z})=\arcsin{\left(\frac{b(\vec{z})}{|\vec{f}(0,0,\vec{z})|}\right)}+\arctan{\left(\frac{f_2(0,0,\vec{z})}{f_1(0,0,\vec{z})}\right)}.$$
			Here we have no equilibria when $|b(\vec{z})|>|\vec{f}(0,0,\vec{z})|$
			and 2 equilibria when $|b(\vec{z})|<|\vec{f}(0,0,\vec{z})|$ (with a saddle node bifurcation when $|b(\vec{z})|=|\vec{f}(0,0,\vec{z})|$).
			}\label{deg:1}
	\item{	At least one of $f_1(0,0,\vec{z})$ or $f_2(0,0,\vec{z})$ is zero and $b(\vec{z})=0$:\\
			Without loss of generality we shall suppose that $f_2(0,0,\vec{z})=0$, then \cref{eq:equator_zeros} becomes
			$$\sin{\theta}\left((a(\vec{z})-d(\vec{z}))\cos{\theta}-f_1(0,0,\vec{z})\right)=0.$$
			In this case we always have two equilibria $\theta=\theta^*(\vec{z})=0,\pi$ on the equator.
			Other equilibria are born in pitchfork bifurcations at $|f_1(0,0,\vec{z})|=|d(\vec{z})-a(\vec{z})|$, given by			
			$$\theta=\theta^*(\vec{z})=\pm\arccos{\left(\frac{f_1(0,0,\vec{z})}{d(\vec{z})-a(\vec{z})}\right)}$$ 
			when $|f_1(0,0,\vec{z})|<|d(\vec{z})-a(\vec{z})|$.
			The analysis when $f_2(0,0,\vec{z})=0$ and $b(\vec{z})=0$ is similar.}\label{deg:2}
\end{enumerate}

\subsection{Scaling chart\texorpdfstring{ $\kappa_2$}{}: dynamics on the discontinuity set}
Following the same procedure as for the general system we study the dynamics in the scaling chart $\kappa_2$, setting $\bar{\varepsilon}=1$, to obtain
\begin{equation}\label{eq:linear_sc}
	\begin{split}
		\varepsilon(\dot{x}_2,\dot{y}_2)^\intercal&=\mat{A}(\vec{z})\vec{e}_\Psi(x_2,y_2;1)+\vec{f}(\varepsilon x_2,\varepsilon y_2,\vec{z}),\\
		\dot{\vec{z}}&=\mat{B}(\vec{z})\vec{e}_\Psi(x_2,y_2;1)+\vec{g}(\varepsilon x_2,\varepsilon y_2,\vec{z}).
	\end{split}
\end{equation}
As expected we obtain a slow-fast system where $x_2$ and $y_2$ are fast.

\begin{theorem}\label{thm:lineps1}Considering the dynamics in the scaling chart $\kappa_2$ given by \cref{eq:linear_sc}, then:
\begin{enumerate}[\rm{(\alph*)}]
	\item{	Where $\mat{A}(\vec{z})$ and $\vec{f}(0,0,\vec{z})$ satisfy 
			\begin{equation}\label{eq:fcondition}
				(\mat{A}(\vec{z})^{-1} \vec{f}(0,0,\vec{z}))^\intercal \mat{A}(\vec{z})^{-1} \vec{f}(0,0,\vec{z}) < 1,
			\end{equation}
			there exists exactly one unique equilibrium of the layer problem in the scaling chart $\kappa_2$ given by
			\begin{equation}\label{eq:lincriticalset}
				C=\left\{(x_2,y_2,\vec{z})|\mat{A}(\vec{z})\vec{e}_\Psi(x_2,y_2;1)+\vec{f}(0,0,\vec{z})=0\right\}
			\end{equation}
			for a fixed $\vec{z}$. On the other hand, if
			\begin{equation}\label{eq:fcondition2}
				(\mat{A}(\vec{z})^{-1} \vec{f}(0,0,\vec{z}))^\intercal \mat{A}(\vec{z})^{-1} \vec{f}(0,0,\vec{z}) > 1,
			\end{equation}
			then there is no such critical set.
			\label{thm:lineps1a}}
	\item{	Where a unique critical set $C$ exists, the slow flow in the reduced problem is given by
			\begin{equation}\label{eq:linslowflow}
				\dot{\vec{z}}=-\mat{B}(\vec{z})\mat{A}^{-1}(\vec{z})\vec{f}(0,0,\vec{z})+\vec{g}(0,0,\vec{z}),
			\end{equation}
			which coincides with the Filippov convention \cite[(5.20)]{Antali2017}.\label{thm:lineps1b}}
\end{enumerate}
\end{theorem}

\begin{figure}[htbp]
				\centering
				\hfil
				\begin{subfigure}{0.42\textwidth}
					\begin{overpic}[width=0.9\textwidth]{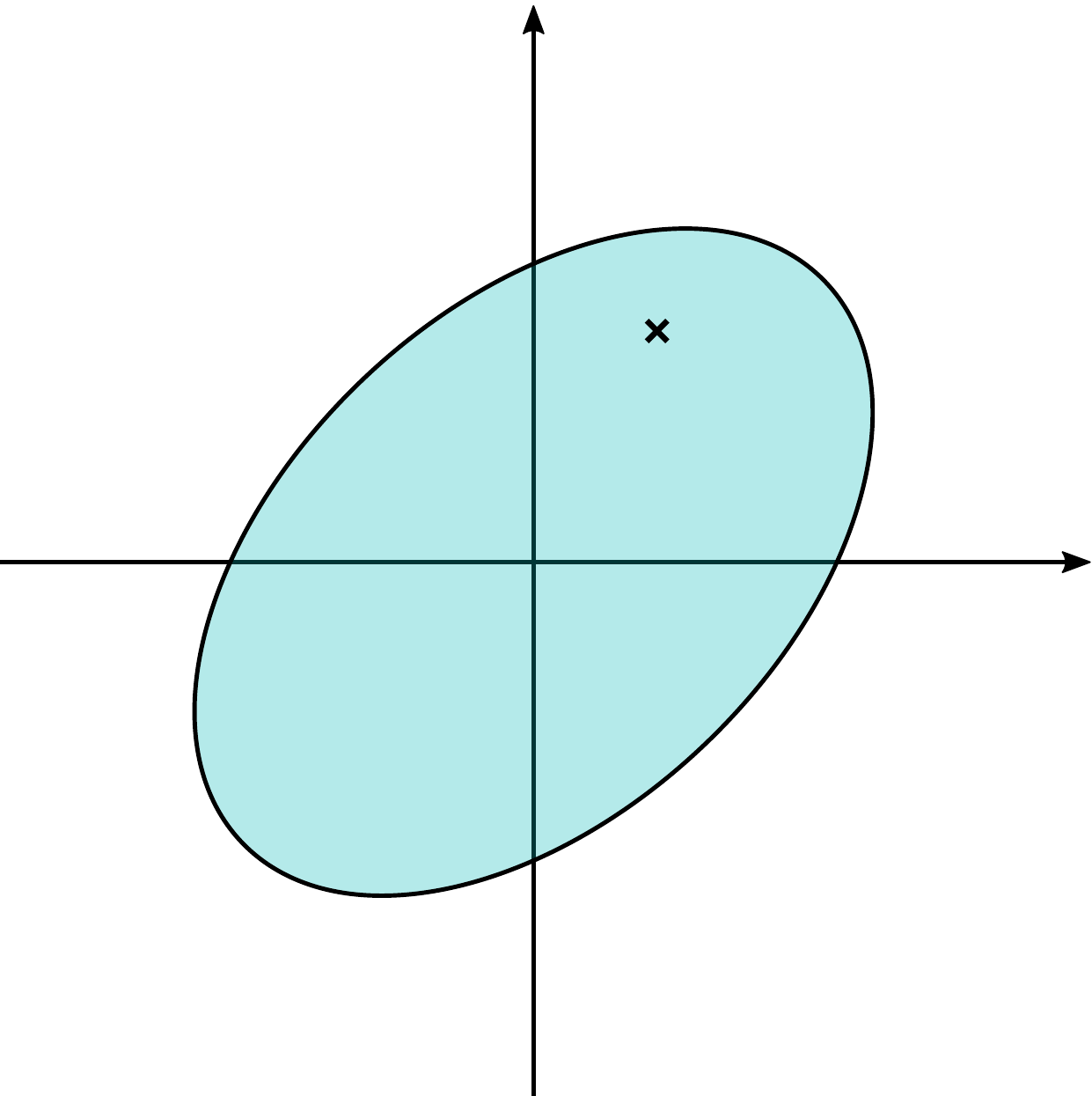}
						\put(92,42){{\rotatebox{0}{$h_1$}}}
						\put(40,92){{\rotatebox{0}{$h_2$}}}
						\put(73,86){$\vec{h}=\vec{f}(0,0,\vec{z})$}
						\put(60,70){\line(2,1){23}}
						\put(24,9){\line(1,2){4}}
						\put(0,4){\colorbox{white}{$( \mat{A}(\vec{z})^{-1} \vec{h})^\intercal \mat{A}(\vec{z})^{-1} \vec{h}=1$}}
					\end{overpic}
					\caption{	$\vec{f}(0,0,\vec{z})$ is within the ellipse and hence there is a critical set within the scaling chart}
				\end{subfigure}
				\hfil
				\begin{subfigure}{0.42\textwidth}
					\begin{overpic}[width=0.9\textwidth]{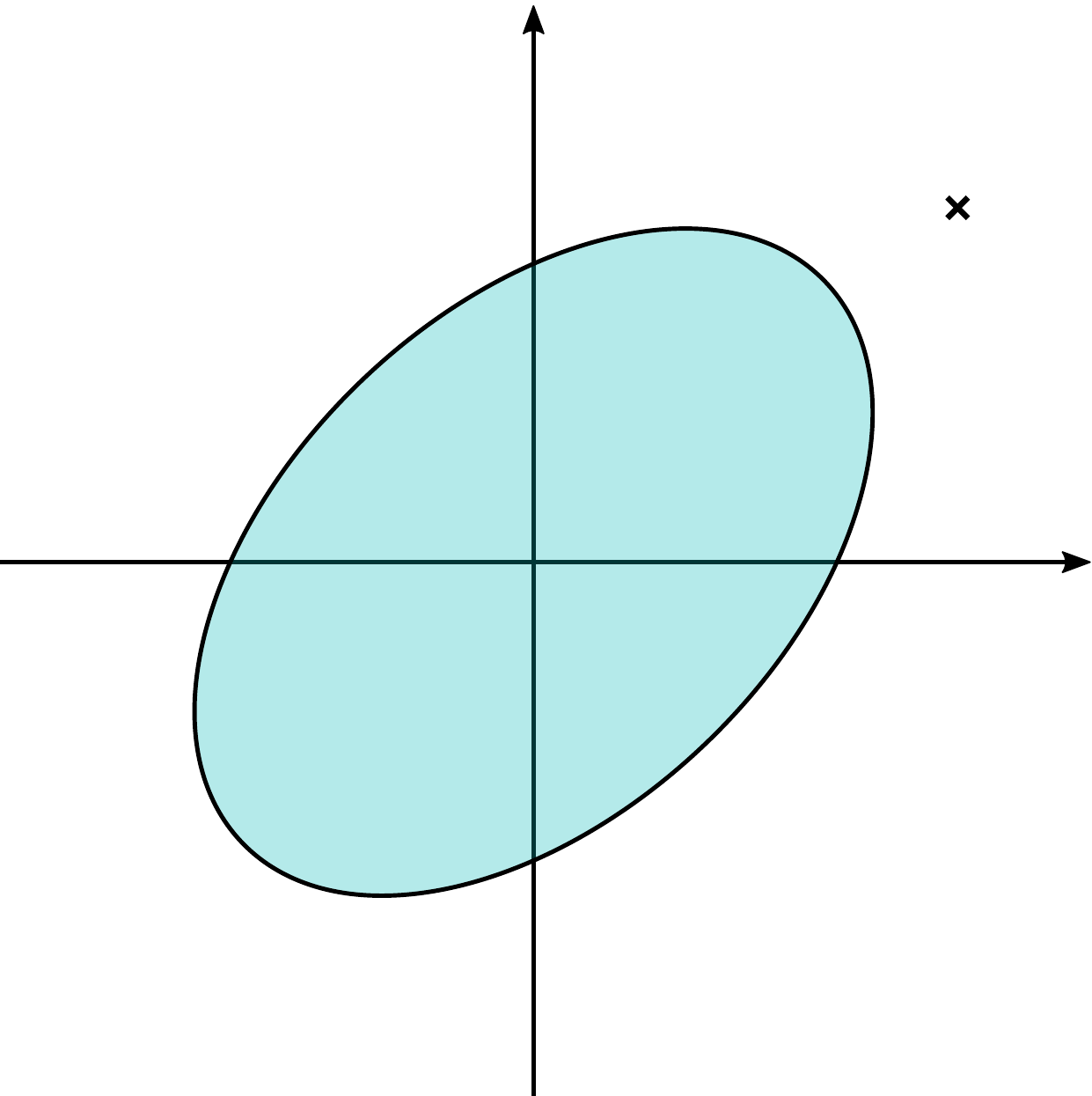}
						\put(92,42){{\rotatebox{0}{$h_1$}}}
						\put(40,92){{\rotatebox{0}{$h_2$}}}
						\put(73,86){$\vec{h}=\vec{f}(0,0,\vec{z})$}
						\put(24,9){\line(1,2){4}}
						\put(0,4){\colorbox{white}{$( \mat{A}(\vec{z})^{-1} \vec{h})^\intercal \mat{A}(\vec{z})^{-1} \vec{h}=1$}}
					\end{overpic}
					\caption{	$\vec{f}(0,0,\vec{z})$ is outside the ellipse and hence there is no critical set within the scaling chart}
				\end{subfigure}
				\hfil
				\caption{	Sketch of a geometric interpretation of the condition on $\vec{f}(0,0,\vec{z})$ for the existence of a critical
				 set in the scaling chart.}\label{fig:ellipse}
			\end{figure}
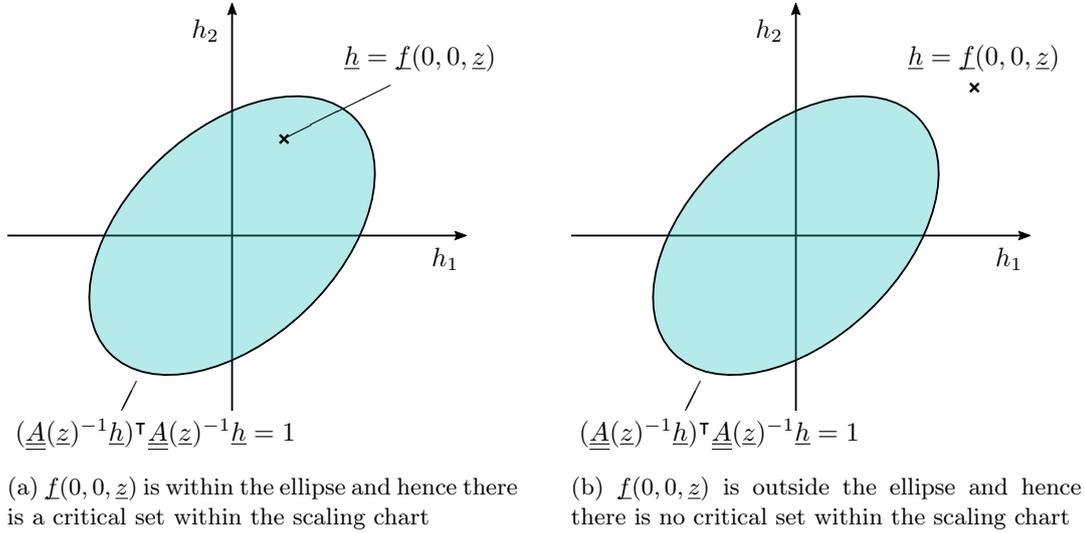

\begin{proof}
Rescaling time in \cref{eq:linear_sc}, we obtain the system in fast time
\begin{equation}
	\begin{split} 
		\left(x_2',y_2'\right)^\intercal
		&=\mat{A}(\vec{z}) \vec{e}_\Psi(x_2,y_2;1)+\vec{f}(\varepsilon x_2,\varepsilon y_2,\vec{z}),\\
		\vec{z}'
		&=\varepsilon \left(\mat{B}(\vec{z})\vec{e}_\Psi(x_2,y_2;1)+\vec{g}(\varepsilon x_2, \varepsilon y_2, \vec{z})\right).
	\end{split}
\end{equation}
The layer problem is given by the limit $\varepsilon\to 0$
\begin{equation}
	\begin{split} 
		\left(x_2',y_2'\right)^\intercal
		&=\mat{A}(\vec{z}) \vec{e}_\Psi(x_2,y_2;1)+\vec{f}(0,0,\vec{z})\\
		\vec{z}'&=0.
	\end{split}\label{eq:linLP}
\end{equation}
\begin{enumerate}[\rm{(\alph*)}]
	\item{	A unique critical set \cref{eq:lincriticalset}
			exists for a given $\vec{z}$ in the scaling chart $\kappa_2$ if there is a solution to
			\begin{equation}\label{eq:linepsisoln}
				\vec{e}_\Psi(x_2,y_2;1)=-\mat{A}^{-1}(\vec{z}) \vec{f}(0,0,\vec{z}).
			\end{equation}
			When \cref{eq:fcondition} holds, solutions to \cref{eq:linepsisoln} exist, since $\vec{e}_\Psi^\intercal \vec{e}_\Psi<1$, from \cref{eq:ePsimod}.
			If the critical set exists, it is unique, from \cref{lem:ePsi}\cref{lem:ePsic}.
			When \cref{eq:fcondition2} holds, no solutions exist to \cref{eq:linepsisoln}, and therefore there can be no critical set $C$. Condition \cref{eq:linepsisoln}, together with \cref{eq:ePsimod}, can be thought of as $\vec{f}(0,0,\vec{z})$ lying within the ellipse 
			\begin{equation}( \mat{A}(\vec{z})^{-1} \vec{h})^\intercal \mat{A}(\vec{z})^{-1} \vec{h}=1\end{equation} since $\det\left(\left(\mat{A}(\vec{z})^{-1}\right)^\intercal\mat{A}(\vec{z})^{-1}\right)=\det(\mat{A}(\vec{z}))^{-2}$ (see \cref{fig:ellipse}).}
	\item{	If the critical set \cref{eq:lincriticalset} exists, the slow flow along it is given by \cref{eq:linslowflow}, found by substituting \cref{eq:linepsisoln}
			into the reduced problem \cref{eq:linear_sc}.}
	\end{enumerate}
\end{proof}
 As with the general system, the slow flow \cref{eq:linslowflow} corresponds to the sliding vector field in Filippov terminology.

Limit cycles may also occur in the scaling chart $\kappa_2$. Hyperbolic limit cycles persist and therefore form a cylinder. The reduced flow along such a cylinder is defined by an average of the slow component of the vector field.

\begin{proposition}\label{prop:LC}
If there exists a limit cycle in layer problem of the scaling chart, then the reduced flow along the limit cycle is given by \cref{eq:linslowflow}.\label{thm:lineps1c}
\end{proposition}
\begin{proof}
If we suppose that there exists a limit cycle where $(x_2(\tau+T),y_2(\tau+T))=(x_2(\tau),y_2(\tau))$,
			then averaging the layer problem \cref{eq:linLP} over the period $T$ we find
			\begin{align}
				\frac{1}{T}\int_0^T 
				\begin{pmatrix}
					{x}_2(\tau)'\\ {y}_2(\tau)'
				\end{pmatrix} \mathrm{d}\tau 
				&=\frac{1}{T} \int_0^T \Big( \mat{A}(\vec{z}) \vec{e}_\Psi(x_2(\tau),y_2(\tau);1)+\vec{f}(0,0,\vec{z})\Big) \mathrm{d}\tau\\
				\implies \vec{0}&=\mat{A}(\vec{z})\, \left(\frac{1}{T}\int_0^T \Big(\vec{e}_\Psi(x_2(\tau),y_2(\tau);1)\Big)\, \mathrm{d}\tau\right)+\vec{f}(0,0,\vec{z})\label{eq:LCstep}
			\end{align}
			Accordingly, the reduced problem along the limit cycle will be given by 
			\begin{equation}
				\dot{\vec{z}}=\mat{B}(\vec{z}) \, \left(\frac{1}{T}\int_0^T \Big(\vec{e}_\Psi(x_2(\tau),y_2(\tau);1)\Big)\, \mathrm{d}\tau\right) + \vec{g}(0,0,\vec{z}),
			\end{equation}			
			which, using \cref{eq:LCstep}, results in \cref{eq:linslowflow}.
\end{proof}

It is also possible to infer some of the dynamics in the scaling chart $\kappa_2$ purely from an examination of the entry chart $\kappa_1$ using arguments from index theory (as was done in \cref{sec:geneps0} when there were no equilibria on the equator).
\subsection{Cases}
\label{sec:cases}
In order to classify the possible global dynamics of \cref{eq:odes_smoothed} we shall now consider 3 cases of our $\vec{e}$-linear system \cref{eq:ode_normal_form}, determined by the elements of $\mat{A}(\vec{z})$ given by \cref{eq:matrix_normal_form}, and described in \cref{tab:cases}. 
\begin{table}[htbp]
	\centering
	\def\arraystretch{1.5}
	\begin{tabular}{|l||l|l|}
		\hline
		Case I  & $\sign(a)=\sign(d)$                       & $\det(\mat{A})>0$  \\ \hline
		Case II & \multirow{2}{*}{$\sign(a)\neq\sign(d)$}   & $\det(\mat{A})<0$      \\ \cline{1-1}\cline{3-3}
		Case III &                                          & $\det(\mat{A})>0$           \\ \hline
	\end{tabular}
	\caption{The 3 cases of the linear normal form \cref{eq:ode_normal_form} with \cref{eq:matrix_normal_form}}\label{tab:cases}
\end{table}

\subsubsection{Case I}
\begin{figure}[htbp]
	\centering
	\hfil
	\begin{subfigure}{0.3\textwidth}
		\centering
		\begin{overpic}[width=0.9\textwidth]{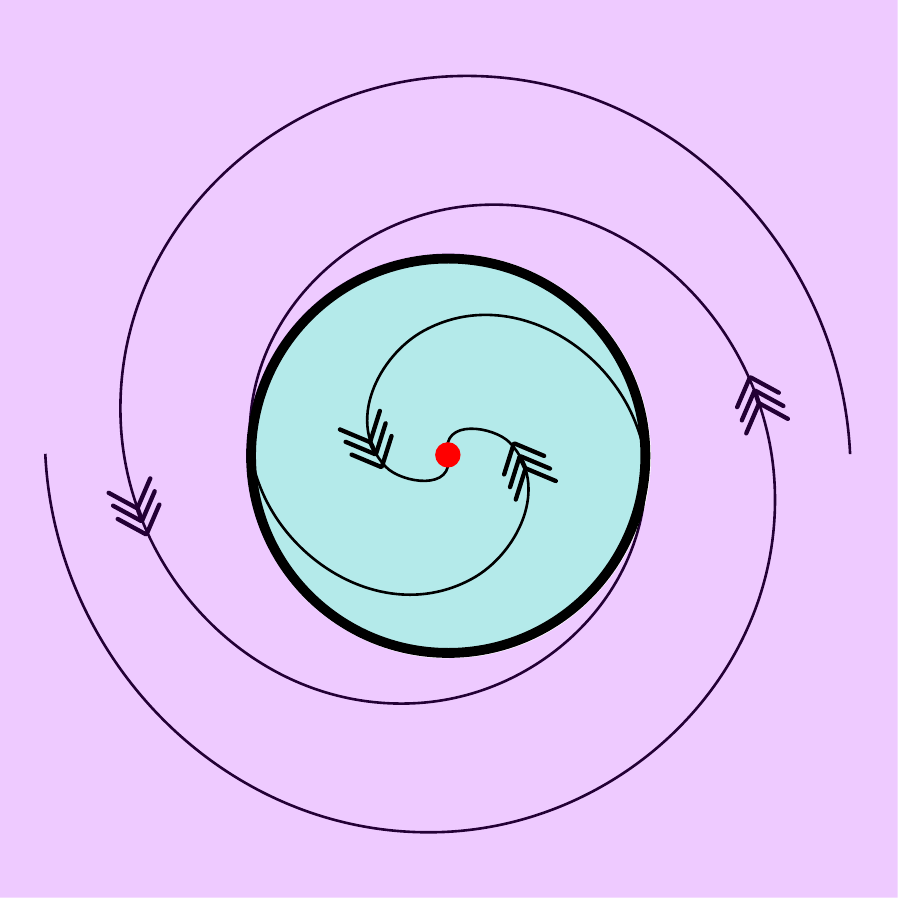}
		\end{overpic}
		\caption{Critical set, no equilibria on the equator}
		\label{fig:CM0}
	\end{subfigure}
	\hfil
	\begin{subfigure}{0.3\textwidth}
		\centering
		\begin{overpic}[width=0.9\textwidth]{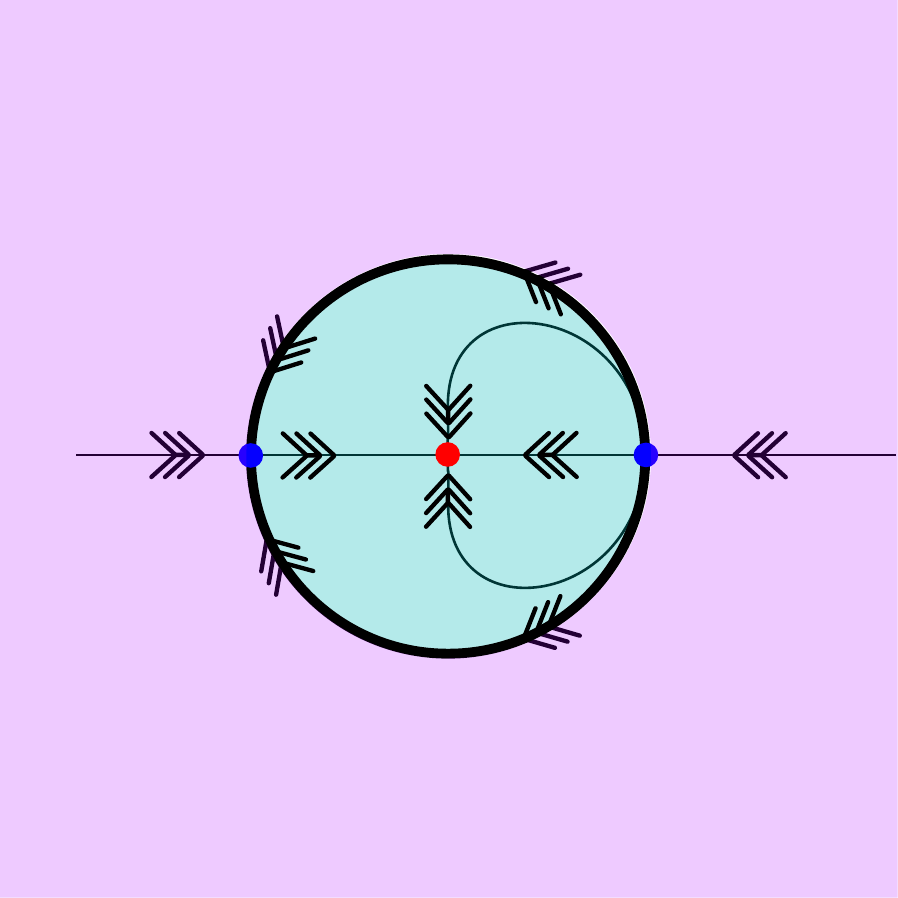}
		\end{overpic}
		\caption{Critical set, 2 equilibria on the equator}
		\label{fig:CM2}
	\end{subfigure}
	\hfil
	\begin{subfigure}{0.3\textwidth}
		\centering
		\begin{overpic}[width=0.9\textwidth]{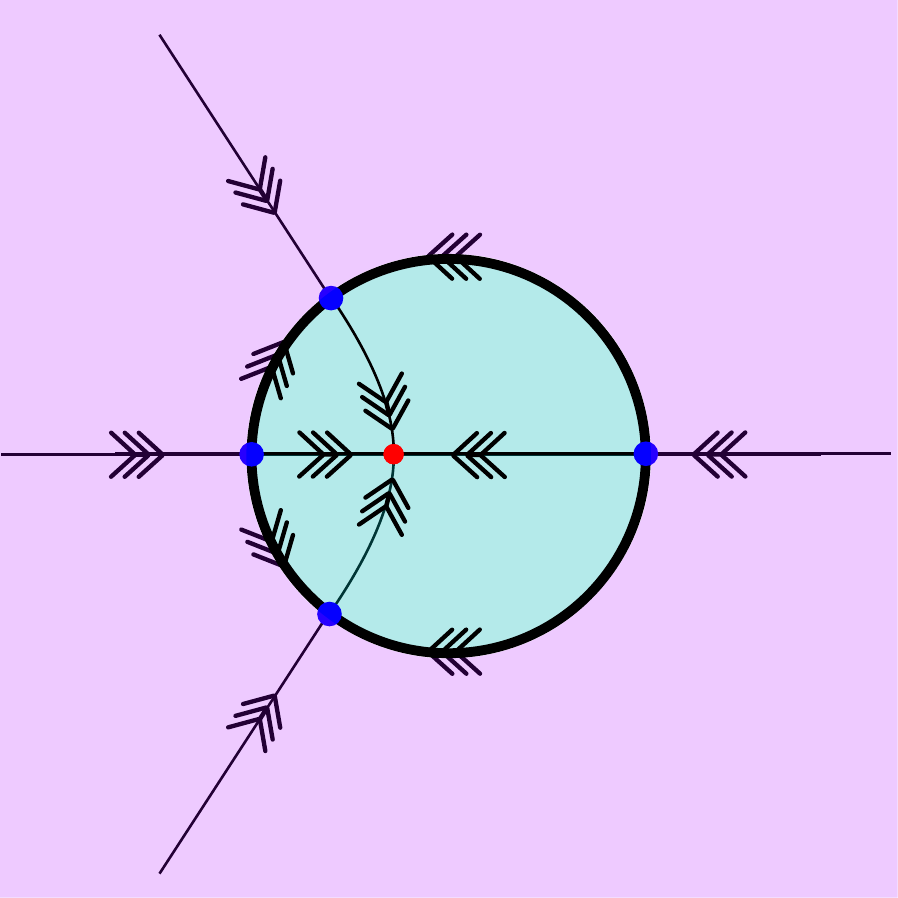}
		\end{overpic}
		\caption{Critical set, 4 equilibria on the equator}
		\label{fig:CM4}
	\end{subfigure}
	\hfil
	\begin{subfigure}{0.3\textwidth}
		\centering
		\begin{overpic}[width=0.9\textwidth]{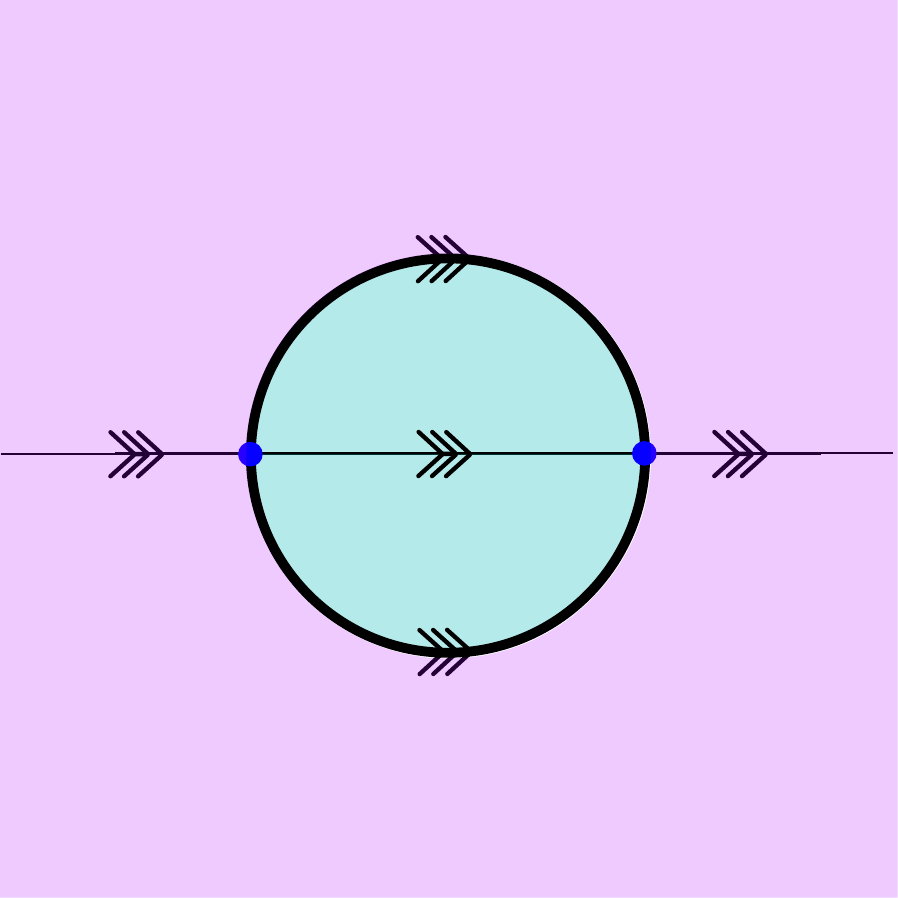}
		\end{overpic}
		\caption{No critical set, 2 equilibria on the equator}
		\label{fig:NCM2}
	\end{subfigure}
	\hfil
	\begin{subfigure}{0.3\textwidth}
		\centering
		\begin{overpic}[width=0.9\textwidth]{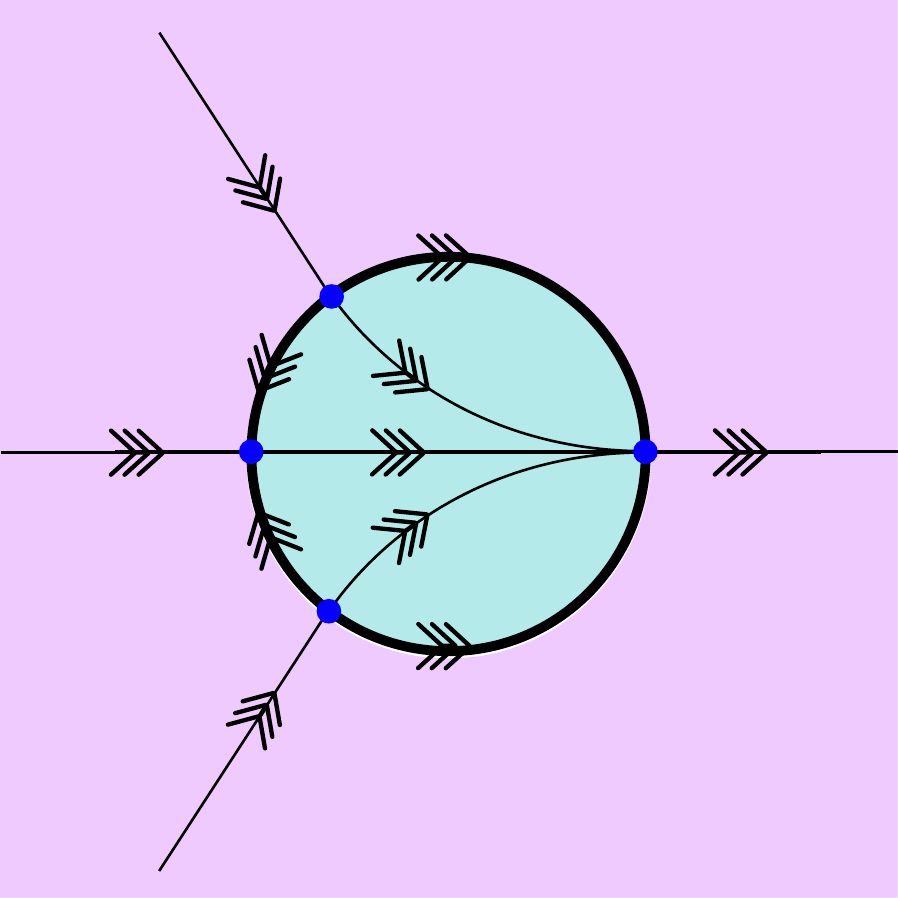}
		\end{overpic}
		\caption{No critical set, 4 equilibria on the equator}
		\label{fig:NCM4}
	\end{subfigure}
	\hfil
	\caption{	Sketches of the 5 possible cases of dynamics for Case I when $a,d<0$.
				The blown up sphere from \cref{fig:insert_sphere} is projected down onto a plane.
				Blue dots are equilibria of system along the equator whilst red dots are critical sets of the layer problem in the scaling chart $\kappa_2$. Sketches for $a,d>0$ can be obtained after a reversal of time.}\label{fig:sketches}
\end{figure}

We consider $a,d<0$ for fixed $\vec{z}$ and hence $\det(\mat{A})$ is necessarily positive (the analysis of $a,d>0$ is equivalent after a reversal of time). 
\begin{theorem}\label{thm:case1}\phantom{urgh}
	\begin{enumerate}[\rm{(\alph*)}]
	\item{	If $(\mat{A}(\vec{z})^{-1} \vec{f}(0,0,\vec{z}))^\intercal \mat{A}(\vec{z})^{-1} \vec{f}(0,0,\vec{z}) <1$,
			there exists a unique critical set in the scaling chart.
			Furthermore, all orbits in the layer problem of the scaling chart limit to the critical set. Generically, there can be 0, 2 or 4 equilibria along the equator for a fixed $\vec{z}$.}
	\item{	If $(\mat{A}(\vec{z})^{-1} \vec{f}(0,0,\vec{z}))^\intercal \mat{A}(\vec{z})^{-1} \vec{f}(0,0,\vec{z}) >1$,
			no critical set exists in the scaling chart and there can be either 2 or 4 equilibria along the equator generically.}
	\item{	There are 5 qualitatively different types of generic phase portraits for systems of the form \cref{eq:ode_normal_form} with \cref{eq:matrix_normal_form} for Case I.}
	\end{enumerate}
\end{theorem}
\begin{proof}
Let us consider \cref{eq:ode_normal_form} with \cref{eq:matrix_normal_form}.
From \cref{thm:lineps1}, if \cref{eq:fcondition} holds, then there exists a critical set in the layer problem in the scaling chart and that if $$\left(\mat{A}^{-1} \vec{f}(0,0,\vec{z})\right)^\intercal \mat{A}(\vec{z})^{-1} \vec{f}(0,0,\vec{z})>1$$ there is no critical set. 
\begin{enumerate}[\rm{(\alph*)}]
	\item{	Let us assume that such a critical set exists. The divergence of vector field of  the layer problem
			\begin{align}
				\mathrm{div}(\dot{x}_2,\dot{y}_2) &= \left(\frac{\partial}{\partial x_2}, \frac{\partial}{\partial x_2}\right) (\mat{A}(\vec{z})\vec{e}_\Psi(x_2,y_2;1)+\vec{f}(0,0,\vec{z}))\\
				&=\frac{(a(\vec{z})+d(\vec{z}))\Psi(\zeta^2)-(a(\vec{z})x_2^2+d(\vec{z})y_2^2)\Psi'(\zeta^2) +(a(\vec{z})y_2^2+d(\vec{z})x_2^2) }{(\zeta^2+\Psi(\zeta^2))}
			\end{align} 
			is negative everywhere for our class of regularisation functions (\cref{def:Psi}), so all orbits in the layer problem tend to it and there can be no limit cycles.
			Furthermore, it is straightforward to show that there can be 0, 2 or 4
			equilibria along the equator generically for a fixed $\vec{z}$ (see \cref{fig:CM0,fig:CM2,fig:CM4,fig:simCM0,fig:simCM2,fig:simCM4}).
			Each equilibria along the equator is radially attracting and the angular attractiveness can be determined using arguments from index theory.}
	\item{	If $\left(\mat{A}(\vec{z})^{-1} \vec{f}(0,0,\vec{z})\right)^\intercal \mat{A}(\vec{z})^{-1} \vec{f}(0,0,\vec{z})>1$,
			we have no critical set and from \cref{thm:geneps0}\cref{thm:geneps0c} there must be at least two equilibria along the equator. Again it is straightforward to show that, genericallly, it is possible to have either 2 or 4
			equilibria on the equator for a fixed $\vec{z}$, whose radial and angular attractiveness can be determined using arguments from index theory
			(see \cref{fig:NCM2,fig:NCM4,fig:simNCM2,fig:simNCM4}).}
	\item{	Thus, given the 3 possible types of dynamics when there exists a critical manifold, and the 2 possible types of dynamics when there
			is not critical set, there are 5 different possible types of dynamics for Case I.}
\end{enumerate}
\end{proof}
Sketches of the 5 possible phase portraits for Case I when $a,d<0$ are given in \cref{fig:sketches}. Numerical examples are given in \cref{fig:sim_5cases} in \cref{sec:numCaseI}.

Coulomb friction falls into our Case I and results in $b(\vec{z})\equiv0$. In that case, the hyperbola \cref{eq:conic11} passes through the origin and so there are at least 2 intersections with the unit circle: at least 2 equilibria along the equator. Hence \cref{fig:CM0} is not possible for Coulomb friction and \cref{fig:CM2,fig:CM4,fig:NCM2,fig:NCM4} correspond to the "four generic cases" in \cite[\S 4.2.4]{antali2019nonsmooth}.

\subsubsection{Case II}
\begin{figure}[htbp]
		\centering
		\begin{subfigure}{0.3\textwidth}
			\begin{overpic}[width=0.9\textwidth]{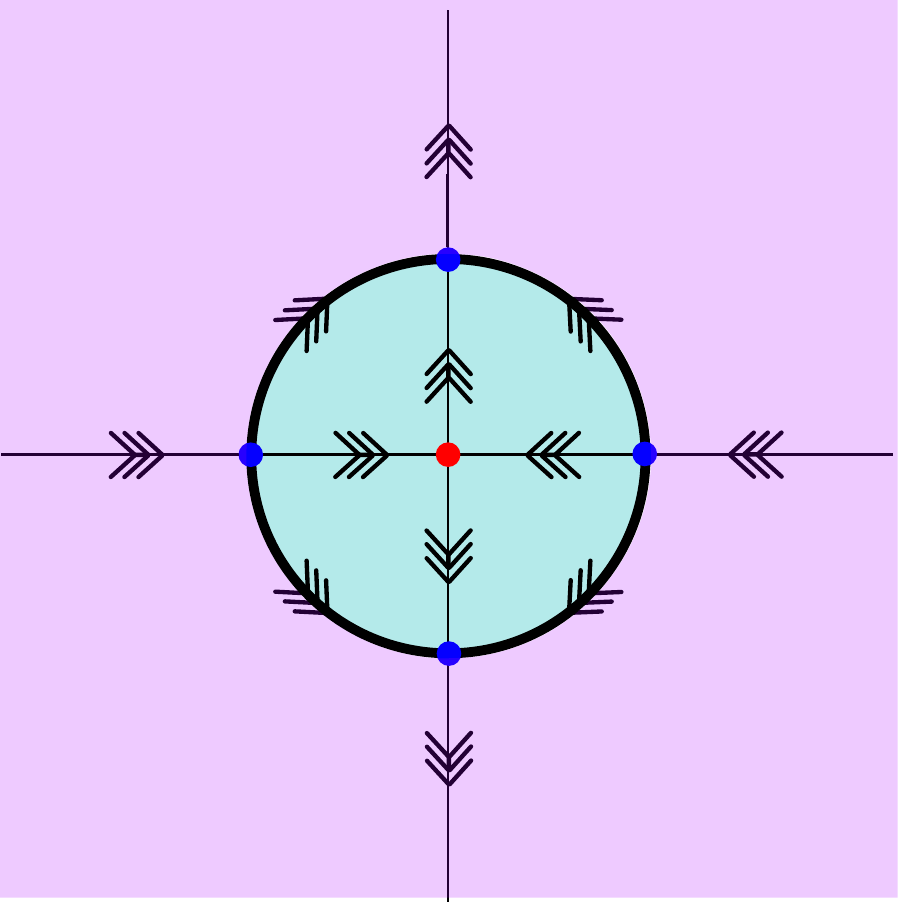}
			\end{overpic}
			\caption{}\label{fig:CM4_saddle}
		\end{subfigure}
		\hfil
		\begin{subfigure}{0.3\textwidth}
			\begin{overpic}[width=0.9\textwidth]{Figs/6cases/NCM2_2}
			\end{overpic}
			\caption{}\label{fig:NCM2_saddle}
		\end{subfigure}
		\hfil
		\begin{subfigure}{0.3\textwidth}
			\begin{overpic}[width=0.9\textwidth]{Figs/6cases/NCM4_2}
			\end{overpic}
			\caption{}\label{fig:NCM4_saddle}
		\end{subfigure}
		\caption{	Sketches of the 3 possible qualitative phase portraits of Case II for the blown up system of \cref{eq:ode_normal_form} when $\sign(a)\neq\sign(d)$ and $\det(\mat{A}) = ad+b^2<0$.
					\cref{fig:CM4_saddle} shows the only possible phase portrait when there is a critical set in the scaling chart.
					\cref{fig:NCM2_saddle,fig:NCM4_saddle} (equivalent to \cref{fig:NCM2,fig:NCM4}) show the two
					possible phase portraits when there is no critical set in the scaling chart.}\label{fig:3cases}
	\end{figure}

Here $a,d$ have opposite signs and $\det(\mat{A})<0$. In this case there are 3 qualitatively different generic phase portraits for fixed $\vec{z}$, as shown in \cref{fig:3cases}.
\begin{theorem}\label{thm:case2}\phantom{urgh}
	\begin{enumerate}[\rm{(\alph*)}]
		\item{	If $( \mat{A}(\vec{z})^{-1} \vec{f}(0,0,\vec{z}))^\intercal \mat{A}(\vec{z})^{-1} \vec{f}(0,0,\vec{z})<1$, then there exists exactly one unique critical set in the layer problem of the scaling chart, which is a saddle regardless of the regularisation function $\Psi$. Hence there are exactly four equilibria along the equator, each of which is a saddle with respect to the radial and angular flows. }
		\item{	If $( \mat{A}(\vec{z})^{-1} \vec{f}(0,0,\vec{z}))^\intercal \mat{A}(\vec{z})^{-1} \vec{f}(0,0,\vec{z})>1$, then there is no critical set in the layer problem
				of the scaling chart and there are generically either 2 or 4 equilibria along the equator as in Case I.}  
	\end{enumerate}
\end{theorem}

\begin{proof}
	Here it is straightforward to show that when $a$ and $d$ are of opposite signs there can be either one unique critical set in the scaling chart or none. 
	\begin{enumerate}[\rm{(\alph*)}]
		\item{	If $\left(\mat{A}(\vec{z})^{-1} \vec{f}(0,0,\vec{z})\right)^\intercal \mat{A}(\vec{z})^{-1} \vec{f}(0,0,\vec{z})<1$
				there is exactly one unique critical set in the scaling chart. Furthermore, if $\det{\mat{A}}=ad+b^2<0$ then,
				due to the positive definiteness of $D\vec{e}_\Psi$, the critical set is a saddle with respect the fast flow.
				Then from the uniqueness of the critical set we can infer that there will be exactly 4 equilibria along the equator:
				one for each of the intersections of the stable and unstable manifolds of the critical set with the equator.
				The radially attractiveness of these equilibria correspond to the stable/unstable manifolds of the equilibria in the scaling chart.
				This is the only possible type of dynamics when there is an equilibrium in the scaling chart. }
		\item{	Conversely, if $\left(\mat{A}(\vec{z})^{-1} \vec{f}(0,0,\vec{z})\right)^\intercal \mat{A}(\vec{z})^{-1} \vec{f}(0,0,\vec{z})>1$
				then there is no critical set in the scaling chart. It follows from index theory that, generically, there must be at least 2
				equilibria along the equator. It is straightforward to verify that there can be either 2 or 4 equilibria in this case.}
	\end{enumerate}
\end{proof}

\subsubsection{Case III}\label{sec:case3}
\begin{figure}[htbp]
	\centering
	\begin{subfigure}[b]{0.3\textwidth}
		\begin{overpic}[width=\textwidth]{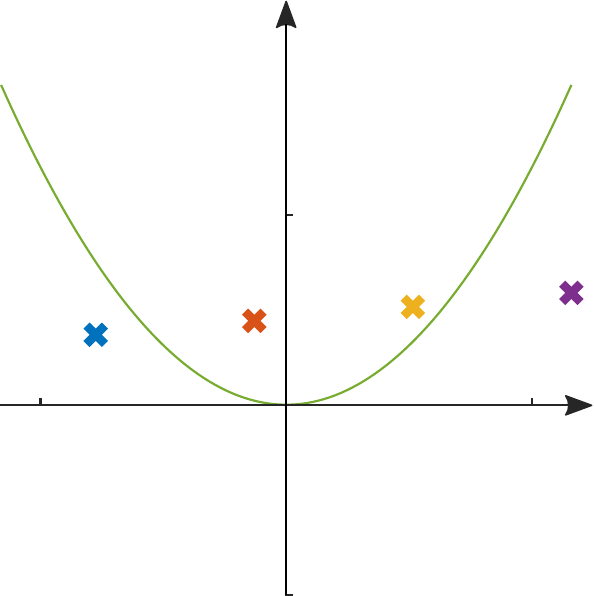}
			\put(92,36){$\mathrm{tr}(\mat{J})$}
			\put(52,94){$\mathrm{det}(\mat{J})$}
			\put(52,70){$\Delta(\mat{J})<0)$}
			\put(52,11){$\Delta(\mat{J})>0$}
			\put(0,24){$-50$}
			\put(83,24){$50$}
			\put(35,62){$500$}
		\end{overpic}
		\caption{Trace and determinant}\label{fig:traces_dets} 
	\end{subfigure}
	\hfil
	\begin{subfigure}[b]{0.3\textwidth}
		\begin{overpic}[width=\textwidth]{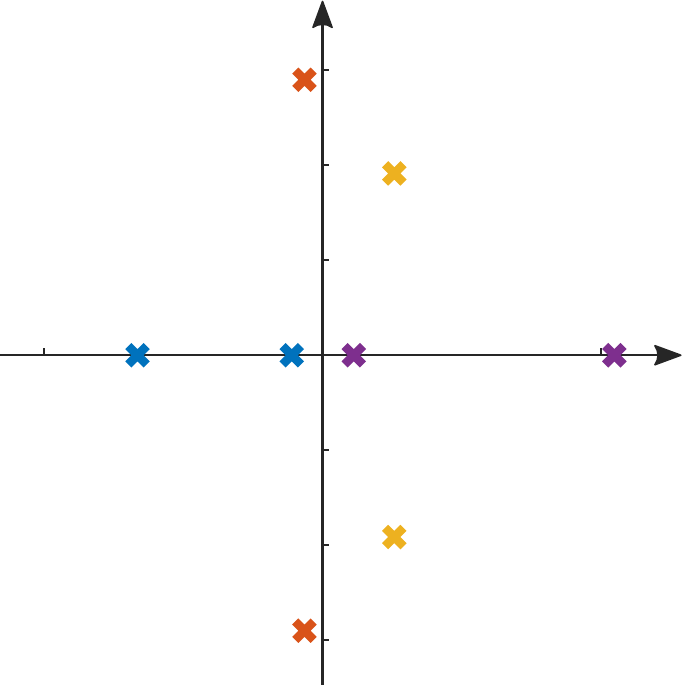}
			\put(82,53){$\mathrm{Re}(\lambda)$}
			\put(50,94){$\mathrm{Im}(\lambda)$}
			\put(32,73){{$\phantom{-}10$}}
			\put(32,18){{$-10$}}
		\end{overpic}
		\caption{Eigenvalues}\label{fig:eigs}
	\end{subfigure}
	\hfil
	\begin{subfigure}[b]{0.3\textwidth}
		\begin{overpic}[width=\textwidth]{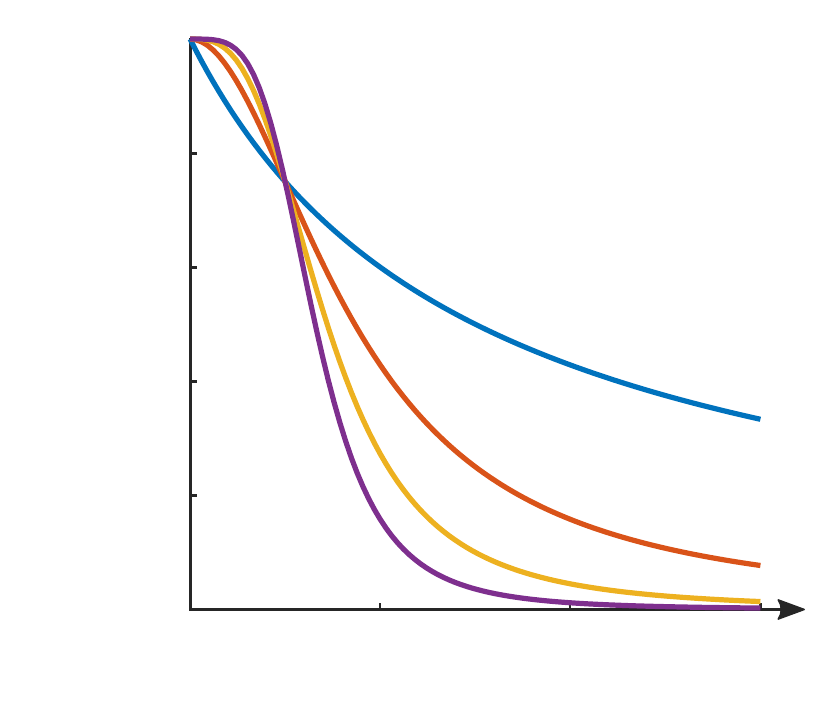}
			\put(16,10){0}
			\put(22,4){0}
			\put(68,4){1}
			\put(50,0){$r^2$}
			\put(17,77){1}
			\put(1,30){\rotatebox{90}{$\Psi(r^2;n)$}}
		\end{overpic}
		\caption{Regularisation function \cref{eq:bad_reg}}\label{fig:Psis}
	\end{subfigure}
	\caption{	Demonstration of the effect the choice of regularisation function can have on the dynamics for Case III.
				For $a(\vec{z})=172$, $b(\vec{z})=186$, $d(\vec{z})=-200$, $f_1(x,y,\vec{z})=-86$, $f_2(x,y,\vec{z})=-93$,
				choosing a different regularisation function \cref{eq:bad_reg} can cause the unique equilibrium to be a stable node, a stable focus, an unstable focus or  an unstable node (for $n=1${\color{matlab1}{$\blacksquare$}}, $n=2${\color{matlab2}{$\blacksquare$}}, $n=3${\color{matlab3}{$\blacksquare$}},
				and $n=4${\color{matlab4}{$\blacksquare$ }} respectively). (a) the trace $\tr(\mat{J})$ and determinant $\det(\mat{J})$ of the Jacobian of the linearisation around
				the unique equilibrium for each $n$, (b) the eigenvalues of the same matrix and (c) the regularisation function $\Psi(r^2;n)$ \cref{eq:bad_reg} for  $n =1\ldots4$. Note the different gradients at $r^2=\frac{1}{4}$ for each $n$. Also note that the sign of the discriminant of $\mat{J}$, $\Delta(\mat{J}):=\tr(\mat{J})^2-4\det(\mat{J})$, is also shown in (a).}\label{fig:bads}
\end{figure}

\begin{figure}[htbp]
	\centering
	\begin{subfigure}{0.4\textwidth}
		\centering
		\begin{overpic}[width=0.675\textwidth]{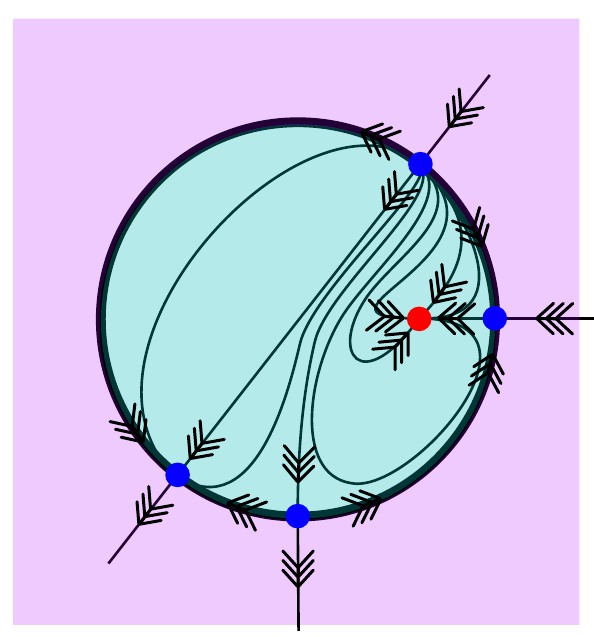}
		\end{overpic}
		\caption{$n=1$: $\mathrm{tr}(\mat{J})<0$, $\Delta(\mat{J})>0$}
	\end{subfigure}
	\hfil
	\begin{subfigure}{0.4\textwidth}
		\centering
		\begin{overpic}[width=0.675\textwidth]{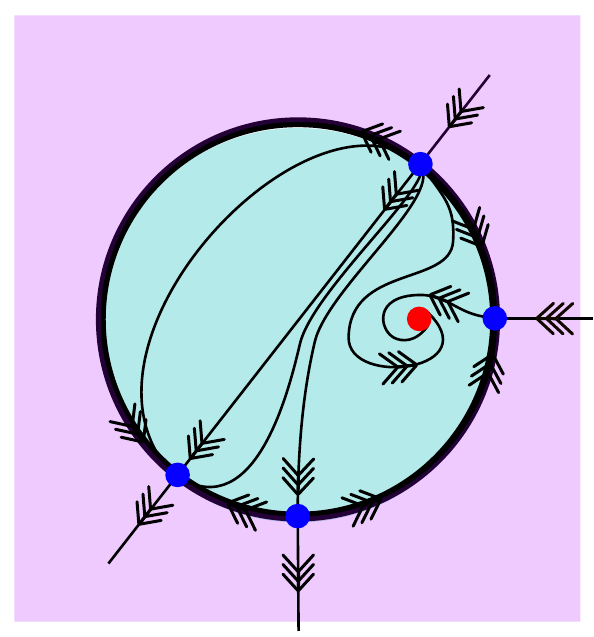}
		\end{overpic}
		\caption{$n=2$: $\mathrm{tr}(\mat{J})<0$, $\Delta(\mat{J})<0$}
	\end{subfigure}
	\begin{subfigure}{0.4\textwidth}
		\centering
		\begin{overpic}[width=0.675\textwidth]{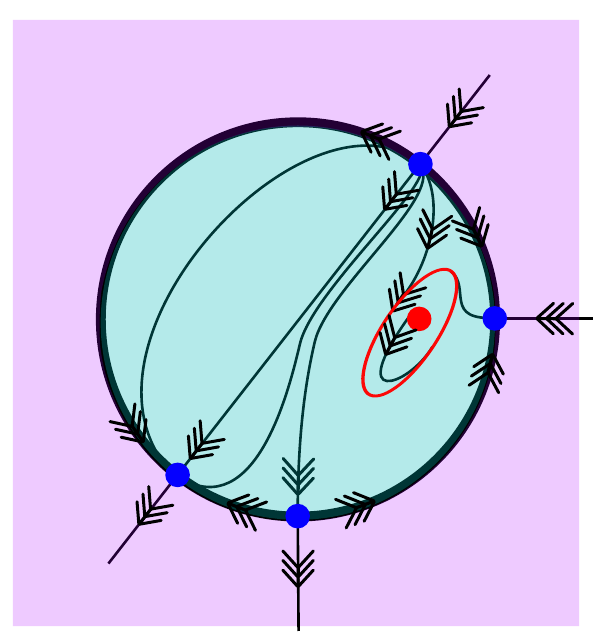}
		\end{overpic}
		\caption{$n=3$: $\mathrm{tr}(\mat{J})>0$, $\Delta(\mat{J})<0$}
	\end{subfigure}
	\hfil
	\begin{subfigure}{0.4\textwidth}
		\centering
		\begin{overpic}[width=0.675\textwidth]{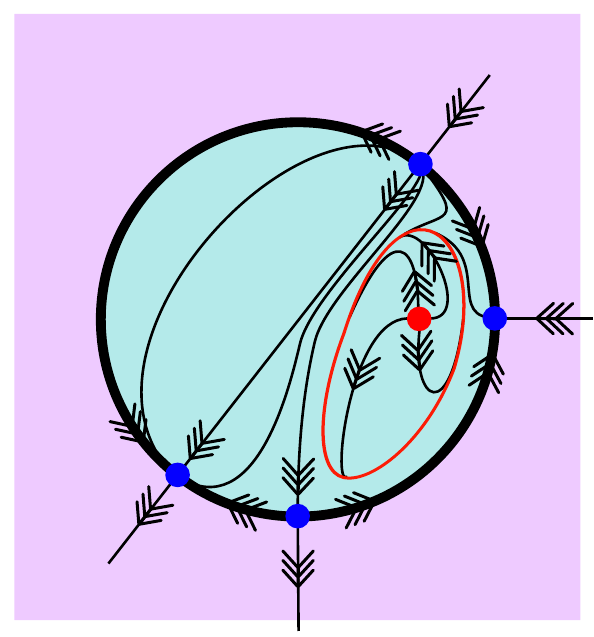}
		\end{overpic}
		\caption{$n=4$: $\mathrm{tr}(\mat{J})>0$, $\Delta(\mat{J})>0$}
	\end{subfigure}
	\caption{	Phase portraits for Case III of the ``\texorpdfstring{$\vec{e}$}{e}-linear'' system \cref{eq:ode_normal_form,eq:matrix_normal_form}, for parameters given in \cref{eq:bad_params}, with regularisation function \cref{eq:bad_reg}, for $n=1\ldots4$ (see also \cref{fig:traces_dets}). In (a) the critical manifold is a stable node, whereas in (b) it is a stable focus. In (c), $\mathrm{tr}(\mat{J})$ has changed sign and the critical manifold undergoes a supercritical Hopf bifurcation; it is now an unstable focus enclosed by a limit cycle (shown in red) in the layer problem. In (d), the limit cycle persists but the critical manifold has become an unstable node.}\label{fig:bad_sketches}
\end{figure}

For a system of the form \cref{eq:ode_normal_form} where $\sign{(a)}\neq\sign{(d)}$ and $\det{\mat{A}}=ad+b^2>0$, the dynamics are sensitive to the regularisation function $\Psi$ in the scaling chart $\kappa_2$. We demonstrate this in the following example.
Consider a system of the form \cref{eq:ode_normal_form} where 
\begin{equation}\label{eq:bad_params}
	a(\vec{z})=172,\quad b(\vec{z})=186, \quad d(\vec{z})=-200, \quad f_1(x,y,\vec{z})=-86, \quad f_2(x,y,\vec{z})=-93.
\end{equation}
Let us also consider a particular class of regularisation functions
\begin{equation}\label{eq:bad_reg}
	\Psi(r^2;n)=\frac{3}{3+(4r^2)^n},
\end{equation}
chosen such that $\Psi((\frac{1}{2})^2;n)=\frac{3}{4}$ is independent of $n$.
It is easy to show that the critical manifold in the scaling chart is a stable node for $n=1$,
a stable focus for $n=2$, an unstable focus for $n=3$ and an unstable spiral for $n=4$ (see \cref{fig:bads,fig:bad_sketches}).
In order to understand the stability of the critical manifold, we study the trace $\mathrm{tr}(\mat{J})$, determinant $\mathrm{det}(\mat{J})$ and  discriminant $\Delta(\mat{J}):=\mathrm{tr}(\mat{J})^2-4\mathrm{det}(\mat{J})$ of the Jacobian of the fast flow $\mat{J}$, given in \cref{eq:scjaco}, around the critical manifold, as shown in \cref{fig:bads}.

As evident from \cref{fig:bad_sketches}, limit cycles can exist in the scaling chart for Case III. The slow flow along these limit cycles is given by \cref{eq:linslowflow} (\cref{thm:lineps1c}). We should also note here that whether trajectories slide or cross depends upon their initial conditions.

We have demonstrated the sensitivity of the dynamics in the scaling chart $\kappa_2$ to the regularisation. Classification of all possible phase portraits for Case III with one up to one limit cycle in the scaling chart is possible using arguments from index theory.

\subsection{Bifurcations}\label{sec:bifurcations}
So far, we have studied the generic cases of \cref{eq:odes_smoothed} for fixed $\vec{z}$. Here we detail some possible bifurcations that are codimension-1 in parameter space, or which can be unfolded by the slow variable $\vec{z}$.
\begin{figure}[htbp]
	\centering
	\begin{subfigure}[t]{0.3\textwidth}
		\centering
		\begin{overpic}[width=0.9\textwidth]{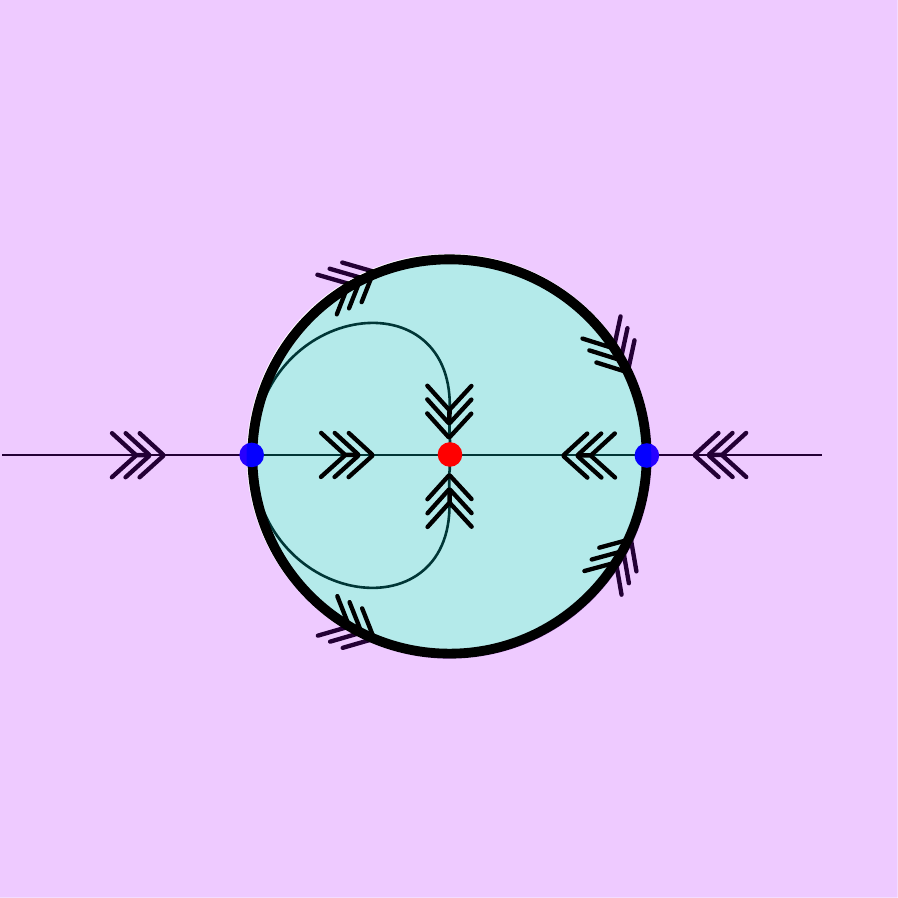}
		\end{overpic}
		\caption{\label{fig:transcrita}}
	\end{subfigure}\hfil
	\begin{subfigure}[t]{0.3\textwidth}
		\centering		
		\begin{overpic}[width=0.9\textwidth]{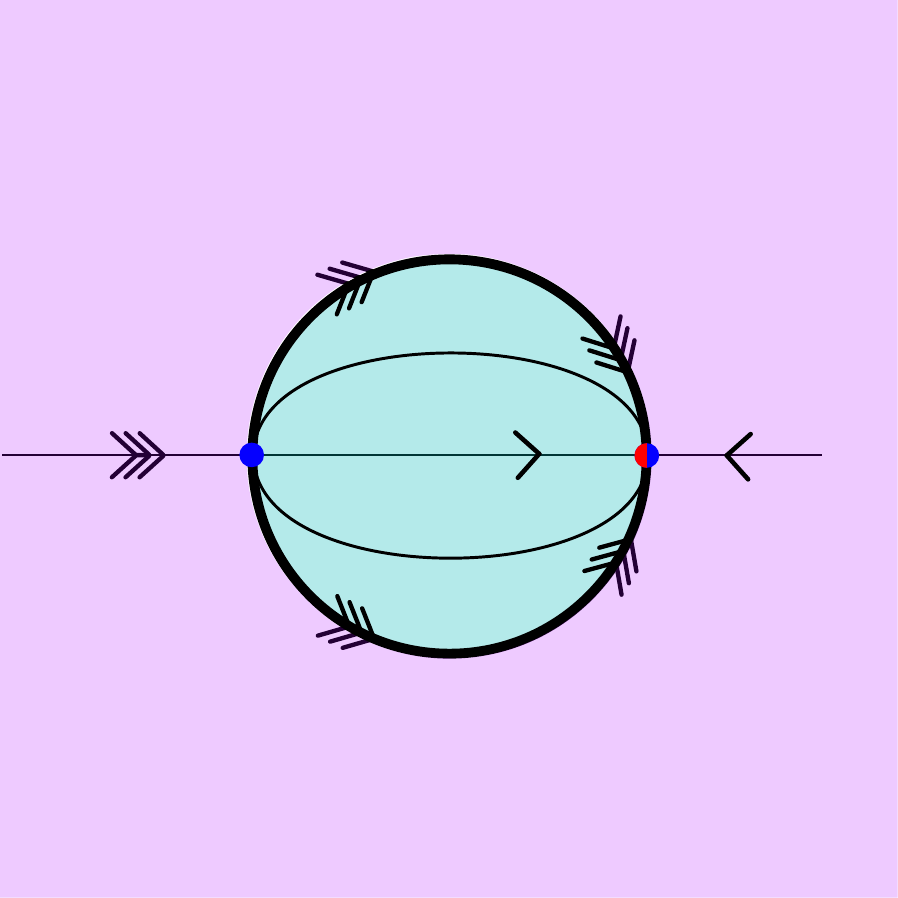}
		\end{overpic}
		\caption{\label{fig:transcritb}}
	\end{subfigure}\hfil
	\begin{subfigure}[t]{0.3\textwidth}
		\centering
		\begin{overpic}[width=0.9\textwidth]{Figs/6cases/NCM2_2}
		\end{overpic}
		\caption{\label{fig:transcritc}}
	\end{subfigure}
	\caption{	(a) there are 2 equilibria on the equator (blue dots), both radially attracting and of opposite angular attractiveness, and a critical set (in red) in the scaling chart $\kappa_2$.
				(b) there are 2 equilibria  on the equator, one of which is radially non-hyperbolic, with a critical set at infinity.
				(c) there are 2 equilibria along the equator, one of which is now radially repelling and no critical set in the scaling chart.
				}\label{fig:transcrit}
\end{figure}

In \cref{fig:transcrit} a critical set is created at infinity in the scaling chart $\kappa_2$, for the case when there are 2 equilibria on the equator (there could also be 4 equilibria on the equator). The radial attractiveness of one equilibrium switches at the bifurcation.  This bifurcation occurs when $$( \mat{A}(\vec{z})^{-1} \vec{f}(0,0,\vec{z}))^\intercal \mat{A}(\vec{z})^{-1} \vec{f}(0,0,\vec{z})=1.$$

\begin{figure}[htbp]
	\centering
	\hfil
	\begin{subfigure}[t]{0.3\textwidth}
		\begin{overpic}[width=0.9\textwidth]{Figs/6cases/CM0_2}
		\end{overpic}
		\caption{}
	\end{subfigure}
	\hfil
	\begin{subfigure}[t]{0.3\textwidth}
		\begin{overpic}[width=0.9\textwidth]{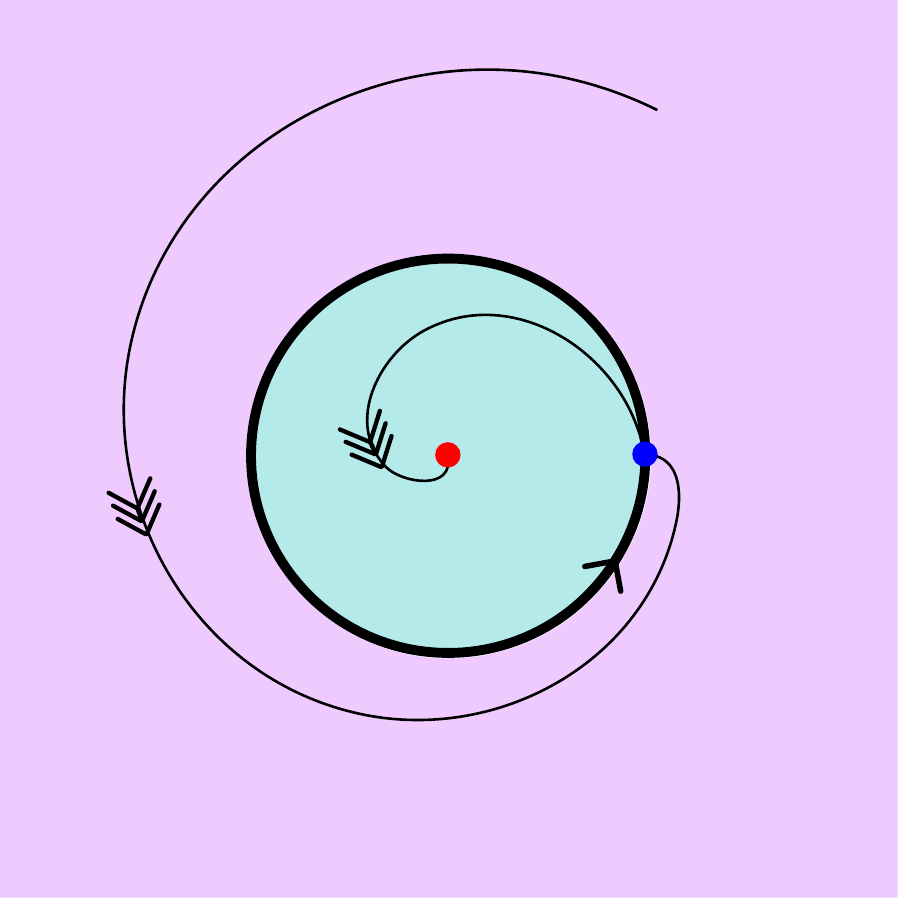}
		\end{overpic}
		\caption{}
	\end{subfigure}
	\hfil
	\begin{subfigure}[t]{0.3\textwidth}
		\begin{overpic}[width=0.9\textwidth]{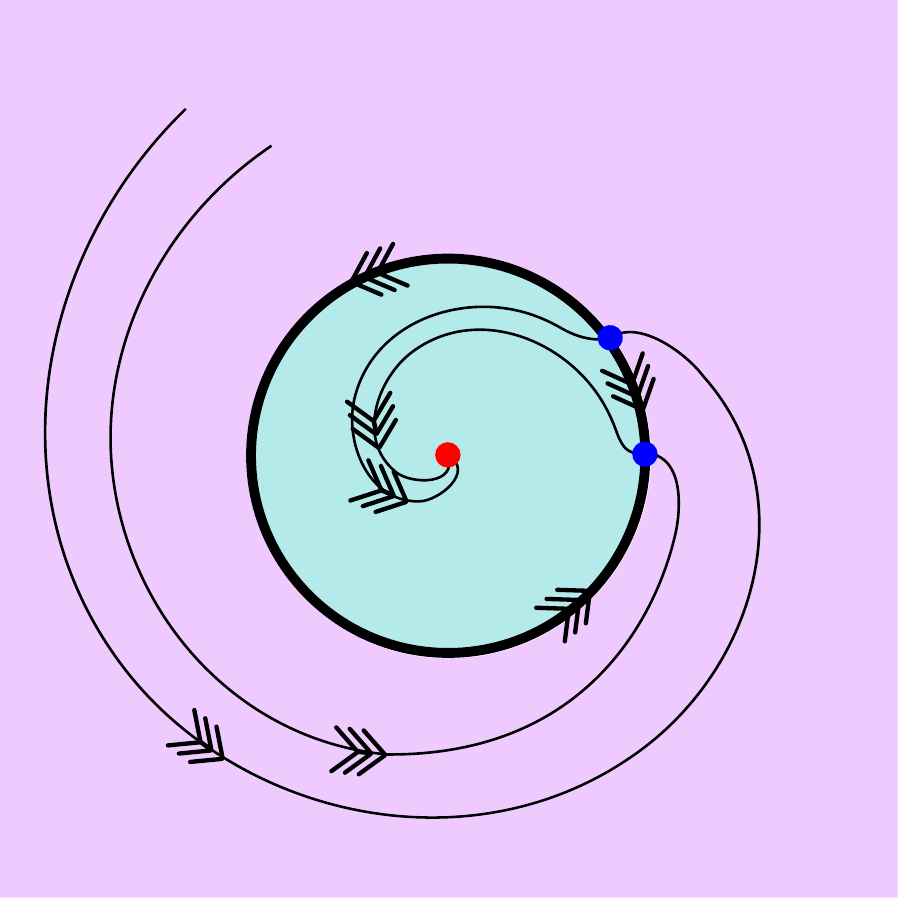}
		\end{overpic}
		\caption{}
	\end{subfigure}
	\hfil
	\caption{	Saddle-node bifurcation of equilibria on the equator. (a) no equilibria along the equator and a critical set in the scaling chart $\kappa_2$; (b) the saddle-node bifurcation: one angularly non-hyperbolic equilibrium appears on the equator with a heteroclinic connection to the critical set and a homoclinic connection; (c) two equilibria on the equator: one angularly attracting, the other repelling and both with heteroclinic connections to the critical set.}\label{fig:SNBIF}
\end{figure}

In \cref{fig:SNBIF}, we consider the saddle-node bifurcation between 0 and 2 equilibria on the equator where there is critical set in the scaling chart (we could also pass between 2 and 4 equilibria). Each pair of equilibria created on the equator has opposite angular attractiveness. This bifurcation occurs when when $\theta\mapsto \Theta(0,\theta,\vec{z})$ has a double root: $\Theta(0,\theta,\vec{z})=\pdiff{}{\theta}\Theta(0,\theta,\vec{z})=0$. Geometrically, this occurs when the hyperbola \cref{eq:conic11} is tangent to unit circle \cref{eq:conic12}). 

\section{Examples}\label{sec:examples}
We now discuss some examples, including some from \cite{Antali2017}.

\subsection{Ball at bottom of a pool\texorpdfstring{ \cite{Antali2017}}{}}
\begin{figure}[htbp]
\centering
\begin{overpic}[width=0.45\textwidth]{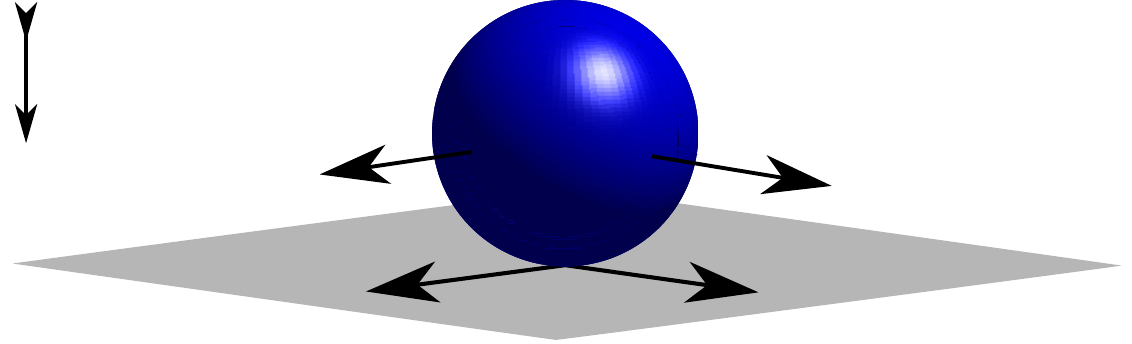}
\put(7,22){$\mathrm{g}$}
\put(30,6){$x$}
\put(68,6){$y$}
\put(27,20){$z_1$}
\put(69,19){$z_2$}
\end{overpic}
\caption{Ball at the bottom of a pool: $(x,y)$ is the relative velocity at the contact point, whilst $\vec{z}=(z_1,z_2)^\intercal$ is the velocity of the centre of mass.}
\label{fig:ballpool}
\end{figure}

\label[Example]{sec:ballpool}
We consider a rigid ball of mass $m$ slipping at the bottom of a pool with Coulomb coefficient of friction $\mu$ (acting on the ball at the point of contact) and viscous coefficient of friction $K>0$ (acting through the centre of mass), as in \cref{fig:ballpool}. 

The dynamics are governed by ordinary differential equations
\begin{equation}
\begin{pmatrix}
\dot{x}\\\dot{y}\\\dot{z}_1\\\dot{z}_2
\end{pmatrix}
=\vec{F}(x,y,z_1,z_2)=\begin{bmatrix}
-\frac{7}{2}\mu g \frac{x}{\sqrt[]{x^2+y^2}}-\frac{K}{m}z_1 \\[0.6em]
-\frac{7}{2}\mu g \frac{y}{\sqrt[]{x^2+y^2}}-\frac{K}{m}z_2 \\[0.6em]
-\mu g \frac{x}{\sqrt[]{x^2+y^2}}-\frac{K}{m}z_1 \\[0.6em]
-\mu g \frac{y}{\sqrt[]{x^2+y^2}}-\frac{K}{m}z_2 
\end{bmatrix}
\end{equation}
where $(x,y)$ is the relative velocity at the point of contact of the ball with the pool and $\vec{z}=(z_1,z_2)$ is the velocity of the centre of mass of the ball in the same axes (for details of the derivation, see \cite[Section 6.2]{Antali2017}). 
These equations can be written in our $\vec{e}$-linear normal form \cref{eq:ode_normal_form,eq:matrix_normal_form}, where
\begin{equation}
\mat{A}=-\frac{7}{2}\mu g\begin{pmatrix}
1&0\\
0&1
\end{pmatrix},\qquad\vec{f}(\vec{z})=-\frac{K}{m}\begin{pmatrix}
z_1\\
z_2
\end{pmatrix},\qquad
\mat{B}=-\mu g\begin{pmatrix}
1 &0 \\ 0 & 1
\end{pmatrix},\qquad
\vec{g}(\vec{z})=-\frac{K}{m}\begin{pmatrix}
z_1\\
z_2
\end{pmatrix}.
\end{equation}


In the entry chart $\kappa_1$, this is degenerate case \cref{deg:1}, since the diagonal elements of $\mat{A}$ are equal. Since the off-diagonal elements of $\mat{A}$ are zero, there are 2 unique equilibria on the equator, given by
\begin{equation}
\theta_1^*=\arctan{\left(\frac{f_2}{f_1}\right)}=\arctan {\left(\frac{z_2}{z_1}\right)},\quad \theta_2^*=\theta_1^*+\pi.
\end{equation}

These equilibria have eigenvalues given by
\begin{align}
&\lambda_\rho(\theta_1^*)=-\frac{7}{2}\mu g-\frac{K}{m}|\vec{z}|,&\lambda_\theta(\theta_1^*)=+\frac{K}{m}|\vec{z}|,\\
&\lambda_\rho(\theta_2^*)=-\frac{7}{2}\mu g+\frac{K}{m}|\vec{z}|,
&\lambda_\theta(\theta_2^*)=-\frac{K}{m}|\vec{z}|,
\end{align}
where $|\vec{z}|=\sqrt{z_1^2+z_2^2}$. Hence, as noted by \cite{Antali2017}, 
\begin{itemize}
\item the equilibrium $\theta_1^*$ is always angularly repelling and radially attracting,
\item the equilibrium $\theta_2^*$ is always angularly attracting, but radially attracting when $|\vec{z}| < z_{\mathrm{\mathrm{crit}}}$ and radially repelling when $|\vec{z}| > z_{\mathrm{crit}}$,
\end{itemize}
where $z_{\mathrm{crit}}=\frac{7\mu m g}{2K}>0$.

In the scaling chart $\kappa_2$ ($\bar{\varepsilon}=1$), the critical manifold \cref{eq:critical_set} is given by
\begin{equation}
C=\left\lbrace(x_2,y_2,\vec{z})\Bigg|\,c(\vec{z})=-\frac{2K}{7\mu m g}z_1,\,s(\vec{z})=-\frac{2K}{7\mu m g}z_2\right\rbrace.
\end{equation}
To exist, $C$ must lie within $D_1$, the unit disc centred on the origin (see \cref{def:MPsi}), and so
\begin{equation}\label{eq:constraint2}
|\vec{z}|<z_{\mathrm{crit}},
\end{equation}
equivalent to the constraint that $\theta_2^*$ is radially attracting.

When $C$ exists, the reduced dynamics along it are given by
\begin{equation}
\begin{split}
\dot{z}_1&=\frac{2K}{7m}z_1-\frac{K}{m}z_1=-\frac{5K}{7m}z_1\\
\dot{z}_2&=\frac{2K}{7m}z_2-\frac{K}{m}z_2=-\frac{5K}{7m}z_2
\end{split}
\end{equation}
and so the speed of the ball tends to zero in this case.
From \cref{thm:geneps1}, the stability of the critical manifold $C$ is given by the eigenvalues of
\begin{equation}
D_{(c,s)}\vec{U}((c,s)^\intercal,0,0,\vec{z})\,D_{(x_2,y_2)}\vec{e}_\Psi(x_2,y_2;1)=\begin{bmatrix}
-\frac{7}{2}\mu g & 0 \\ 0 &-\frac{7}{2}\mu g 
\end{bmatrix}D\vec{e}_\Psi.
\end{equation}
Since
\begin{equation}
\mathrm{tr}(D_{(c,s)}\vec{U}((c,s)^\intercal,0,0,\vec{z})\,D_{(x_2,y_2)}\vec{e}_\Psi(x_2,y_2;1))=-\frac{7}{2}\mu g \mathrm{tr}(D\vec{e}_\Psi)<0 
\end{equation}
and 
\begin{equation}
\mathrm{det}(D_{(c,s)}\vec{U}\,D\vec{e}_\Psi)=\frac{49}{4}\mu^2 g^2 \mathrm{det}(D\vec{e}_\Psi)>0,
\end{equation}
(using properties from \cref{rem:jaco_properties}) the critical manifold $C$, when it exists, is stable, independent of the regularisation $\Psi$, as expected from \cref{thm:case1}.

\begin{figure}[htbp]
\centering
\begin{overpic}[width=.5\textwidth]{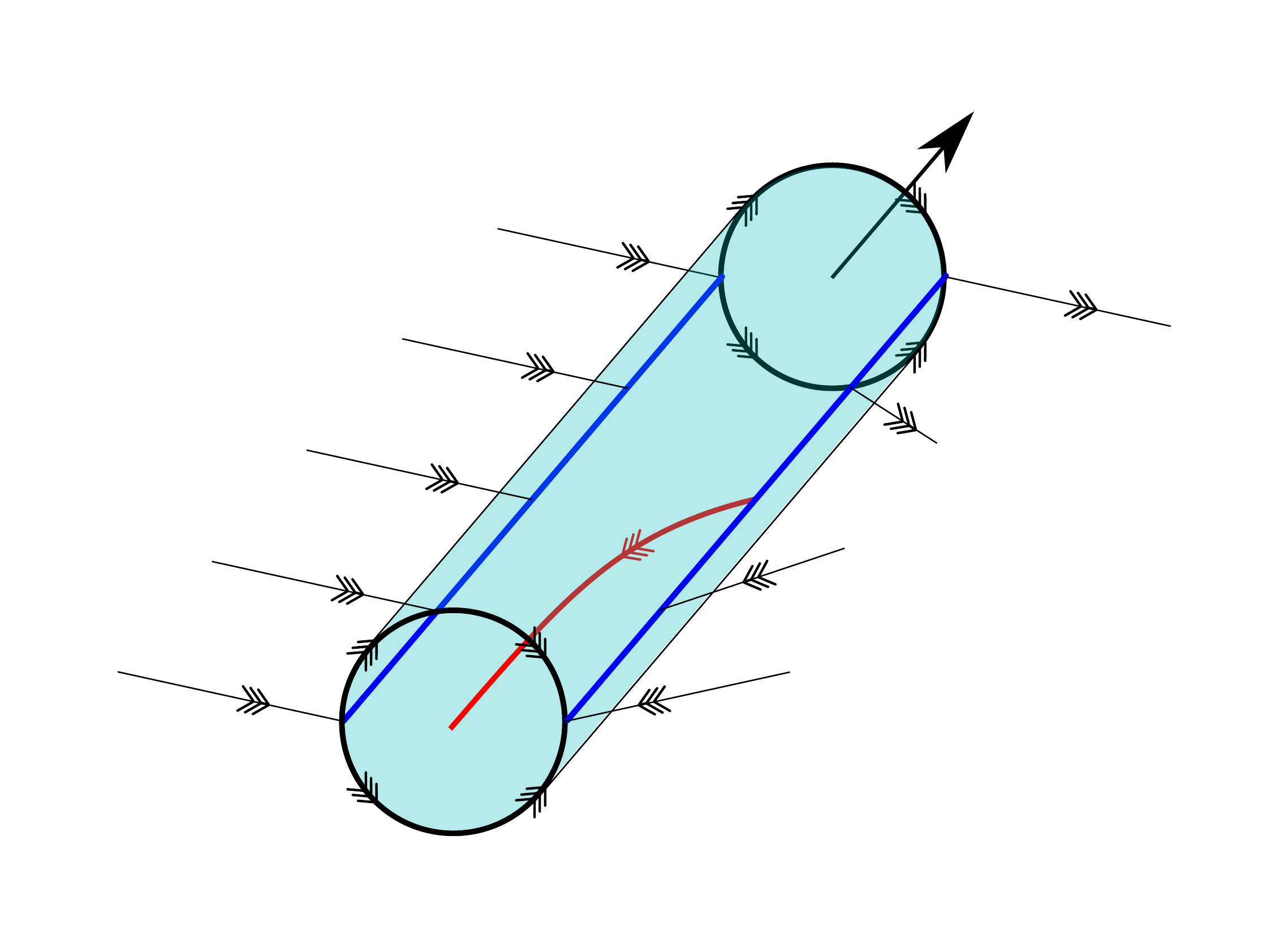}
\put(70,65){$\vec{z}$}
\put(67,34){$|\vec{z}|=z_{\mathrm{crit}}$}
\put(59,35){\line(2,0){7}}
\end{overpic}
\caption{Ball at the bottom of the pool. There is a tangency of the vector field to the discontinuity at $|\vec{z}|=z_{\mathrm{crit}} \equiv \frac{7\mu m \mathrm{g}}{2K}$. This could be considered to be the extension of a PWS invisible fold tangency to an isolated codimension-2 discontinuity \cite{filippov2013differential,kristiansen2019resolution}. Taking sections at $|\vec{z}|<z_{\mathrm{crit}}$, $|\vec{z}|=z_{\mathrm{crit}}$ and $|\vec{z}|>z_{\mathrm{crit}}$ gives \cref{fig:transcrita,fig:transcritb,fig:transcritc}, respectively. The cylindrical surface corresponds to the equator of the blowup for each $\vec{z}$ whilst the cylinder's volume corresponds to the scaling chart. The blue lines along the cylinder are therefore set of equilibria in the entry chart, whilst the red curve corresponds to the critical manifold in the scaling chart.}\label{fig:ballpool3d}
\end{figure}

This problem contains the bifurcation shown in \cref{fig:transcrit}, when the critical manifold $C$ appears on the equator at $\theta_2^*$ (scaled such that $\theta_2^*=0$). Physically, we can interpret the bifurcation as follows. In \cref{fig:transcrita}, all trajectories end up on the critical manifold $C$; the speed of the centre of mass $|\vec{z}| < z_{\mathrm{crit}}$, so the ball sticks and starts rolling, slowing due to viscous friction. In  \cref{fig:transcritb}, at the critical speed $|\vec{z}| = z_{\mathrm{crit}}$, the ball is on the point of slipping, as the equilibrium $\theta_2^*$ becomes radially non-hyperbolic at the bifurcation. Finally in \cref{fig:transcritc}, for $|\vec{z}| > z_{\mathrm{crit}}$ all trajectories are slipping except for the singular trajectory across the unit circle where the ball instantaneously sticks, only to slip again.
This bifurcation can occur dynamically as shown in \cref{fig:ballpool3d}.


\subsection{System with attracting and repelling directions}\label[Example]{sec:mateslidecross}
Consider the following system \cite[Example 4.5]{Antali2017}, which is nonlinear in $\vec{e}$.
\begin{equation}
\begin{split}
\dot{x}&=-\frac{x}{\sqrt{x^2+y^2}}+\frac{2x^2}{{x^2+y^2}}-\frac{1}{2}\\
\dot{y}&=-\frac{y}{\sqrt{x^2+y^2}}\\
\dot{z}&=-z.
\end{split}
\label{eq:ex4p5}
\end{equation}
Using our approach, we recover the results in  \cite{Antali2017} and quantify the effects of regularisation on the dynamics. 

Regularising \cref{eq:ex4p5}, we obtain
\begin{equation}\label{eq:ardir_example_reg}
\begin{split}
\dot{x}&=-\frac{x}{\sqrt{r^2+\varepsilon^2 \Psi(r^2/\varepsilon^2)}}+\frac{2x^2}{{r^2+\varepsilon^2 \Psi(r^2/\varepsilon^2)}}-\frac{1}{2}\\
\dot{y}&=-\frac{y}{\sqrt{r^2+\varepsilon^2 \Psi(r^2/\varepsilon^2)}}\\
\dot{z}&=-z
\end{split}
\end{equation}
where $r^2=x^2+y^2$.

In the entry chart $\kappa_1$, written in desingularised polar coordinates, \cref{eq:ardir_example_reg} becomes
\begin{equation}
\begin{split}
\frac{\mathrm{d}}{\mathrm{d}\mathcal{T}}\rho&=\rho\xi\left(2\cos^3\theta - \cos \theta -\frac{1}{2}\right),\\
\frac{\mathrm{d}}{\mathrm{d}\mathcal{T}}\theta&=\xi^{-1}\frac{\sin{\theta}}{2}\left(1-2\cos{\theta}\right)\left(1+2\cos{\theta}\right),\\
\frac{\mathrm{d}}{\mathrm{d}\mathcal{T}}z&=-\rho z,\\
\frac{\mathrm{d}}{\mathrm{d}\mathcal{T}}\varepsilon_1&=-\frac{\varepsilon_1}{\rho}\frac{\mathrm{d}}{\mathrm{d}\mathcal{T}}\rho.
\end{split}
\end{equation}

Equilibria on the equator occur at $Q_i=\{(\rho,\theta,z,\varepsilon_1)|\rho=0,\theta=\theta_i,\varepsilon_1=0\}$ where $\theta_i=(i-1)\frac{\pi}{3}$. Equilibria $Q_i,i\in\{1,3,5\}$ are angularly attracting ($\lambda_\theta<0$ for $\theta_i,i\in\{1,3,5\}$), whilst the others are angular repelling. Furthermore, $Q_i,i\in\{2,3,4,5,6\}$ are all radially attracting ($\lambda_\rho<0$ for $\theta_i,i\in\{2,3,4,5,6\}$), whilst $Q_1$ is radially repelling.

In the scaling chart $\kappa_2$ ($\bar{\varepsilon}=1$) critical sets are given implicitly by 

\begin{equation}\label{eq:criticalman}
C_{\pm}=\left\{(x_2,y_2)\left|y_2=0,x_2=x_2^*,x_2^*=\pm\sqrt{\left(1\pm\frac{2}{\sqrt{5}}\right)\Psi\left({x_2^*}^2\right)} \right.\right\},
\end{equation}
where $C_{+}$ is always a saddle, whilst $C_{-}$ is always a stable node.

\begin{figure}[hbtp]
\centering
\begin{overpic}[width=0.3\textwidth]{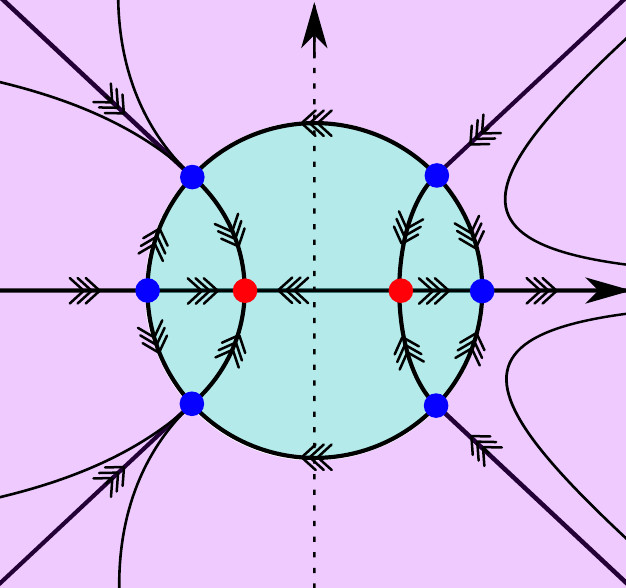}
\put(90,49){$\bar{x}$}
\put(53,88){$\bar{y}$}
\put(40,41){$C_-$}
\put(53,50){$C_+$}
\put(77.5,41.5){$Q_1$}
\put(59,63){$Q_2$}
\put(34,63){$Q_3$}
\put(26.5,40){$Q_4$}
\put(34,27){$Q_5$}
\put(58,28){$Q_6$}
\end{overpic}
\caption{Figure showing the blown-up discontinuity $x=y=\varepsilon=0$ of \cref{eq:ardir_example_reg}. Equilibria $Q_i$, $i \in\{1,3,5\}$, are angularly attracting whilst $Q_i$, $i \in\{2,4,6\}$, are angularly repelling. $Q_1$ is radially repelling whilst the other equilibria along the equator are radially attracting. Within the scaling chart $\kappa_2$, there are two critical manifolds $C_\pm$ \cref{eq:criticalman}, a saddle and a stable node respectively with respect to the fast flow.}
\label{fig:ex4p5}
\end{figure}
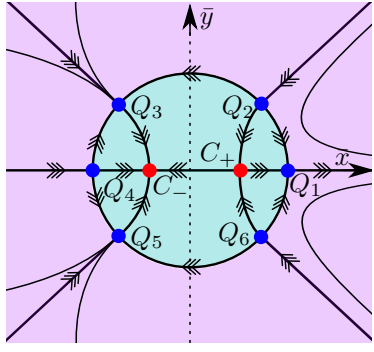

The combined dynamics projected into the ($\bar{x},\bar{y}$) plane is shown in \cref{fig:ex4p5}.
This example demonstrates the ambiguity of the terms {\it crossing} and {\it sliding} for problems with isolated codimension-2 discontinuities.
For initial conditions $(x_0,y_0)$ where $\left|\arctan\left({y_0}/{x_0}\right)\right|<\frac{\pi}{3}$ there is crossing, otherwise there is sliding.

%
Note that in this example \cref{eq:ex4p5}, $\dot{z}$ is independent of $\vec{e}$ and so the slow flow in the regularised problem is the same for the two critical sets $C_{\pm}$ \cref{eq:criticalman}. 
In the next example $\dot{\vec{z}}$ depends on $\vec{e}$. We will show that when there are multiple stable critical sets in the scaling chart with different corresponding slow flows, this can result in non-unique sliding in the nonsmooth limit. This sliding can then depend upon the direction along which trajectories approach the discontinuity set $\Sigma$.

\subsection{System with non-unique sliding vector field}\label[Example]{sec:nonunique}
Consider the following system, which is nonlinear in $\vec{e}$ and where $\dot{z}$ depends on $\vec{e}$.
\begin{equation}
\begin{split}
\dot{x}&=-\frac{x}{\sqrt{x^2+y^2}}\left(\frac{x}{\sqrt{x^2+y^2}}-\frac{1}{2}\chi(z)\right)\left(\frac{x}{\sqrt{x^2+y^2}}+\frac{1}{2}\chi(z)\right),\\
\dot{y}&=-\frac{1}{2}\frac{y}{\sqrt{x^2+y^2}},\\
\dot{z}&=-2\frac{x}{\sqrt{x^2+y^2}}+\frac{1}{4}
\end{split}
\label{eq:exmult}
\end{equation}
where $\chi(z)=(1+2z^2)/(1+z^2)$. Regularising as before, \cref{eq:exmult} becomes
\begin{equation}
\begin{split}
\dot{x}&=-\frac{x}{\sqrt{r^2+\varepsilon^2\Psi\left(\frac{r^2}{\varepsilon^2}\right)}}\left(\frac{x}{\sqrt{r^2+\varepsilon^2 \Psi\left(\frac{r^2}{\varepsilon^2}\right)}}-\frac{1}{2}\chi(z)\right)\left(\frac{x}{\sqrt{r^2+\varepsilon^2\Psi\left(\frac{r^2}{\varepsilon^2}\right)}}+\frac{1}{2}\chi(z)\right),\\
\dot{y}&=-\frac{1}{2}\frac{y}{\sqrt{r^2+\varepsilon\Psi\left(\frac{r^2}{\varepsilon^2}\right)}},\\
\dot{z}&=-2\frac{x}{\sqrt{r^2+\varepsilon\Psi\left(\frac{r^2}{\varepsilon^2}\right)}}+\frac{1}{4}. 
\end{split}\label{eq:pre_blowup_cubic_example}
\end{equation}

Following the procedure \cref{eq:polarising_desingularising} of using polar coordinates in the plane $\bar{\varepsilon}=0$ (the entry chart $\kappa_1$), we find that \cref{eq:pre_blowup_cubic_example} becomes
\begin{equation}
\begin{split}
\frac{\mathrm{d}}{\mathrm{d}\mathcal{T}}{\rho}&=-\rho\left(\cos^4\theta-\frac{\chi(z)^2+2}{4}\cos^2\theta+\frac{1}{2}\right),\\
\frac{\mathrm{d}}{\mathrm{d}\mathcal{T}}{\theta}&=\sin\theta\cos\theta\left(\cos^2\theta-\frac{\chi(z)^2+2}{4}\right),\\
\frac{\mathrm{d}}{\mathrm{d}\mathcal{T}}{z}&=-2\rho\cos\theta 	z+\frac{1}{4},\\
\frac{\mathrm{d}}{\mathrm{d}\mathcal{T}}\varepsilon_1&=-\frac{\varepsilon_1}{\rho}\frac{\mathrm{d}}{\mathrm{d}\mathcal{T}}\rho.
\end{split}
\label{eq:exmultdesing}
\end{equation}
There are 8 equilibria (limit directions) $Q_i=\{(\rho,\theta,z)|\rho=0,\theta=\theta_i(z), i=1\dots8\}$ of \cref{eq:exmultdesing} along the equator, where
\begin{equation}
\begin{split}
\theta_1=0,\quad \theta_{2}=-\theta_8=\arccos\left(\frac{\sqrt{\chi^2+2}}{2}\right),\quad \theta_{3}=-\theta_7=\frac{\pi}{2},\quad  \theta_{4}=-\theta_6=\pi-\theta_2,\quad\theta_5=\pi.
\end{split}
\end{equation} 
They are all radially attracting ($\lambda_\rho(\theta_i(z),z)<0,\, \forall i$). 
The equilibria $\theta=\theta_i, i\in\{1,3,5,7\}$ are angularly repelling ($\lambda_\theta(\theta_i(z),z)>0, i\in\{1,3,5,7\}$) whilst the equilibria $\theta=\theta_i, i\in\{2,4,6,8\}$ are angularly attracting ($\lambda_\theta(\theta_i(z),z)<0, i\in\{2,4,6,8\}$). 
 
 We now proceed to study the dynamics in the scaling chart $\kappa_2$ of \cref{eq:pre_blowup_cubic_example} after the blowup. 
 We find the slow-fast system
 \begin{equation}
\begin{split}
\varepsilon\dot{x}_2&=-\frac{x_2}{\sqrt{\zeta^2+\Psi\left({\zeta^2}\right)}}\left(\frac{x_2}{\sqrt{\zeta^2+\Psi\left({\zeta^2}\right)}}-\frac{1}{2}\chi(z)\right)\left(\frac{x_2}{\sqrt{\zeta^2+\Psi\left({\zeta^2}\right)}}+\frac{1}{2}\chi(z)\right),\\
\varepsilon\dot{y}_2&=-\frac{1}{2}\frac{y_2}{\sqrt{\zeta^2+\Psi\left({\zeta^2}\right)}},\\
\dot{z}&=-2\frac{x_2}{\sqrt{\zeta^2+\Psi\left({\zeta^2}\right)}} 	+\frac{1}{4},
\end{split}\label{scaling_chart_cubic_example}
\end{equation}
where $\zeta=x_2^2+y_2^2$ as before.
There are three critical sets \cref{eq:critical_set} given by
\begin{subequations}
\begin{align}
C_{1}&=\left\{(x_2,y_2,z)|y_2=0,x_2/\sqrt{x_2^2+\Psi\left(x_2^2\right)}=-\frac{1}{2}\chi(z)\right\},\\
C_{2}&=\left\{(x_2,y_2,z)|y_2=0,x_2/\sqrt{x_2^2+\Psi\left(x_2^2\right)}=0 \right\},\\
C_{3}&=\left\{(x_2,y_2,z)|y_2=0,x_2/\sqrt{x_2^2+\Psi\left(x_2^2\right)}=+\frac{1}{2}\chi(z) \right\}.
\end{align}
\end{subequations}
For a fixed $z$, $C_1$ and $C_3$ are stable nodes in the layer problem for all regularisation functions $\Psi$ in our class (\cref{def:Psi}), whilst $C_2$ is a saddle. 

The slow flow along along $C_1$ is $\dot{z}=\chi(z)+\frac{1}{4}$, and $\dot{z}=-\chi(z)+\frac{1}{4}$ along $C_3$. We can see the dynamics of this example in \cref{fig:eps0_cubic,fig:eps1cubic_combi,fig:epsboth_cubic,fig:epsboth_cubic_3D}. 

These critical sets correspond to different sliding vector fields in the nonsmooth limit. Hence we have constructed an example where the sliding vector field is dependent upon the direction of  approach to the codimension-2 discontinuity set $\Sigma$. 

\section{Discussion and conclusions}


In this paper, we have proposed an approach for the study of general dynamical systems with isolated codimension-2 discontinuity sets \cref{eq:codim2}, by using GSPT and blowup to study a regularised version of the system, which also helps with the understanding of the robustness of these discontinuous problems to smoothing perturbations. 

For our general system \cref{eq:vf}, we have demonstrated that the methods and terminology of Filippov \cite{filippov2013differential} can be ambiguous. In particular, sliding and crossing can depend upon the approach to the discontinuity set (as in \cref{fig:bad_sketches} and \cref{sec:mateslidecross}), and it is  possible to have more than one sliding vector field (see \cref{sec:nonunique}). 

We have also proposed a natural extension of standard Filippov systems to codimension-2 problems in \cref{sec:lin}: the $\vec{e}$-linear system. For such systems, we have classified the dynamics into three cases and found analogues to sliding, crossing and the sliding vector field.

Within our framework there is potential for sliding if there is a critical set \cref{eq:RPCM} in the scaling chart $\kappa_2$. In particular, the sliding vector field corresponds to the slow flow along the critical set.
In the $\vec{e}$-linear system, this flow will be unique, even if there are limit cycles \cref{prop:LC}.
Whether or not there is sliding can depend upon the direction of approach to the discontinuity set (see \cref{fig:bad_sketches}), and we have the possibility of multiple sliding vector fields in \cref{sec:nonunique}.



In this paper, we consider the class of regularisation functions $\Psi$ that lead to monotonic smoothed step functions (see \cref{fig:component}). Hence within our current framework, we study Coulomb friction. Other types of regularisation would be needed to study other friction laws, such as stiction \cite{bossolini2017canards,kristiansen2021stiction}. 
Using our regularisation approach, certain physical phenomena that are not covered here may be facilitated. For example, we can expect stick-slip oscillations to be given by heteroclinic connections between equilibria along the equator, where one of the connections occurs in the chart $\kappa_1$ and the other in the chart $\kappa_2$ (slip and stick respectively).

\appendix

\section{Proof of \texorpdfstring{\cref{lem:ePsi}}{Lemma \ref*{lem:ePsi}}}\label{sec:lemePsi}
\begin{proof}
\begin{enumerate}[(a)]
\item{From the definition of $\mat{R}(\phi)$ \cref{eq:rotation_matrix}, we have that
\begin{align}
\vec{e}_\Psi\left(\mat{R}(\phi)(x,y)^\intercal;\varepsilon\right)&=\vec{e}_\Psi(x\cos\phi-y\sin\phi,x\sin\phi+y\cos\phi;\varepsilon).\\
&=\frac{\left(\left(x\cos\phi-y\sin\phi\right),\left(x\sin\phi+y\cos\phi\right)\right)^\intercal}{\sqrt{x^2+y^2+\varepsilon^2 \Psi\left(\frac{x^2+y^2}{\varepsilon^2}\right)}},\\
&=\mat{R}(\phi)\frac{\left(x,y\right)^\intercal}{\sqrt{x^2+y^2+\varepsilon^2 \Psi\left(\frac{x^2+y^2}{\varepsilon^2}\right)}},\\
&\equiv\mat{R}(\phi)\vec{e}_\Psi\left(x,y;\varepsilon\right)
\end{align}
since $\left|\mat{R}(\phi)(x,y)^\intercal\right|^2=\left|(x,y)^\intercal\right|^2$}
\item{From the definition of $\vec{e}_\Psi$ \cref{eq:ePsi}, we have that
\begin{align}
\vec{e}_\Psi(k x,k y;k \varepsilon )&=\frac{(kx,ky)^\intercal}{\sqrt{k^2x^2+k^2y^2+k^2\varepsilon^2\Psi\left(\frac{k^2x^2+k^2y^2}{k^2\varepsilon^2}\right)}}\\
&=\frac{k(x,y)^\intercal}{\sqrt{k^2}\sqrt{x^2+y^2+\varepsilon^2\Psi\left(\frac{x^2+y^2}{\varepsilon^2}\right)}}.\\
\intertext{Therefore, for $k>0$,}
\vec{e}_\Psi(k x,k y;k \varepsilon )&=\frac{(x,y)^\intercal}{\sqrt{x^2+y^2+\varepsilon^2\Psi\left(\frac{x^2+y^2}{\varepsilon^2}\right)}},\\
&\equiv \vec{e}_\Psi(x,y;\varepsilon ).
\end{align}}
\item{Writing
\begin{align}
\vec{e}_\Psi(\tilde{\rho}\cos{\tilde{\theta}},\tilde{\rho}\cos{\tilde{\theta}};\varepsilon)&=\vec{e}_\Psi({\rho}\cos{{\theta}},{\rho}\cos{{\theta}};\varepsilon)\\
\intertext{in terms of polar coordinates,}
\left\langle\frac{\tilde{\rho}}{\sqrt{\tilde{\rho}^2+\varepsilon^2\Psi(\tilde{\rho}^2/\varepsilon^2)}},\tilde{\theta}\right\rangle&=\left\langle\frac{{\rho}}{\sqrt{{\rho}^2+\varepsilon^2\Psi({\rho}^2/\varepsilon^2)}},{\theta}\right\rangle.
\end{align}
Comparing components, we have $\theta\equiv\tilde{\theta}$, and from the monotonicity of $\Psi$ (\cref{def:Psi}), we have $\rho\equiv\tilde{\rho}$.}
\item{Straightforward calculation shows that
\begin{equation}
    \vec{e}_\Psi^\intercal\vec{e}_\Psi=\frac{x^2+y^2}{x^2+y^2+\varepsilon^2\Psi\left(\frac{x^2+y^2}{\varepsilon^2}\right)}<1.
\end{equation}
}
\end{enumerate}
\end{proof}

\section{Proof of \texorpdfstring{\cref{lem:normalform}}{Lemma \ref*{lem:normalform}}}\label{sec:lemnormalform}

\begin{lemma}\label{lem:normalform}
Consider \cref{eq:ode_normal_form}, which is linear in $\vec{e}$. Then there is transformation such that the system can be written as in \cref{eq:ode_normal_form}, with $\mat{A}$ in the form
\begin{equation}\label{eq:matrix_normal_form_app}
	\mat{A}(\vec{z})=
	\begin{pmatrix}
		a(\vec{z})& -b(\vec{z})\\
		b(\vec{z})& \phantom{-}d(\vec{z})
	\end{pmatrix}.
\end{equation}
\end{lemma}
\begin{proof}
	Let us consider the system
	\begin{equation*}
		\begin{pmatrix}
			{\dot{x}}\\{\dot{y}}\\\dot{\vec{z}}
		\end{pmatrix}=
		\begin{pmatrix}
			\mat{A}(\vec{z})\\ \mat{B}(\vec{z})
		\end{pmatrix}
		\vec{e}({x},{y})+
		\begin{pmatrix}
			\vec{f}(x,y,\vec{z})\\
			\vec{g}(x,y,\vec{z})
		\end{pmatrix}
\end{equation*}
where $\mat{A}$ is any real matrix 
$$\mat{{A}}({\vec{z}})=\begin{pmatrix}{a}(\vec{z})&{b}(\vec{z})\\{c}(\vec{z})&{d}(\vec{z})\end{pmatrix}.$$
We can write $\mat{{A}}(\vec{z})$ as
\begin{equation*}
	\mat{{A}}(\vec{z})=\mat{R}^{\intercal}\left(\phi(\vec{z})\right)\mat{\tilde{A}}(\vec{z})\mat{R}\left(\phi(\vec{z})\right)
\end{equation*}
where $\mat{R}$ is a rotation matrix through the angle 
$$\phi(\vec{z})=\frac{1}{2}\arctan\left(\frac{\tilde{b}(\vec{z})+\tilde{c}(\vec{z})}{\tilde{d}(\vec{z})-\tilde{a}(\vec{z})}\right),$$
and $\mat{\tilde{A}}$ is a matrix in the same form as \eqref{eq:matrix_normal_form_app}
$$\mat{\tilde{A}}:=\begin{pmatrix}
	\tilde{a}(\vec{z})&-\tilde{b}(\vec{z})\\\tilde{b}(\vec{z})&\tilde{d}(\vec{z})
\end{pmatrix}.$$
Let us now change to coordinates $(\tilde{x},\tilde{y},z)$ using the coordinate transformation given by 
\begin{equation}\label{eq:transform}
	(\tilde{x},\tilde{y})^\intercal=\mat{R}\left(\phi(\vec{z})\right)({x},{y})^\intercal.
\end{equation}
Differentiating \eqref{eq:transform}
\begin{align}
	\begin{pmatrix}
		\dot{\tilde{x}}\\\dot{\tilde{y}}
	\end{pmatrix}
	&=\mat{R}(\phi(\vec{z})) 
	\begin{pmatrix}
		\dot{x}\\\dot{y}
	\end{pmatrix} 
	+\pdiff{\phi}{\vec{z}}\dot{\vec{z}}\pdiff{}{\phi}\left(\mat{R}(\phi(\vec{z}))\right)
	\begin{pmatrix}
		{x}\\{y}
	\end{pmatrix},\nonumber\\
	\intertext{then substituting and using the equivariance of $\vec{e}$ from \cref{lem:ePsi}\cref{lem:ePsia}}
	&=\mat{R}\left(\phi(\vec{z})\right)\left({\mat{A}}(\vec{z}) \mat{R}^\intercal\left(\phi(\vec{z})\right)
	{\vec{e}}({{\tilde{x},\tilde{y}}})+{\vec{f}}(\mat{R}^\intercal\left(\phi(\vec{z})\right)
	(\tilde{x},\tilde{y})^\intercal,\vec{z}),\right)+\nonumber\\
	&\qquad\quad\pdiff{\phi}{\vec{z}} \dot{\vec{z}}\pdiff{}{\phi} \mat{R}\left(\phi(\vec{z})\right)
	\mat{R}^\intercal\left(\phi(\vec{z})\right)(\tilde{x},\tilde{y})^\intercal\nonumber\\
	&=\mat{\tilde{A}}(\vec{z}) \vec{e}(\tilde{x},\tilde{y}) +\mat{R}(\phi(\vec{z}))\vec{f}(\mat{R}^\intercal\left(\phi(\vec{z})\right)
	(\tilde{x},\tilde{y})^\intercal,\vec{z})+ \pdiff{\phi}{\vec{z}} \dot{\vec{z}} \mat{R}\left(\pi/2\right) (\tilde{x},\tilde{y})^\intercal \nonumber\\
	(\dot{\tilde{x}},\dot{\tilde{y}})^\intercal&:= \mat{\tilde{A}}(\vec{z}) \vec{e}(\tilde{x},\tilde{y})+\vec{\tilde{f}}(\tilde{x},\tilde{y},\vec{z}).
	\intertext{Similarly}
	\dot{\vec{z}}&=\mat{B}(\vec{z})\mat{R}^\intercal(\phi(\vec{z})) \vec{e}(\tilde{x},\tilde{y})+\vec{g}(\mat{R}^\intercal (\phi(\vec{z}))(\tilde{x},\tilde{y})^\intercal,\vec{z}),\nonumber\\
	\dot{\vec{z}}&:= \mat{\tilde{B}}(\vec{z}) \vec{e}(\tilde{x},\tilde{y})+\vec{\tilde{g}}(\tilde{x},\tilde{y},\vec{z})
\end{align}
\end{proof}

\section{Numerical examples for Case I}\label{sec:numCaseI}
We present numerical examples of the 5 possible phase portraits for Case I of the ``\texorpdfstring{$\vec{e}$}{e}-linear'' system in \cref{sec:cases}, when $a,d<0$, to be compared with \cref{fig:sketches}.
\begin{figure}[htbp]
	\centering
	\begin{subfigure}{0.3\textwidth}
		\begin{overpic}[width=\textwidth]{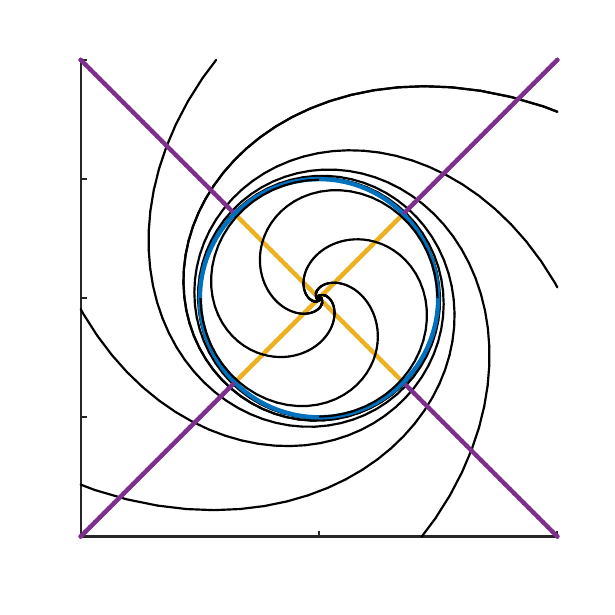}
			\put(89,6){2}
			\put(50,6){0}
			\put(10,6){-2}
			\put(8,87){2}
			\put(8,49){0}
			\put(6,11){-2}
		\end{overpic}
		\caption{$a=-1$, $b=1$, $d=-1$, $f_1=0$ \& $f_2=0$}\label{fig:simCM0}
	\end{subfigure}
	\hfil
	\begin{subfigure}{0.3\textwidth}
		\begin{overpic}[width=\textwidth]{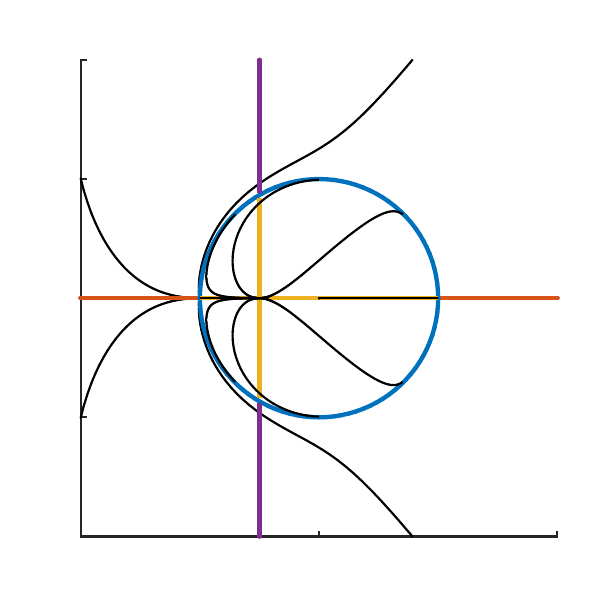}
			\put(89,6){2}
			\put(50,6){0}
			\put(10,6){-2}
			\put(8,87){2}
			\put(8,49){0}
			\put(6,11){-2}
		\end{overpic}
		\caption{$a=-1$, $b=0$, $d=-1.5$, $f_1=-0.5$ \& $f_2=0$}\label{fig:simCM2}
	\end{subfigure}
	\hfil
	\begin{subfigure}{0.3\textwidth}
		\begin{overpic}[width=\textwidth]{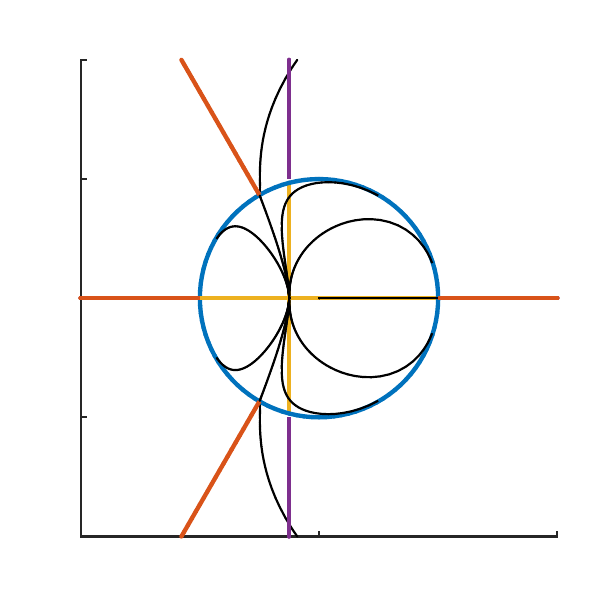}
			\put(89,6){2}
			\put(50,6){0}
			\put(10,6){-2}
			\put(8,87){2}
			\put(8,49){0}
			\put(6,11){-2}
		\end{overpic}
		\caption{$a=-2$, $b=0$, $d=-1$, $f_1=-0.5$ \& $f_2=0$}\label{fig:simCM4}
	\end{subfigure}
	\begin{subfigure}{0.3\textwidth}
		\begin{overpic}[width=\textwidth]{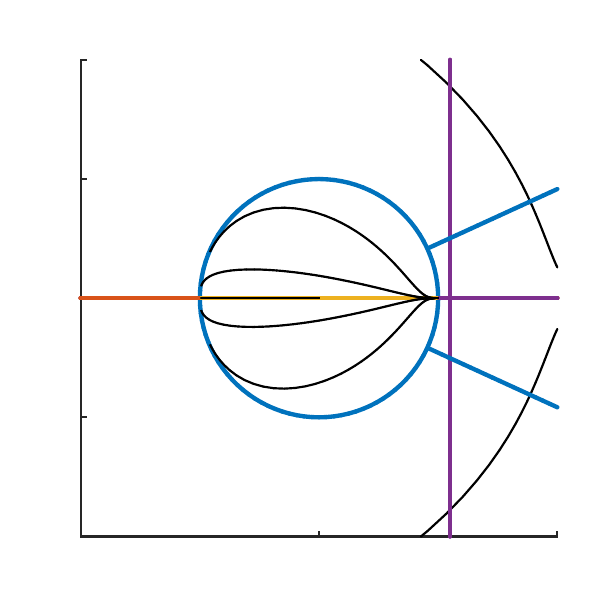}
			\put(89,6){2}
			\put(50,6){0}
			\put(10,6){-2}
			\put(8,87){2}
			\put(8,49){0}
			\put(6,11){-2}
		\end{overpic}
		\caption{$a=-1$, $b=0$, $d=-1$, $f_1=1.1$ \& $f_2=0$}\label{fig:simNCM2}
	\end{subfigure}
	\hfil
	\begin{subfigure}{0.3\textwidth}
		\begin{overpic}[width=\textwidth]{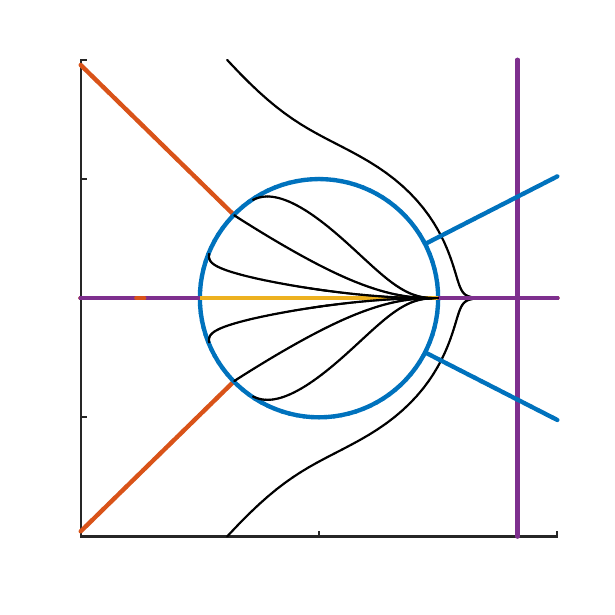}
			\put(89,6){2}
			\put(50,6){0}
			\put(10,6){-2}
			\put(8,87){2}
			\put(8,49){0}
			\put(6,11){-2}
		\end{overpic}
		\caption{$a=-0.3$, $b=0$, $d=-1$, $f_1=0.5$ \& $f_2=0$}\label{fig:simNCM4}
	\end{subfigure}
	\caption{	Numerical examples of  the 5 possible cases of dynamics of Case I when $a,d<0$.
				The blown up sphere from \cref{fig:insert_sphere} is projected down onto the $(\bar{x},\bar{y})$ plane, and the equator is the unit circle.
				Trajectories are plotted in black $\blacksquare$. In $\bar{\varepsilon}=0$, nullclines of $\dot{\rho}$ are shown in
				{\color{matlab1}{$\blacksquare$}} and nullclines of $\dot{\theta}$ are shown in {\color{matlab2}{$\blacksquare$}}.
				The lines $a\bar{x}-b\bar{y}=f_1$ and $b\bar{x}+d\bar{y}=f_2$ are also shown ({\color{matlab3}{$\blacksquare$}} when inside the the unit circle and {\color{matlab4}{$\blacksquare$}} outside).}
				\label{fig:sim_5cases}
\end{figure}
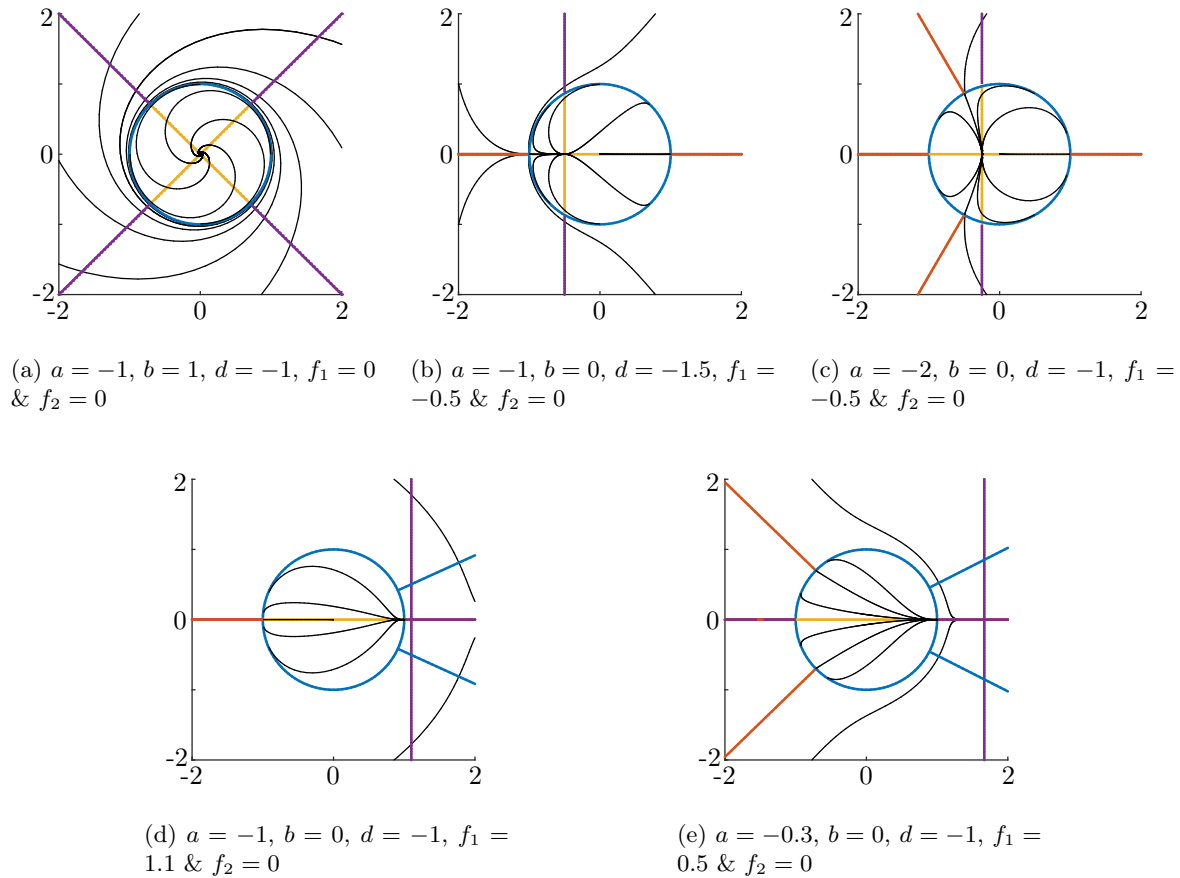

\section*{Acknowledgments}
We would like to thank Panagiotis Kaklamanos (University of Edinburgh), and M\'{a}t\'{e} Antali and G\'{a}bor St\'{e}p\'{a}n (Budapest University of Technology and Economics) for the time and useful discussions they kindly offered regarding their work. We would also like to thank Mario di Bernardo (University of Naples Federico II), Antonella Ferrara (University of Pavia), and Josep Olm (Technical University of Catalonia) for their illuminating correspondence on isolated codimension-2 discontinuities in sliding mode control.

\bibliographystyle{siamplain}

\bibliography{codim_bibliography}

\begin{thebibliography}{10}

\bibitem{Antali2017}
{\sc M.~Antali and G.~St\'ep\'an}, {\em Sliding and crossing dynamics in
  extended {F}ilippov systems}, SIAM J. Appl. Dyn. Sys, 17 (2018),
  pp.~823--858, \url{https://doi.org/10.1137/17M1110328}.

\bibitem{antali2019nonsmooth}
{\sc M.~Antali and G.~St\'ep\'an}, {\em Nonsmooth analysis of three-dimensional
  slipping and rolling in the presence of dry friction}, Nonlinear Dynamics, 97
  (2019), pp.~1799--1817, \url{https://doi.org/10.1007/s11071-019-04913-x}.

\bibitem{bernardo2008piecewise}
{\sc M.~Bernardo, C.~Budd, A.~R. Champneys, and P.~Kowalczyk}, {\em
  Piecewise-smooth dynamical systems: theory and applications}, vol.~163,
  Springer Science \& Business Media, 2008,
  \url{https://doi.org/10.1007/978-1-84628-708-4}.

\bibitem{bossolini2017canards}
{\sc E.~Bossolini, M.~Br{\o}ns, and K.~U. Kristiansen}, {\em Canards in
  stiction: on solutions of a friction oscillator by regularization}, SIAM J.
  Appl. Dyn. Sys, 16 (2017), pp.~2233--2258,
  \url{https://doi.org/10.1137/17M1120774}.

\bibitem{carmona2008existence}
{\sc V.~Carmona, F.~Fern{\'a}ndez-S{\'a}nchez, and A.~E. Teruel}, {\em
  Existence of a reversible {T}-point heteroclinic cycle in a piecewise linear
  version of the {M}ichelson system}, SIAM J. Appl. Dyn. Sys., 7 (2008),
  pp.~1032--1048, \url{https://doi.org/10.1137/070709542}.

\bibitem{dieci2017moments}
{\sc L.~Dieci and F.~Difonzo}, {\em The moments sliding vector field on the
  intersection of two manifolds}, J. Dyn. Diff. Eqs, 29 (2017), pp.~169--201,
  \url{https://doi.org/10.1007/s10884-015-9439-9}.

\bibitem{dieci2013filippov}
{\sc L.~Dieci, C.~Elia, and L.~Lopez}, {\em A {F}ilippov sliding vector field
  on an attracting co-dimension 2 discontinuity surface, and a limited
  loss-of-attractivity analysis}, J. Diff. Eqs, 254 (2013), pp.~1800--1832,
  \url{https://doi.org/10.1016/j.jde.2012.11.007}.

\bibitem{dieci2011sliding}
{\sc L.~Dieci and L.~Lopez}, {\em Sliding motion on discontinuity surfaces of
  high co-dimension. {A} construction for selecting a {F}ilippov vector field},
  Numerische Mathematik, 117 (2011), pp.~779--811,
  \url{https://doi.org/10.1007/s00211-011-0365-4}.

\bibitem{dumortier1996canard}
{\sc F.~Dumortier, R.~Roussarie, and R.~H. Roussarie}, {\em Canard cycles and
  center manifolds}, vol.~577, American Mathematical Soc., 1996,
  \url{http://doi.org/10.1090/memo/0577}.

\bibitem{fenichel1979geometric}
{\sc N.~Fenichel}, {\em Geometric singular perturbation theory for ordinary
  differential equations}, J. Diff. Eqs, 31 (1979), pp.~53--98,
  \url{https://doi.org/10.1016/0022-0396(79)90152-9}.

\bibitem{filippov2013differential}
{\sc A.~F. Filippov}, {\em Differential equations with discontinuous righthand
  sides}, vol.~18, Springer, 1988,
  \url{https://doi.org/10.1007/978-94-015-7793-9}.

\bibitem{Jeffrey2014}
{\sc M.~R. Jeffrey}, {\em Dynamics at a switching intersection: Hierarchy,
  isonomy, and multiple sliding}, SIAM J. Appl. Dyn. Sys, 13 (2014),
  pp.~1082--1105, \url{https://doi.org/10.1137/13093368X}.

\bibitem{jeffrey2018hidden}
{\sc M.~R. Jeffrey}, {\em Hidden dynamics: the mathematics of switches,
  decisions and other discontinuous behaviour}, Springer, 2018,
  \url{https://doi.org/10.1007/978-3-030-02107-8}.

\bibitem{kaklamanos2019regularization}
{\sc P.~Kaklamanos and K.~U. Kristiansen}, {\em Regularization and geometry of
  piecewise smooth systems with intersecting discontinuity sets}, SIAM J. Appl.
  Dyn. Sys, 18 (2019), pp.~1225--1264,
  \url{https://doi.org/10.1137/18M1214470}.

\bibitem{kristiansen2020regularized}
{\sc K.~U. Kristiansen}, {\em The regularized visible fold revisited}, J.
  Nonlin. Sci., 30 (2020), pp.~2463--2511,
  \url{https://doi.org/10.1007/s00332-020-09627-8}.

\bibitem{kristiansen2021stiction}
{\sc K.~U. Kristiansen}, {\em A stiction oscillator under slowly varying
  forcing: Uncovering small scale phenomena using blowup}, 2021,
  \url{https://arxiv.org/abs/2102.10658}.

\bibitem{kristiansen2015use}
{\sc K.~U. Kristiansen and S.~J. Hogan}, {\em On the use of blowup to study
  regularizations of singularities of piecewise smooth dynamical systems in
  $\mathbb{R}^{3}$}, SIAM J. Appl. Dyn. Sys, 14 (2015), pp.~382--422,
  \url{https://doi.org/10.1137/140980995}.

\bibitem{kristiansen2019resolution}
{\sc K.~U. Kristiansen and S.~J. Hogan}, {\em On the interpretation of the
  piecewise smooth visible--invisible two-fold singularity in {$\mathbb{R}^3$}
  using regularization and blowup}, J. Nonlin. Sci., 29 (2019), pp.~723--787,
  \url{https://doi.org/10.1007/s00332-018-9502-x}.

\bibitem{krupa2001extending}
{\sc M.~Krupa and P.~Szmolyan}, {\em Extending geometric singular perturbation
  theory to nonhyperbolic points---fold and canard points in two dimensions},
  SIAM J. Math. Anal., 33 (2001), pp.~286--314,
  \url{https://doi.org/10.1137/S0036141099360919}.

\bibitem{levant2017quasi}
{\sc A.~Levant and B.~Shustin}, {\em Quasi-continuous mimo sliding-mode
  control}, IEEE Transactions on Automatic Control, 63 (2017), pp.~3068--3074,
  \url{https://doi.org/10.1109/TAC.2017.2778251}.

\bibitem{michelson1986steady}
{\sc D.~Michelson}, {\em Steady solutions of the {K}uramoto-{S}ivashinsky
  equation}, Physica D, 19 (1986), pp.~89--111,
  \url{https://doi.org/10.1016/0167-2789(86)90055-2}.

\bibitem{perko}
{\sc L.~Perko}, {\em Differential equations and dynamical systems}, vol.~7,
  Springer Science \& Business Media, 2013.

\bibitem{teixeira2012regularization}
{\sc M.~A. Teixeira and P.~R. da~Silva}, {\em Regularization and singular
  perturbation techniques for non-smooth systems}, Physica D, 241 (2012),
  pp.~1948--1955, \url{https://doi.org/10.1016/j.physd.2011.06.022}.

\bibitem{utkin1999sliding}
{\sc V.~Utkin, J.~Guldner, and M.~Shijun}, {\em Sliding mode control in
  electro-mechanical systems}, vol.~34, CRC press, 1999,
  \url{https://doi.org/10.1201/9781420065619}.

\end{thebibliography}

\end{document}